\documentclass[a4paper,11pt]{article}

\usepackage{amsthm}
\usepackage{hyperref}
\usepackage{amssymb}
\usepackage{latexsym}
\usepackage{amsmath}
\usepackage{amsfonts}
\usepackage{mathabx}

\usepackage{mathrsfs} 

\usepackage[dvips]{graphicx}
\usepackage[utf8]{inputenc}
\usepackage{bbm}

\usepackage[english]{babel}
\usepackage{url}
\usepackage{enumerate}
\usepackage{hyperref}
\usepackage{color}
\usepackage{stmaryrd}

\usepackage{fancyhdr}
\newtheorem{Theorem}{Theorem}[part]
\newtheorem{Definition}{Definition}[part]
\newtheorem{Proposition}{Proposition}[part]

\newtheorem{Corollary}{Corollary}[part]
\newtheorem{Remark}{Remark}[part]

\makeatletter \@addtoreset{equation}{section}

\@addtoreset{Definition}{section}

\@addtoreset{Theorem}{section}

\@addtoreset{Proposition}{section}

\@addtoreset{Property}{section}

\@addtoreset{Assumption}{section}

\@addtoreset{Corollary}{section}

\@addtoreset{Lemma}{section}

\@addtoreset{Remark}{section}

\@addtoreset{Example}{section}

\newcommand{\ess}{\operatorname{ess}}
\newcommand{\esssup}{\ess\sup}

\newcommand{\N}{\mathbb{N}}

\newcommand{\Fc}{\mathcal{F}}

\newcommand{\Ic}{\mathcal{I}}
\newcommand{\Lc}{\mathcal{L}}

\newcommand{\tip}{{t_{i+1}}}
\newcommand{\rj}{{r_j}}
\newcommand{\rjp}{{r_{j+1}}}
\newcommand{\rjm}{{r_{j-1}}}

\newcommand{\D}{\mathbb{D}^{1,2}}

\newcommand{\tY}{{\widetilde Y}}

\newcommand{\Xp}{X^{\pi}}

\newcommand{\Yp}{Y^{\Re,\pi}}

\newcommand{\abs}[1]{\left\vert#1\right\vert}

\newcommand{\px}[1]{\partial_{x^{#1}}}
\newcommand{\py}[1]{\partial_{y^{#1}}}
\newcommand{\pz}[1]{\partial_{z^{#1}}}

\newcommand{\fA}{ \mathfrak{A} }
\newcommand{\fX}{ \mathfrak{X} }
\newcommand{\fY}{ \mathfrak{Y} }

\newcommand{\fV}{ \mathfrak{V} }
\newcommand{\fU}{ \mathfrak{U} }
\newcommand{\fG}{ \mathfrak{G} }

\def\sqw{\hbox{\rlap{\leavevmode\raise.3ex\hbox{$\sqcap$}}$%
\sqcup$}}
\def\sqb{\hbox{\hskip5pt\vrule width4pt height6pt depth1.5pt%
\hskip1pt}}

\def\qed{\ifmmode\hbox{\hfill\sqb}\else{\ifhmode\unskip\fi%
\nobreak\hfil
\penalty50\hskip1em\null\nobreak\hfil\sqb
\parfillskip=0pt\finalhyphendemerits=0\endgraf}\fi}
\def\cqfd{\ifmmode\sqw\else{\ifhmode\unskip\fi\nobreak\hfil
\penalty50\hskip1em\null\nobreak\hfil\sqw
\parfillskip=0pt\finalhyphendemerits=0\endgraf}\fi}
\DeclareUnicodeCharacter{00A0}{~}

\newcommand{\cA}{\mathcal{A}}

\newcommand{\cE}{\mathcal{E}}
\newcommand{\cF}{\mathcal{F}}

\newcommand{\cI}{\mathcal{I}}

\newcommand{\cK}{\mathcal{K}}
\newcommand{\cL}{\mathcal{L}}
\newcommand{\cM}{\mathcal{M}}
\newcommand{\cN}{\mathcal{N}}

\newcommand{\cP}{\mathcal{P}}
\newcommand{\cQ}{\mathcal{Q}}

\newcommand{\cU}{\mathcal{U}}
\newcommand{\cV}{\mathcal{V}}

\newcommand{\cY}{\mathcal{Y}}
\newcommand{\cZ}{\mathcal{Z}}

\newcommand{\E}{\mathbb{E}}

\renewcommand{\P}{\mathbb{P}}

\newcommand{\R}{\mathbb{R}}

\def \proof{{\noindent \bf Proof. }}
\def \eproof{\hbox{ }\hfill$\Box$}

\newcommand{\ud}{\mathrm{d}}

\newcommand{\HYP}[1]
    {\ensuremath{({H#1} ) }}

\newcommand{\1}{{\bf 1}}
    
\newcommand{\set}[1]
    {\ensuremath{\{ #1 \}}}
    
\newcommand{\HP}[1] 
    {\ensuremath{\mathscr{H}^{#1}}}

\newcommand{\esp}[1]{\ensuremath{\mathbb{E}  \left[#1\right] }}

\newcommand{\EFp}[2]
    {\ensuremath{
     \mathbb{E}_{#1} \left[#2\right] }}
     
\newcommand{\ti}[1]{t_{i #1}}
\newcommand{\tk}[1]{t_{k #1}}

\renewcommand{\Xi}[1]{X_{i #1}}

\newcommand{\YR}{Y^\Re}

\newcommand{\ZR}{Z^\Re}

\newcommand{\tYR}{\tilde{Y}^\Re}
\newcommand{\cYR}{\cY^\Re}

\newcommand{\cZR}{\cZ^\Re}
\newcommand{\cKR}{\cK^\Re}
\newcommand{\tcYR}{\tilde{\cY}^\Re}

\newcommand{\ckY}{\check{Y}}
\newcommand{\ckf}{\check{f}}
\newcommand{\ckK}{\check{K}}
\newcommand{\ckZ}{\check{Z}}
\newcommand{\ckU}{\check{U}}
\newcommand{\ckV}{\check{V}}
\newcommand{\cka}{\check{a}}
\newcommand{\ckG}{\check{\Gamma}}

\title{Rate of convergence for the discrete-time approximation of reflected  BSDEs  arising in switching problems}

\author{ 
          Jean-Fran\c{c}ois Chassagneux
         \\\small Laboratoire  de Probabilit\'{e}s  et Mod\`{e}les  Al\'{e}atoires
         \\\small CNRS, UMR 7599,  Universit\'{e} Paris Diderot
         \\\small \sf jean-francois.chassagneux@univ-paris-diderot.fr
		\and
             Adrien Richou
             \\\small  \small Univ. Bordeaux, IMB, UMR 5251, F-33400 Talence, France.
              \\\small  \sf adrien.richou@math.univ-bordeaux.fr 
              }

\begin{document}
\maketitle

\begin{abstract}
In this paper, we prove new convergence results improving the ones by Chassagneux, \'Elie and Kharroubi [\textit{Ann. Appl. Probab.} \textbf{22} (2012) 971--1007] for the discrete-time approximation of multidimensional obliquely reflected BSDEs. These BSDEs, arising in the study of switching problems, were considered by Hu and Tang [\textit{Probab. Theory Related Fields} \textbf{147} (2010) 89--121] and generalized by Hamad\`ene and Zhang [\textit{Stochastic Process. Appl.} \textbf{120} (2010) 403--426] and Chassagneux, \'Elie and Kharroubi [\textit{Electron. Commun. Probab.} \textbf{16} (2011) 120--128]. Our main result is a rate of convergence obtained in the Lipschitz setting and under the same structural conditions on the generator as the one required for the existence and uniqueness of a solution to the obliquely reflected BSDE. 

\end{abstract}

\vspace{5mm}


\vspace{5mm}

\noindent {\bf MSC Classification (2000):}  93E20, 65C99, 60H30.


 %

 %
 %
 

\section{Introduction}
In this paper, we study the discrete-time approximation of the following system of reflected backward stochastic differential equations
\begin{equation}\label{eqBSDECORIntro}\qquad
\begin{cases}
\displaystyle
Y_t = g(X_T)+\int_t^T f(X_s,Y _s, Z _s)\,\ud s-\int
_t^T Z_s \,\ud W_s  + K_T - K_t, &\quad 0 \le t \le
T,\vspace*{2pt}\cr
\displaystyle Y_t^\ell \ge \max_{j \in \cI} \set{Y_t^j -c^{\ell j}(X_t)} ,  &\quad
0\le t\le T,\; \ell \in \cI, \vspace*{2pt}\cr
\displaystyle\int_0^T\Bigl[Y _t^\ell-\max_{\textcolor{black}{j\in \Ic \setminus \set{\ell} } }\{Y ^j_t -
c^{\ell j}(X_t)\}\Bigr]
\,\ud K ^\ell_t=0,
&\quad \ell \in
\cI,
\end{cases}
\end{equation}
where $\Ic:=\{1,\ldots,d\}$, $f$, $g$ and $(c^{ij})_{i,j\in\Ic}$ are
Lipschitz functions and $X$ is  solution to the following forward stochastic
differential equation (SDE) with Lipschitz coefficients
\begin{equation}
\label{eqEDS}
 X_t = x +\int_0^t b(X_s)ds + \int_0^t \sigma(X_s)dW_s.
\end{equation}

An important motivation for this study  comes from economics applications, especially to energy markets. Indeed, it has been shown that the solution to
the above equations  allows to compute the solution of optimal switching problems which are linked to real option pricing (see e.g. \cite{carlud08}). This motivated a huge literature on switching problems both on the financial economics and applied mathematics sides, as pointed out in the introduction of \cite{hamzha10}.
The theoretical study of equation \eqref{eqBSDECORIntro} has started in dimension $2$ in the paper  \cite{hamjea07}
and was latter extended in higher dimension in \cite{djeham09, carlud08, portou09}. These studies are related to optimal 
switching problem and, in terms of existence and uniqueness result to \eqref{eqBSDECORIntro}, impose really strong conditions
on the driver $f$ of the BSDEs. These conditions were then weakened successively in \cite{huytan10, hamzha10, Chassagneux-Elie-Kharroubi-11}.
It is quite important to  notice that contrary to normally reflected BSDEs \cite{Gegout-Petit-Pardoux-96}, the best existence and uniqueness result available in the literature requires  structural conditions, see below, both on the driver $f$ and the function $c$. To the best of our knowledge, it can be found in the paper \cite{hamadene2013viscosity}.

The numerical study of \eqref{eqBSDECORIntro} by probabilistic methods has attracted much less attention \cite{portou09,elie2010probabilistic,Chassagneux-Elie-Kharroubi-10}. The first rate of convergence for a numerical scheme associated to \eqref{eqBSDECORIntro}
was proved in \cite{Chassagneux-Elie-Kharroubi-11} but under quite restrictive condition on the driver $f$. The main goal of our work is actually to prove a rate of convergence for a discrete-time scheme to obliquely reflected BSDEs under the same conditions on $f$ required to have existence and uniqueness
to \eqref{eqBSDECORIntro} and {only a} Lipschitz condition on the function $c$.

\vspace{10pt}
As in \cite{boucha08,majzha05,Chassagneux-Elie-Kharroubi-10}, we first introduce a discretely reflected version of \eqref{eqBSDECORIntro}, where the reflection occurs only on a deterministic
grid $\Re =\{r_{0} := 0,\ldots, r_{\kappa}:=T\}$:
 $Y^{\Re}_T = \tY^{\Re}_T := g(X_{T}) \in\cQ(X_{T})$, and, for $j
\le
\kappa-1$ and
$t \in[\rj,\rjp)$,
%
%
\begin{equation} \label{BSDEDORintro}
\begin{cases}
\displaystyle \tY^{\Re}_{t } = Y^{\Re}_\rjp+ \int_{t}^\rjp f(X_{u},
\tY^{\Re}_u, Z^{\Re}_u)\,\ud u - \int_{t}^\rjp Z^{\Re}_u \,\ud W_u,
\vspace*{2pt}\cr
\displaystyle Y^{\Re}_{t}= \tY^{\Re}_t\mathbf{1}_{ \{ t \notin\Re\} } +
\cP(X_t,\tY^{\Re}_t)\mathbf{1}_{ \{ t \in\Re\} } ,
\end{cases}
\end{equation}
where $\cP(x,.)$ is the oblique projection operator on the closed convex domain 
\begin{equation*}
 \mathcal{Q}(x):=\left\{ y \in \mathbb{R}^d | y^i \ge \max_{j \in \cI} (y^j - c^{ij}(x)), \forall i \in \cI\right\},
\end{equation*}
defined by
\begin{equation*}
 \cP:(x,y) \in \mathbb{R}^d \times \mathbb{R}^d \mapsto \left(\max_{j \in \cI}\set{y^j -c^{ij}(x)}\right)_{1 \le i \le d}.
\end{equation*}
{Let us remark that \eqref{BSDEDORintro} rewrites equivalently for $t \in [0,T]$ as 
\begin{equation} \label{BSDEDORintro equivalente}
\begin{cases}
\displaystyle \tY^{\Re}_{t } = Y^{\Re}_\rjp+ \int_{t}^\rjp f(X_{u},
\tY^{\Re}_u, Z^{\Re}_u)\,\ud u - \int_{t}^\rjp Z^{\Re}_u \,\ud W_u+(K^{\Re}_T - K^{\Re}_t),
\vspace*{2pt}\cr
\displaystyle K^{\Re}_{t}:= \sum_{r \in \Re \setminus \{0\}} \Delta K_r^{\Re} \mathbf{1}_{\{ r \leqslant t\} }, \text{ with } \Delta K_t^{\Re}:= Y_t^{\Re} - \tY^{\Re}_t = -(\tY^{\Re}_t - \tY^{\Re}_{t^-}). 
\end{cases}
\end{equation}}
We denote $|\Re|$ the modulus of $\Re$ given by $ |\Re|:=\max_{0 \le i \le \kappa-1} |r_{i+1}-r_i|$.

An important step in our study is to prove that these discretely reflected BSDEs are a good approximation
of the continuously reflected ones \eqref{eqBSDECORIntro}. In Section 4, we are able to control the error in terms of
$|\Re|$ under a {simple} Lipschitz condition for the cost functions $c$, which is new in the literature. {Finally, the main result of this section is given by Theorem \ref{th convergence continuement ref vers discret refl}. It improves in particular the results of \cite{Chassagneux-Elie-Kharroubi-10} where a less general structure on $f$ is imposed, $f$ is assumed to be bounded with respect to $z$ and some extra smoothness assumptions are needed on $c$. Moreover error bounds obtained in \cite{Chassagneux-Elie-Kharroubi-10} are slightly improved.}

{Let us remark that we were able to get away from some structural contraints on $f$ imposed in \cite{Chassagneux-Elie-Kharroubi-10} by adapting arguments of \cite{Chassagneux-Elie-Kharroubi-11} in our context (see in particular the proof of Proposition \ref{pr theorem 5.2}). Moreover, we get ride of the boundedness of $f$ by obtaining the following strong estimate on $Z$ in Corollary \ref{co control Z}:
\begin{align*}
|Z_t| \le C(1 + |X_t|) \quad \ud \P\otimes \ud t \; \textrm{a.e.}
\end{align*}
Let us point out that this kind of estimate is interesting for its own sake, that it should be easily generalised in a path-dependent framework and that the proof strongly rest on the discretely reflected BSDE approximation of \eqref{eqBSDECORIntro}.
}

\vspace{10pt}
We then consider a Euler type approximation scheme associated to the
BSDE~(\ref{BSDEDORintro}) defined on a grid $\pi=\{t_0,\ldots,t_n\}$ by
$\Yp_n:=g(\Xp_T)$ and, for $i\in\{n-1,\ldots,0\}$,
%
%
\begin{equation}\label{schemeintro}
\begin{cases}
Z_i^{\Re,\pi} := \mathbb{E}[Y^{\Re,\pi}_{i+1}H_i \mid\mathcal{F}_{t_i}] ,
\vspace*{2pt}\cr
\widetilde{Y}^{\Re,\pi}_{i} := \mathbb{E}[Y^{\Re,\pi}_{i+1}\mid\mathcal{F}_{t_i
}] + h_i
f(\Xp_{t_i},\widetilde{Y}^{\Re,\pi}_{i}, Z_i^{\Re,\pi}) ,
\vspace*{2pt}\cr
Y_i^{\Re,\pi} := \widetilde{Y}_i^{\Re,\pi} \mathbf{1}_{ \{ t_i\notin\Re\} } + \cP
(\Xp_{t_i},
\widetilde Y_i^{\Re,\pi})\mathbf{1}_{ \{ t_i\in\Re\} } ,
\end{cases}
\end{equation}
where $\Xp$ is the Euler scheme associated to $X$, $h_i:=t_{i+1}-t_i$ and weights $(H_i)_{0 \le i \le n-1}$ are matrices in $\mathcal{M}^{1,d}$ given by
\begin{align}
\label{eq def HiR intro}
 (H_i)^{\ell} = \frac{-R}{h_i} \vee \frac{W^\ell_{t_{i+1}} -W^\ell_{t_i}}{h_i} \wedge \frac{R}{h_i}, \quad 1 \le \ell \le d,
\end{align}
 with $R$ a positive parameter.
We denote $|\pi|$ the modulus of $\pi$ given by $ |\pi|:=\max_{0 \le i \le n-1} h_i$ and we assume that we always have $\Re \subset \pi$.

\vspace{10pt}
To obtain our convergence results, we work, throughout this paper, under the following assumption:\\
\HYP{f}  \begin{enumerate}[(i)]
 \item The functions $\sigma : \mathbb{R}^d \rightarrow \mathcal{M}^{d,d}$ and $b : \mathbb{R}^d \rightarrow \mathbb{R}^d$ are Lipschitz-continuous functions.
 \item The functions $f : \mathbb{R}^d \times \mathbb{R}^d \times \mathcal{M}^{d,d} \rightarrow \mathbb{R}^d $, $g : \mathbb{R}^d \rightarrow \mathbb{R}^d$ and $(c^{ij}: \mathbb{R}^d \rightarrow \mathbb{R})_{i,j \in \cI}$ are Lipschitz-continuous functions and $f^j(x,y,z)=f^j(x,y,z^{j.})$ {where $z^{j.}$ stands for the $j$th row of $z$}. We denote by $L^Y$ and $L^Z$ the Lipschitz constants of $f$ with respect to $y$ and $z$.
 \item $g(x) \in \mathcal{Q}(x)$, for all $x \in \mathbb{R}^d$.
 \item The cost functions $(c^{ij})_{i,j \in \cI}$ satisfy the following structure condition
 \begin{equation}
  \label{condition de structure c}
  \begin{cases}
   c^{ii}=0, & \textrm{for } 1 \le i \le d;\\
   \inf_{x \in \mathbb{R}^d} c^{ij}(x) >0, & \textrm{for } 1 \le i,j \le d \textrm{ with } i \neq j;\\
   \inf_{x \in \mathbb{R}^d} \{c^{ij}(x)+c^{jl}(x) - c^{il}(x)\} >0, & \textrm{for } 1 \le i,j \le d \textrm{ with } i \neq j, j \neq l.
  \end{cases}
 \end{equation}
\end{enumerate}
 Let us emphasize here the fact that our results are obtained without any assumption on the non-degeneracy of the volatility matrix $\sigma$.
 We also point out that $\HYP{f}$(ii) is the best condition -- up to now -- for existence and uniqueness to \eqref{eqBSDECORIntro} to hold. {Let us remark that all the results contained in [8] need a stronger structure assumption on $f$ since authors also assume that, for all $j \in \cI$, $f^j(x,y,z) =  f^j(x,y^j,z^{j.})$.}
 
 \vspace{5pt}
 A fundamental result to obtain convergence for continuously reflected BSDEs is first to prove that the scheme given in
 \eqref{schemeintro} approximates efficiently discretely reflected BSDEs. This result is interesting in itself if one is only interested
 in the approximation of Bermudan switching problem (i.e. when the switching times are restricted to lie in the grid $\Re$). It is discussed in Section 3 below and requires, in particular, the use of a new representation result for the scheme  \eqref{schemeintro}. {More precisely, the new representation result is given by Proposition \ref{prop snell representation schema} and an upper bound on the discrete-time approximation error is obtained in Theorem \ref{th convergence scheme vers discret refl}. Let us remark that Proposition \ref{prop snell representation schema} is crucially based on the comparison result for approximation schemes for BSDEs obtained in \cite{Chassagneux-Richou-14}, its explains why there are some bounded weights $(H_i)_{0 \leqslant i \leqslant n-1}$ in the scheme \eqref{schemeintro}. Moreover, Proposition \ref{prop snell representation schema} allows to obtain a general stability result for schemes (see Proposition \ref{prop stabilite schema gene aleatoire}) which is the 
keystone result that was missing in \cite{Chassagneux-Elie-Kharroubi-10} and which allows to strongly improve their upper bound on the approximation error when the generator depends on $z$.}

\vspace{10pt}
Combining the fact that discretely reflected BSDEs \eqref{BSDEDORintro} are a good approximation of continuously reflected BSDEs and that the scheme
 \eqref{schemeintro}  is also a good approximation of \eqref{BSDEDORintro}, we obtain our new convergence result, which is the main result of this paper and is summarised in the following Theorem.

\begin{Theorem}
\label{th convergence schema vers solution continuement reflechie}
Let us assume that \HYP{f} is in force. Set $R$ such that $L^Z R \le 1$, $\pi$ such that $L^Y |\pi| < 1$ and define $\alpha(|\pi|)=\log (2T/|\pi|)$. Then the following holds, for some positive constant $C$:
\begin{enumerate}[(i)]
 \item Taking $|\Re| \sim |\pi|^{1/2}$, we have
$$  \sup_{0 \le i \le n} \esp{|Y_{t_i} -\widetilde{Y}^{\Re,\pi}_i|^2 + |{Y}_{t_i} -{Y}^{\Re,\pi}_i|^2}\le C |\pi|^{1/2} \alpha(|\pi|).$$
 \item Taking $|\Re| \sim |\pi|^{1/3}$, we have
$$  \sup_{0 \le i \le n} \esp{|Y_{t_i} -\widetilde{Y}^{\Re,\pi}_i|^2 + |{Y}_{t_i} -{Y}^{\Re,\pi}_i|^2}\le C |\pi|^{1/3}\alpha(|\pi|),$$
and
$$ \esp{ \sum_{i=0}^{n-1}\int_{t_i}^{t_{i+1}} |Z_s - Z_i^{\Re,\pi}|^2ds}\le C |\pi|^{1/6} \sqrt{\alpha(|\pi|)}.$$
\end{enumerate}
Moreover, if the cost functions $c$ are constant, then the previous estimates remain true with $\alpha(|\pi|):=1$.
\end{Theorem}

It is important to compare the previous result with Theorem 5.4 in \cite{Chassagneux-Elie-Kharroubi-10} which gives also
rates of convergence for the discrete-time approximation of obliquely reflected BSDEs. Up to a slight modification of the scheme {(introducing the
truncation of the Brownian increments by \eqref{eq def HiR intro} whereas in \cite{Chassagneux-Elie-Kharroubi-10} authors consider the same kind of scheme but where weights $(H_i)_{0 \leqslant i \leqslant n-1}$ are given by the unbounded Brownian increments)}, we see that we are able to obtain the convergence rate $1/4$ for $Y$. 
{However,  Theorem 5.4 in \cite{Chassagneux-Elie-Kharroubi-10} only gives a logarithmic convergence rate by assuming more restrictive assumptions on $f$ and $c$: for all $j \in \cI$ they need to have $f^j(x,y,z) = f^j(x,y^j,z^{j.})$, $f$ is assumed to be bounded with respect to $z$ and $c$ must be a difference of two $C^2$ functions with some boundedness assumptions on derivatives. Nevertheless, they are also able to get a convergence rate $1/2$ for $Y$ when, for all $j \in \cI$, $f^j(x,y,z) = f^j(x,y^j)$.}


\vspace{10pt}
The rest of the paper is organised as follows.
In Section 2, we present preliminary results that will be useful in the rest of the paper. We discuss the representation property of obliquely reflected BSDEs in terms of auxiliary one-dimensional BSDEs. We also give new regularity results for the discretely reflected BSDEs which are key tools to obtain our convergence results. Section 3 is devoted to the study of the Euler numerical scheme, in particular its fundamental stability property. Using this stability property and the regularity results given in Section 2, {we obtain in Theorem \ref{th convergence scheme vers discret refl}} a control of the error between the scheme and the discretely reflected BSDEs. Section 4 is concerned with the approximation of continuously reflected BSDEs by the discretely reflected ones. A convergence rate is obtained in {Theorem \ref{th convergence continuement ref vers discret refl}} that makes possible to prove, using the result of Section 3, our main result, Theorem \ref{th convergence schema 
vers solution continuement reflechie} above.
For the reader convenience, some technical proofs are postponed in an Appendix Section.

\paragraph{Notations}
Throughout this paper we are given a finite time horizon $T$ and a probability space $(\Omega, \mathcal{F},\mathbb{P})$ endowed with a $d$-dimensional standard Brownian motion $(W_t)_{t \ge 0}$. The filtration $(\mathcal{F}_t)_{t \le T}$ is the Brownian filtration. $\mathscr{P}$ denotes the $\sigma$-algebra on $[0,T] \times \Omega$ generated by progressively measurable processes. Any element $x \in \mathbb{R}^n$ will be identified to a column vector with $i$th component $x^i$ and Euclidean norm $|x|$. For $x,y \in \mathbb{R}^n$, $x.y$ denotes the scalar product of $x$ and $y$. We denote by $ \preccurlyeq $ the component-wise partial ordering relation on vectors. $\cM^{n,m}$ denotes the set of real matrices with $n$ lines and $m$ columns. For a matrix $M \in \cM^{n,m}$, $M^{ij}$ is the component at row $i$ and column $j$, $M^{i.}$ is the $i$th row and $M^{.j}$ the $j$th column. 

We denote by $\mathcal{C}^{k,b}$ the set of functions with continuous and bounded derivatives up to order $k$. 
For a function $f:\R^n \rightarrow \R,\; x \mapsto f(x)$, we denote by $\px{}f = (\px{1} f, \dots,\px{n} f)$. 
If $f:\R^n\times\R^d \rightarrow \R,\, (x,y) \mapsto f(x,y)$ we denote $\px{}f$ (resp. $\py{} f$) the derivatives with respect to the variable $x$ (resp. $y$). For $g:\R^n \mapsto \R^d,\, x \rightarrow g(x)$, $\px{} g$ is a matrix and $(\px{} g)^{i.}= \px{} g^i$.

For ease of notation, we will sometimes write $\mathbb{E}_t[.]$ instead of $\mathbb{E}[.|\mathcal{F}_t]$, $t \in [0,T]$. 
Finally, for any $p \geq 1$, we introduce the following: 
\begin{itemize}
\item[-] $\mathscr{L}^p$ the set of $\cF_T$-measurable random variables $G$ satisfying $|G|_{\mathscr{L}^p}:= \esp{|G|^p}^\frac1p<+\infty$,
\item[-]  $\mathscr{S}^p$ the set of càdlàg adapted processes $U$ satisfying 
$$|U|_{\mathscr{S}^p}:=\esp{\sup_{t \in [0,T]} |U_t|^p }^\frac1p< \infty,$$ 
and $\mathscr{S}^p_c$ the subset of continuous processes in $\mathscr{S}^p$,
\item[-] $\mathscr{H}^p$ the set of progressively measurable processes $V$ satisfying 
$$|V|_{\mathscr{H}^p}:=\esp{ \left(\int_0^T |V_t|^2 \ud t \right)^{\frac{p}2} }^\frac1p< \infty,$$
\item[-]  $\mathscr{K}^p$ the set of continuous non-decreasing processes in $\mathscr{S}^p$,
\item[-] $\mathscr{K}^{\Re,p}$ the set of pure jump non-decreasing processes in $\mathscr{S}^p$ with jump times in $\Re$.
\end{itemize}

In the sequel, we denote by $C$ a constant whose value may change from line to line but which never depends on $|\pi|$ nor $|\Re|$. The notation $C_{\alpha}$ is used to stress the fact that the constant depends on some parameter $\alpha$.

 \section{Preliminary results}

In this section, we present key properties of continuously and discretely reflected BSDEs. 
We start by recalling the representation property in terms of ``switched'' BSDEs of the multidimensional
systems of reflected BSDEs \eqref{eqBSDECORIntro}  or \eqref{BSDEDORintro}. 

In a second part, we study the regularity properties of the solution to discretely reflected BSDEs in a Markovian setting.
These results are key tools to obtain a convergence rate for the numerical approximation. They are new in the framework of
this paper but their proofs rely on arguments that are now quite well understood.

\subsection{Representation  of obliquely reflected BSDEs}
\label{section representation property}

As mentioned in the introduction, the motivation to work on the above class of obliquely reflected BSDEs comes 
from the study of ``switching problems'' in the financial economics literature.
Indeed, RBSDEs provide a characterization of the solution to these switching problems.
Interestingly, the interpretation of the RBSDE in term of the solution of a ``switching problem'' is
a key tool in our work. We now recall the link between the two objects, which takes the form of a representation
theorem for the solution of the RBSDEs in terms of ``switched BSDEs''. This link has been established before, see e.g. \cite{huytan10}.
We state it here in a generic framework as this will be useful latter on.

\vspace{10pt}
We consider a matrix valued process $C=(C^{ij})_{1 \le i,j \le n}$ such that $C^{ij}$ belongs to $\mathscr{S}^2$ for $i,j \in \cI$ and satisfies the structure condition 
 \begin{equation}
  \label{condition de structure C}
  \begin{cases}
   C_t^{ii}=0, & \textrm{for } 1 \le i \le d \textrm{ and } 0 \le t \le T;\\
   \inf_{t \in [0,T]} C^{ij}_t \ge \varepsilon >0, & \textrm{for } 1 \le i,j \le d \textrm{ with } i \neq j;\\
   \inf_{t \in [0,T]} \{C^{ij}_t+C^{jl}_t - C^{il}_t\} >0, & \textrm{for } 1 \le i,j \le d \textrm{ with } i \neq j, j \neq l.
  \end{cases}
 \end{equation}
We introduce a random closed convex set family associated to $C$:
$$\mathcal{Q}_t := \left\{ y \in \mathbb{R}^d | y^i \ge \max_{j\in \cI} (y^j -C^{ij}_t), 1 \le i \le d \right\}, \quad 0 \le t \le T,$$
and the oblique projection operator onto $\mathcal{Q}_t$, denoted $\mathcal{P}_t$ and defined by
\begin{equation}
\label{eq def projection cout C}
\cP_t : y \in \R^d \mapsto \left(\max_{j \in \cI} \{ y^j - C_t^{ij} \}\right)_{1 \le i \le d}.
\end{equation}
\begin{Remark}
\label{re monotonie proj}
{It follows from the structure condition \eqref{condition de structure C} that, for all $t \in [O,T]$, $0  \in \mathcal{Q}_t \neq \emptyset$. It implies that $\mathcal{P}_t$ is well defined. Moreover, we can easily check that $\cP_t$ is increasing with respect to the partial ordering relation $\preccurlyeq$.}
\end{Remark}

\vspace{2mm}
A switching strategy $a$ is a nondecreasing sequence of stopping times $(\theta_j)_{j\in\N}$ , combined with a sequence of random variables $(\alpha_j)_{j\in\N}$ valued in $\Ic$, such that $\alpha_j$ is $\Fc_{\theta_j}-$measurable, for any $j\in\N$. 
We denote by $ \mathscr{A}$ the set of such strategies.
For $a= (\theta_{j},\alpha_{j})_{j\in\N}\in\mathscr{A}$, we introduce $\cN^a$ the (random) number of switches before $T$:
\begin{align}
\cN^a  =  \#\{k\in\N^*~:~\theta_{k} \leq T \}\;. \label{eq def Na}
\end{align}
To any switching strategy $a=(\theta_{j},\alpha_{j})_{j\in\N}\in\mathscr{A}$, we
associate the current state process $(a_t)_{t\in[0,T]}$ and the cumulative cost process $(\cA^a_t)_{t\in[0,T]}$ defined respectively by
  \begin{align*}
  a_t \;:=\; \alpha_{0}\1_{\{0\leq t<\theta_0\}}+\sum_{j=1}^{\cN^a}  \alpha_{j-1} \1_{\{\theta_{j-1}\leq t < \theta_{j}\}}  \;\text{ and }\;  
  \cA^a_t \;:=\; \sum_{j=1}^{\cN^a} C^{\alpha_{j-1}\alpha_{j}}_{\theta_{j}} \1_{\{\theta_j\le t \le T\}}\,, 
  \end{align*}
for $0\le t\le T$. {Since $\mathcal{N}^a$ appears in the definition of processes $(a_t)_{t \in [0,T]}$ and $(\cA^a_t)_{t\in[0,T]}$, it is not clear at first sight that these two processes are adapted. Nevertheless, we can remark that we have equivalently
  \begin{align*}
  a_t \;:=\; \alpha_{0}\1_{\{0\leq t<\theta_0\}}+\sum_{j=1}^{\cN^a_t}  \alpha_{j-1} \1_{\{\theta_{j-1}\leq t < \theta_{j}\}}  \;\text{ and }\;  
  \cA^a_t \;:=\; \sum_{j=1}^{\cN^a_t} C^{\alpha_{j-1}\alpha_{j}}_{\theta_{j}} \1_{\{\theta_j\le t \le T\}}\,, \quad 0\le t\le T,
  \end{align*}
with $\cN^a_t  =  \#\{k\in\N^*~:~\theta_{k} \leq t \}$.
}  
  
\begin{Remark} \label{re swithing strategy} 
\begin{enumerate}[(i)]
 \item The sequence of stopping times is only supposed to be non-decreasing, but the assumptions
 on the cost processes \eqref{condition de structure C} imply that any reasonable strategy uses a sequence of increasing stopping times, apart from a possible instantaneous switch at initial time.
 This is specially the case for the optimal strategies. 
 \item Note that the cumulative cost process will keep track of all the switching times, even the instantaneous ones;
 whereas the state process will keep track of the last state when instantaneous switches occur.
 \end{enumerate}
\end{Remark}

For $(t,i)\in[0,T]\times\cI$, the set $\mathscr{A}_{t,i}$ of admissible strategies starting from state $i$ at time $t$ is defined by
\begin{align*}
\mathscr{A}_{t,i} = \set{a=(\theta_{j},\alpha_{j})_{j}\in \mathscr{A}~| \theta_{0}= t,~\alpha_{0}= i,~\esp{|\mathcal{A}_{T}^a|^2}< \infty}\;,
\end{align*}

similarly we introduce $\mathscr{A}^\Re_{t,i}$ the restriction to $\Re-$admissible strategies
 \begin{align*}
\mathscr{A}^\Re_{t,i} := \left\{~ a=(\theta_j,\alpha_j)_{j\in\N}  \in \mathscr{A}_{t,i}\quad|\quad  \theta_{j}\in\Re~, ~\forall j\leq \cN^a   ~\right\} \;,
 \end{align*}
 and denote $\mathscr{A}^\Re:=\bigcup_{i\le d} \mathscr{A}^\Re_{0,i}$.

For a strategy $a \in \mathscr{A}_{t,\ell}$, we introduce  the one-dimensional \emph{switched BSDE} whose solution $(\cU^a,\cV^a)$ satisfies
\begin{align}\label{eq de switched BSDE}
\cU^a_t = \xi^{a_T} + \int_t^TF^{a_s}(s,\cV^a_s)\ud s -\int_t^T\cV^a_s \ud W_s - \cA^a_T+ \cA^a_t\;
\end{align}
where the terminal condition $\xi$ and the random driver $F$ satisfy following assumptions, for some $p\geqslant 2$:

\HYP{F_p}
\begin{enumerate}[(i)]
 \item $F:\Omega \times [0,T]\times \mathcal{M}^{d, d} \rightarrow \mathbb{R}^d$ is $\mathscr{P} \otimes \mathcal{B}(\mathbb{R}^d) \otimes \mathcal{B}(\mathcal{M}^{d,d})$-measurable,
 \item $F^j(\cdot,z) =  F^j(\cdot,z^{j.})$ for all $j \in \cI$,
 \item $|F(s,z)-F(s,z')| \le C|z-z'|$ for all $s \in [0,T]$, $z,z' \in \mathcal{M}^{d, d}$,
 \item $\xi$ is $\mathcal{F}_T$-measurable and is valued in $\mathcal{Q}_T$,
 \item $\mathbb{E}\left[|\xi|^p+ \int_0^T \abs{F(s,0)}^p ds \right] \le C_p$.
\end{enumerate}
{Let us remark that $\HYP{F_p} \Rightarrow \HYP{F_{p'}}$ when $p \geqslant p'$.}
We now define multidimensional processes $\bar{\cY}$ and $\bar{\cY}^\Re$ as follows, for $\ell \in \set{1, \dots,d}$ 
{
\begin{align*}
(\bar{\cY}_t)^\ell := \esssup_{a \in \mathscr{A}_{t,\ell} } \left(\cU^a_t -\cA^a_t\right)  \; \text{ and } \;
(\bar{\cY}^\Re_t)^\ell := \esssup_{a \in \mathscr{A}^\Re_{t,\ell} } \left(\cU^a_t - \cA^a_t\right)  \;.
\end{align*}
}
The process $\bar{\cY}$ represents the optimal value that can be obtained from the switched BSDEs following strategies in $\mathscr{A}$.
The process $\bar{\cY}^\Re$ can be seen as a ``Bermudan'' version of it, i.e.  when the switching times are restricted to lie in $\Re$.
Both processes enjoy a representation in terms of reflected BSDEs, the main difference lying into the reflecting process that for the latter
will be a pure jump process with jump times in $\Re$.

\vspace{2mm}

Let $(\cY,\cZ,\cK)$   be the solution to the following BSDE
\begin{equation}\label{eq th rep switch generic}\qquad
\begin{cases}
\displaystyle
\cY^\ell_t = \xi^\ell +\int_t^T F^\ell(s, \cZ _s)\,\ud s-\int
_t^T \cZ
^\ell_s \,\ud W_s  + \cK ^\ell_T - \cK ^\ell_t, &\quad 0 \le t \le
T,\;\; \ell \in \cI,\vspace*{2pt}\cr
\displaystyle \cY ^\ell_t \geq\max_{j\in\Ic} \{ \cY ^j_t - C^{\ell j}_t\}, &\quad
0\leq
t\leq
T,\;\; \ell \in \cI,\vspace*{2pt}\cr
\displaystyle\int_0^T\Bigl[ \cY _t^\ell-\max_{j\in\Ic\setminus \set{\ell}}\{ \cY ^j_t -
C^{\ell j}_t\}\Bigr]
\,\ud \cK^\ell_t=0,
&\quad \ell \in
\cI,
\end{cases}
\end{equation}
and 
$(\tcYR, \cYR, \cZR,\cKR)$ with $\cYR_{t}=\tcYR_{t -} $, $t \in (0,T]$ be the solution of following discretely reflected BSDEs,
\begin{equation}\label{eq th rep switch generic disc}\qquad
\begin{cases}
\displaystyle
\tcYR_t = \xi +\int_t^T F(s, \cZR _s)\,\ud s-\int
_t^T \cZR_s \,\ud W_s  + \cKR_T - \cKR_t, &\quad 0 \le t \le
T,\vspace*{2pt}\cr
\displaystyle
\cYR_r \in \cQ_r, &\quad
r \in \Re
,\vspace*{2pt}\cr
\displaystyle\int_0^T\Bigl[ (\cYR_{t})^\ell-\max_{j\in\Ic\setminus \set{\ell} }\{ (\cYR_{t})^\ell -
C^{\ell j}_t\}\Bigr]
\,\ud (\cKR_t)^\ell=0,
&\quad \ell \in
\cI,
\end{cases}
\end{equation}
Existence and uniqueness of a solution for equation \eqref{eq th rep switch generic} has been addressed in \cite{huytan10,hamzha10} and in \cite{Chassagneux-Elie-Kharroubi-10} (Proposition 2.1) for equation \eqref{eq th rep switch generic disc}. For the reader convenience we recall here these results.
\begin{Proposition}
 \label{prop existence unicite EDSRs obliqu refle cont et discr avec generateur aleatoire}
 Assume that \HYP{F_p} holds for some $p \geqslant 2$. There exists a unique solution $( \cY,\cZ,\cK)\in \mathscr{S}^2_c \times \mathscr{H}^2 \times \mathscr{K}^2$ to \eqref{eq th rep switch generic} and a unique solution $(\tcYR,\cYR,\cZR, \cKR)$ with $(\tcYR,\cZR, \cKR)\in \mathscr{S}^2 \times \mathscr{H}^2 \times \mathscr{K}^{\Re,2}$ to \eqref{eq th rep switch generic disc}. They also satisfy
 $$|\cY|_{\mathscr{S}^p}+ |\cZ|_{\mathscr{H}^p}+ |\cK_T|_{\mathscr{L}^p} \leq C_p 
 \;\text{ and }\;
 |\tcYR|_{\mathscr{S}^p}+ |\cZR|_{\mathscr{H}^p}+ |\cKR_T|_{\mathscr{L}^p} \leq C_p   .$$
\end{Proposition}

Gathering Proposition 3.2 in \cite{Chassagneux-Elie-Kharroubi-11} and Theorem 2.1 in \cite{Chassagneux-Elie-Kharroubi-10}, we have the following key representation result.

\begin{Proposition}\label{th rep switch generic}
 Assume that \HYP{F_2} is in force. The following holds:
 \begin{enumerate}[(i)]
  \item For all $\ell \in \set{1, \dots, d}$, $t \in [0,T]$,
   \begin{align*}
(\cY_t)^\ell = (\bar{\cY}_t)^\ell=  \cU^{\bar{a}}_t {-\cA_t^{\bar{a}}} \quad \text{ and } \quad {(\cYR_t)^\ell} = (\bar{\cY}^\Re_t)^\ell=  \cU^{\bar{a}^\Re}_t {-\cA_t^{\bar{a}^\Re}}
\end{align*}
for some $\bar{a} \in \mathscr{A}_{t,\ell} $ and $\bar{a}^\Re \in \mathscr{A}^\Re_{t,\ell} $.
 \item The strategy $\bar{a}= (\bar{\theta}_j,\bar{\alpha}_j)_{j \ge 0}$ can be defined recursively by $(\bar{\theta}_0,\bar{\alpha}_0):=(t,\ell)$ and, for $j \ge 1$,
 $$\bar{\theta}_j := \inf \left\{ s \in [\bar{\theta}_{j-1},T] \bigg| (\tilde{\cY}_s)^{\bar{\alpha}_{j-1}} \le \max_{k \neq \bar{\alpha}_{j-1}} \{(\tilde{\cY}_s)^k - C_s^{\bar{\alpha}_{j-1} k }\} \right\},$$
 $$\bar{\alpha}_{j} := \min \left\{ \ell \neq \bar{\alpha}_{j-1} \bigg| (\tilde{\cY}_{\bar{\alpha}_{j}})^{\ell} - C_{\bar{\theta}_j}^{\bar{\alpha}_{j-1} \ell} = \max_{k \neq \bar{\alpha}_{j-1}} \{(\tilde{\cY}_{\bar{\theta}_j})^k - C_{\bar{\theta}_j}^{\bar{\alpha}_{j-1} k }\} \right\}.$$
 \item The strategy $\bar{a}^\Re= (\bar{\theta}^\Re_j,\bar{\alpha}^\Re_j)_{j \ge 0}$ can be defined recursively by $(\bar{\theta}^\Re_0,\bar{\alpha}^\Re_0):=(t,\ell)$ and, for $j \ge 1$,
 $$\bar{\theta}^\Re_j := \inf \left\{ s \in [\bar{\theta}_{j-1}^\Re,T] \cap \Re \bigg| (\tilde{\cY}_s^\Re)^{\bar{\alpha}^\Re_{j-1}} \le \max_{k \neq \bar{\alpha}^\Re_{j-1}} \{(\tilde{\cY}_s^\Re)^k - C_s^{\bar{\alpha}^\Re_{j-1} k }\} \right\},$$
 $$\bar{\alpha}^\Re_{j} := \min \left\{ \ell \neq \bar{\alpha}^\Re_{j-1} \bigg| (\tilde{\cY}_{\bar{\theta}^\Re_{j}}^\Re)^{\ell} - C_{\bar{\theta}^\Re_j}^{\bar{\alpha}^\Re_{j-1} \ell} = \max_{k \neq \bar{\alpha}^\Re_{j-1}} \{(\tilde{\cY}_{\bar{\theta}^\Re_j}^\Re)^k - C_{\bar{\theta}^\Re_j}^{\bar{\alpha}^\Re_{j-1} k }\} \right\}.$$
 \end{enumerate}
\end{Proposition}

\begin{Remark} \label{re rep switched bsdes}
 If ${Y}_t^\ell \in \partial \mathcal{Q}_t$ then there is an instantaneous jump, i.e. $\bar{\theta}_1=t$. In the same way, if $t \in \Re$ and $({Y}_t^\Re)^\ell \in \partial\mathcal{Q}_t$ then $\bar{\theta}^\Re_1=t$.
\end{Remark}

{
\begin{Remark}
  Proposition \ref{th rep switch generic} does not correspond exactly to Proposition 3.2 in \cite{Chassagneux-Elie-Kharroubi-11} and Theorem 2.1 in \cite{Chassagneux-Elie-Kharroubi-10}. Indeed, Sa\"id Hamad\`ene pointed us out a mistake in the representation theorem in \cite{huytan10} (Theorem 3.1) where the term $\cA_t^{\bar{a}}$ in (i) is missing. Let us emphasise that this missing term is due to the possible instantaneous jump when the solution is on the boundary of the domain (c.f. Remark \ref{re rep switched bsdes}). Since Proposition 3.2 in \cite{Chassagneux-Elie-Kharroubi-11} strongly relies on the proof of Theorem 3.1 in \cite{huytan10}, the same error occurs. The same applies also for the proof of Theorem 2.1 in \cite{Chassagneux-Elie-Kharroubi-10}: the term $\cA_t^{\bar{a}^\Re}$ is also missing in \cite{Chassagneux-Elie-Kharroubi-10} and moreover the representation theorem is incorrectly obtained for $\tcYR$ instead of $\cYR$. Lastly, let us remark that the representation for $\mathcal{Y}$ was 
obtained in \cite{huytan10,Chassagneux-Elie-Kharroubi-11} for switching strategies with an increasing sequence of switching times. For more details on proofs and necessary corrections, we refer the reader to the proof of Proposition \ref{prop snell representation schema} where we have obtained the same kind of representation for obliquely reflected backward schemes.
\end{Remark}
}

 %
 %

 
\subsection{Discretely obliquely reflected BSDEs in a Markovian setting}
We will now study the discretely obliquely reflected BSDEs \eqref{eq th rep switch generic disc} in a Markovian setting, namely the solution to \eqref{BSDEDORintro}. We will in particular prove regularity results for this process.
The main difference with Section 3 in \cite{Chassagneux-Elie-Kharroubi-10} comes from the assumption on $f$, in particular the full dependence in the $y$-variable, recall \HYP{f}(ii). 

\vspace{10pt} \noindent
Let us recall that under assumption \HYP{f}(i), there exists a unique strong solution to the SDE \eqref{eqEDS} which satisfies
\textcolor{black}{
\begin{align} \label{eq control X} \EFp{t}{\sup_{s \in
[t,T]}|X_s|^p} \le C_p(1+|X_t|^p) \;,\quad \text{ for all } p \ge 2\;, \; t\in[0,T]\;.
\end{align}
}

\subsubsection{Basic properties}

The following proposition gives some useful estimates on the solution to \eqref{BSDEDORintro}. Its proof is postponed to the Appendix.

 \begin{Proposition}
 \label{prop estimee Y Z K}
 Assume that \HYP{f} is in force. There exists a unique solution $(\tilde{Y}^{\Re},Y^{\Re},Z^{\Re}) \in \mathscr{S}^2 \times \mathscr{S}^2 \times \mathscr{H}^{2}$ to \eqref{BSDEDORintro}, or equivalently \eqref{BSDEDORintro equivalente},  and it satisfies, for all $p \geq 2$, 
 $$|\tilde{Y}^{\Re}|_{\mathscr{S}^p}+ |Z^{\Re}|_{\mathscr{H}^p}+ |K_T^{\Re}|_{\mathscr{L}^p} \leq C_p.$$
 \end{Proposition}

\vspace{10pt}
We now  precise the results of Proposition \ref{th rep switch generic}, in the setting of this section. In particular, we describe the
optimal strategy and some of its properties that will be useful in the sequel.

\begin{Corollary}\label{co rep switch cor and dor}
(i) The following equalities hold, for all $\ell \in \set{1, \dots, d}$, $t \in [0,T]$, 
\begin{align*}
(\tYR_t)^\ell = \esssup_{a \in \mathscr{A}^{\Re}_{t,\ell}} \left(U^{\Re,a}_t{-A^{\Re,a}_t}\right) = U^{\Re,\bar{a}^{\Re}}_t {-A^{\Re,\bar{a}^{\Re}}_t} \; \text{ for some } \bar{a}^{\Re} \in \mathscr{A}^\Re_{t,\ell}\;,
\end{align*}
where $(U^{\Re,a},V^{\Re,a})$ is solution of the switched BSDE \eqref{eq de switched BSDE} 
with random driver $F(s,z) := f(s,X_s,\tYR_s,z)$ for $(s,z)\in [0,T]\times \mathcal{M}^{d,d}$, terminal condition $\xi := g(X_T)$ and costs $C_s^{ij}=c^{ij}(X_s)$. {We denote $N^{\Re,a}$ the number of switches and $A^{\Re,a}$ the cumulative cost process associated to the strategy $a$.}

\vspace{5pt}

(ii) The optimal strategy  $\bar{a}^{\Re}=(\theta_j,\alpha_j)_{j\ge 0}$  can be defined recursively by $(\theta_0,\alpha_0):=(t,\ell)$ and, for $j\ge 1$,
 \begin{align*}
 \theta_j &:= \inf \left\{  s\in [\theta_{j-1},T]\cap\Re \;\Big|\;\;\;  (\tY_s^{\Re})^{\alpha_{j-1}} \le \max_{k\neq \alpha_{j-1}} \left\{ (\tY^\Re_s)^{k} - c^{\alpha_{j-1}k}(X_s) \right\}  \right\}\,,\\
 \alpha_j &:= \min  \left\{ q \neq \alpha_{j-1} \;\Big|\;\;\;  
  (\tY^\Re_{\theta_j})^q -c^{\alpha_{j-1}q}(X_{\theta_{j}}) =
 \max_{k\neq \alpha_{j-1}} \left\{ (\tY^{\Re}_{\theta_j})^k -
 c^{\alpha_{j-1} k}(X_{\theta_{j}}) \right\} \right\}\;. 
 \end{align*}

(iii) Moreover, for all $\ell \in \set{1, \dots, d}$, $t \in [0,T]$, the optimal strategy $\bar{a}^\Re \in \mathscr{A}^\Re_{t,\ell}$ satisfies
\textcolor{black}
{
\begin{align}
\label{eq estimee U V A N}
 \EFp{t}{\sup_{s \in [t,T]} \abs{U^{\Re,\bar{a}^\Re}_s}^p + \left( \int_t^T \abs{V^{\Re,\bar{a}^\Re}_s}^2 ds \right)^{p/2} + \abs{A_T^{\Re,\bar{a}^\Re}}^p+\abs{N^{\Re,\bar{a}^\Re}}^p} \le C_p(1+|X_t|^p).
\end{align}
}
\end{Corollary}

\proof
Thanks to Proposition \ref{prop estimee Y Z K}, we can apply Proposition \ref{th rep switch generic} with the random driver $F$, the terminal condition $\xi$ and costs $C_s^{ij}$ defined above, which gives us the representation result. The first estimate in \eqref{eq estimee U V A N} is a direct application of this representation result and Proposition \ref{prop estimee Y Z K}. 
\textcolor{black}{
Other estimates in \eqref{eq estimee U V A N} are obtained by using standard arguments for BSDEs combined with the estimate \eqref{eq control YRe markov}, see proof of Proposition 2.2 in \cite{Chassagneux-Elie-Kharroubi-10} for details.
}
\eproof

\subsubsection{Fine estimates on $(Y^\Re, \tilde{Y}^\Re, Z^\Re)$}
In this section, we prove regularity results on the solution $(Y^\Re, \tilde{Y}^\Re, Z^\Re)$ of the
discretely reflected BSDEs. To do that, we will use techniques already exposed in \cite{boucha08, cha09, Chassagneux-Elie-Kharroubi-10}, based essentially on
a representation of $Z^\Re$, obtained by using Malliavin's calculus or differentiability. For a general presentation of Malliavin
Calculus, we refer to \cite{nua96}. We now introduce some notations and recall some known results on Malliavin differentiability of SDEs solution.

At the beginning of this subsection we will work under the following assumption.
\\
\HYP{r}  The coefficients $b$, $\sigma$, $g$ $f$, and $(c^{ij})_{i,j}$ are $C^{1,b}$ in all their variables, with the Lipschitz constants dominated by $L$.
\\
\textcolor{black}{This assumption is classically relieved using a kernel regularisation argument, see e.g. the proofs of Proposition 4.2 in \cite{cha09} or Proposition 3.3 in \cite{boucha08}.} { Let us emphasize that assumption \HYP{r} is purely technical. The aim of the subsection is to obtain at the end two regularity properties for $Y^{\Re}$ and $Z^{\Re}$ (see Proposition \ref{prop regularite YR} and Proposition \ref{prop regularite ZR}) where we do not need \HYP{r}. }

\vspace{10pt}

We denote by $\D$ the set of random variables $G$ which are differentiable in the Malliavin sense 
and such that $\|G\|_{\D}^2:=\|G\|^2_{\mathscr{L}^2} + \int_0^T \|D_t G\|_{\mathscr{L}^2}^2dt<\infty$, 
where $D_t G$ denotes the Malliavin derivative of $G$ at time $t\le T$. 
After possibly passing to a suitable version, an adapted process belongs to the subspace $\Lc_a^{1,2}$ of $\mathscr{H}^2$ 
whenever $V_s\in\D$ for all $s\le T$ and $\|V\|_{\cL_a^{1,2}}^2:=\|V\|^2_{\mathscr{H}^2} + \int_0^T \|D_t V\|_{\mathscr{H}^2}^2dt<\infty$.

 \begin{Remark}\label{Remark Malliavin Diff X}{\rm
 Under  \HYP{r}, the solution of \eqref{eqEDS} is Malliavin differentiable and its derivative satisfies
 \begin{align} \label{eq majo DX}
 \| \sup_{s\le T} | D_s X | \|_{_{\mathscr{S}^p}} < \infty \; \text{ and } \EFp{r}{\sup_{r \le s \le T}|D_u X_s|^p} \le C(1+|X_r|^p)\,,\quad u \le r \le T\,.
 \end{align}
 Moreover, we have
 \begin{align}\label{Control Regu DsXt}
 \sup_{s\le u } \| D_sX_t - D_s X_u \|_{_{\mathscr{L}^p}} +
 \|\sup_{t \le s \le T } | D_tX_s - D_u X_s | \; \|_{_{\mathscr{L}^p}}
 &\le C_L^p  | t-u |^{1/2} \;,\quad
 \end{align}
 for any $0\le u\le t\le T$.}\end{Remark}

\paragraph{Malliavin derivatives of $(Y^\Re, \tilde{Y}^\Re, Z^\Re)$.}

We now study the Malliavin differentiability of $(Y^\Re, \tilde{Y}^\Re, Z^\Re)$. The
techniques used are classical by now, see \cite{boucha08, cha09}. In this paragraph,
we will follow the presentation of \cite{Chassagneux-Elie-Kharroubi-10}. Once again, the main difference
with this paper is the assumption $\HYP{f}$ made on the driver $f$. In the setting of  \cite{Chassagneux-Elie-Kharroubi-10},
$f$ has to satisfy $f^i(x,y,z)=f^i(x,y^i,z^i)$ whereas $\HYP{f}$ does not impose such restriction on the $y$ variable.
This implies that the representation of $Z$, see Corollary \ref{co rep Z} below, is slightly more complicated. Namely, it  contains the term
$D \tilde{Y}$, compare to Proposition 3.2 in \cite{Chassagneux-Elie-Kharroubi-10}. To obtain the regularity results on $(Y^\Re, \tilde{Y}^\Re, Z^\Re)$,
we need thus to prove estimates on $D \tilde{Y}$, which is the main result of the next proposition.

\begin{Proposition}
 Under \HYP{f}-\HYP{r}, $(\tilde{Y}^\Re, Z^\Re)$ is Malliavin differentiable and its derivative satisfies, for all $r \in [0,T]$, $u \leq r$, $i \in \cI$, 
  \begin{align}\label{eq rep DY}
 D_u (\tY^{\Re}_r)^i &=\EFp{r}{ \px{}g^{a_T}(X_T) D_uX_T 
 + \int_r^{T} \px{}f^{a_s} (\Theta^\Re_s) D_uX_s\ud s 
 \right.
 \nonumber\\
  & \left. + \int_r^{T} \py{} f^{a_s} (\Theta^\Re_s)D_u\tY^{\Re}_s \ud s
  + \int_r^{T} \sum_{\ell=1}^d \pz{a_s \ell} f^{a_s} (\Theta^\Re_s) D_u(Z^{\Re}_s)^{a_s\ell} \ud s
 \right.
 \nonumber\\
 &  \left. -\sum_{j=1}^{N^{a}} \px{} c^{\alpha_{j-1}\alpha_j}(X_{\theta_j}) D_u X_{\theta_j}
 }
 \end{align}
 where  $a := \bar{a}^{\Re}$ is the optimal strategy associated with the representation
 in terms of switched BSDEs, recall Corollary \ref{co rep switch cor and dor}, and $\Theta^{\Re} := (X,\tilde{Y}^{\Re},Z^{\Re})$.
Moreover, following estimates hold true: for all $r \in [0,T]$, $0\leq u \leq r$, $0 \leq v \leq r$,
\begin{align}\label{eq estim DY}
 |D_u \tilde{Y}_r|^2 \le C_L (1+|X_r|^2)
 \end{align}
and
\begin{align}\label{eq estim DY-DY}
|D_u \tilde{Y}_r - D_v \tilde{Y}_r|^2 \le C_L (1+|X_r|)\EFp{r}{\sup_{r \le s \le T}|D_u X_s-D_v X_s|^4}^\frac12\;.
 \end{align}

\end{Proposition}

\proof

Let $G \in \D(\R^d)$. Since $X$ belongs to
$\cL_a^{1,2}$ under \HYP{r}, and $\cP$ is a Lipschitz continuous function,
we deduce that $\cP(X_{t},G) \in \D (\R^{d})$. Using  Lemma 5.1 in \cite{boucha08}, we compute
\begin{align}
\label{eq mal deriv cP}
& D_{s}(\cP(X_{t},G))^{i} = \\
& \sum_{j = 1}^{d}  (D_{s}G^{j}\!-\!D_{s}c_{ij}(X_{t}))\1_{G^{j} - c^{ij}(X_{t}) > \max_{\ell<j}(G^{\ell}-c^{i \ell}(X_{t}))}\1_{G^{j} - c^{ij}(X_{t}) 
\geq \max_{\ell>j}(G^{\ell}-c^{i \ell}(X_{t}))}. \nonumber
\end{align}

Combining \eqref{eq mal deriv cP}, Proposition 5.3 in \cite{ElKaroui-Peng-Quenez-97} and an induction argument, we obtain that
$(Y^{\Re},\tY^{\Re},Z^{\Re})$ is Malliavin differentiable and that a version of $(D_{u} \tY^{\Re}, D_{u}Z^{\Re})$ is given by, for all $i \in \cI$, $0 \leqslant j \leqslant \kappa-1$,  $t \in [r_j,r_{j+1})$, $0 \leq u \leq t$,

 \begin{align}\label{eq DY}
 D_u (\tY^{\Re}_t)^i &=\! D_u (Y^{\Re}_{\rjp})^i
 -  \sum_{k=1}^d  \int_t^{r_{j+1}} D_u (Z^{\Re}_s)^{ik} \ud W^k_s
 + \int_t^{\rjp} \px{}f^i (\Theta^\Re_s) D_uX_s\ud s \nonumber\\
 &\quad+ \int_t^{\rjp} \py{} f^i (\Theta^\Re_s) D_u\tY^{\Re}_s \ud s
  + \int_t^{\rjp} \sum_{\ell=1}^d \pz{i\ell} f^i (\Theta^\Re_s) D_u(Z^{\Re})^{i\ell}_s \ud s
 \,
 \end{align}
 recall \HYP{f}. 

 Now, we consider the optimal strategy $a:=\bar{a}^{\Re}$ defined in Corollary \ref{co rep switch cor and dor} (ii) above
and fix $j < \kappa$. Observing that the process $a$ is constant on the interval $[\theta_{j},\theta_{j+1})$, we deduce from  \eqref{eq DY} 
  \begin{align} \label{eq DY temp theta}
   D_u(\tY^{\Re}_r)^{\alpha_j} & =
    D_u (Y^{\Re}_{\theta_{j+1}})^{\alpha_j}
 -  \sum_{k=1}^d  \int_r^{\theta_{j+1}} D_u (Z^{\Re}_s)^{\alpha_j k} \ud W^k_s   
    +  \int_r^{\theta_{j+1}} \px{}f^{\alpha_j} (\Theta^{\Re}_s) D_uX_s \ud s 
   \\
 &\quad 
 + \int_r^{\theta_{j+1}} \py{} f^{\alpha_j} (\Theta^\Re_s)  D_u\tY^{\Re}_s \ud s \nonumber
  + \int_r^{\theta_{j+1}} \sum_{\ell=1}^d \pz{\alpha_j\ell} f^{\alpha_j} (\Theta^{\Re}_s) D_u(Z^{\Re})^{\alpha_j\ell}_s \ud s
   \end{align}
   for $ r \in [\theta_j,\theta_{j+1}]$ and $0 \leq u \leq t$.
 Combining \eqref{eq mal deriv cP} and the definition of $a$ given in Corollary \ref{co rep switch cor and dor} (ii), we compute, for $u\leq \theta_{j+1}$ and $j < \kappa$,
  \begin{align*}
 D_u (Y^{\Re}_{\theta_{j+1}})^{\alpha_j}
  \! =  D_u (\tY^{\Re}_{\theta_{j+1}})^{\alpha_{j+1}}
  - \px{} c^{\alpha_{j}\alpha_{j+1}}\!(X_{\theta_{j+1}}) D_u X_{\theta_{j+1}}.
   \end{align*}
 Inserting the previous equality into \eqref{eq DY temp theta} and summing up over $j$ 
 we obtain, for all $t\le r \le T$,
 \begin{align}\label{eq DY temp}
 D_u (\tY^{\Re}_r)^i &=\! \px{}g^{a_T}(X_T) D_uX_T 
 -    \int_r^{T} \sum_{k=1}^d D_u (Z^{\Re}_s)^{a_s k}  \ud W_s
 + \int_r^{T} \px{}f^{a_s} (\Theta^\Re_s) D_uX_s\ud s \nonumber\\
  & + \int_r^{T} \py{} f^{a_s} (\Theta^\Re_s)D_u\tY^{\Re}_s \ud s
  + \int_r^{T} \sum_{\ell=1}^d \pz{a_s \ell} f^{a_s} (\Theta^\Re_s) D_u(Z^{\Re}_s)^{a_s\ell} \ud s
 \nonumber\\
 &  -\sum_{j=1}^{N^{a}} \px{} c^{\alpha_{j-1}\alpha_j}(X_{\theta_j}) (D_u X)_{\theta_j}
 \,.
 \end{align}
Taking conditional expectation on both sides of the previous equality
proves \eqref{eq rep DY}. 
\\
Moreover, we are in the framework of section \ref{se app apriori switched rep} in the Appendix by setting
$\mathfrak{Y}=D_u Y$ and $\mathfrak{X}= D_u X$. Condition \eqref{eq app estim Na} is satisfied here by $N^{\bar{a}^\Re}$ with $\beta := C_L(1+|X|)$, recall
\eqref{eq estimee U V A N}. Using Proposition \ref{pr estim gene}  and \eqref{eq majo DX},
we then obtain \eqref{eq estim DY}.

From equation \eqref{eq DY temp}, we easily deduce the dynamics of $D_u (\tY^{\Re}) - D_v (\tY^{\Re})$, which leads, using again Proposition \ref{pr estim gene}, to \eqref{eq estim DY-DY}.
\eproof

\vspace{10pt}
The representation result for $Z^\Re$ is then an easy consequence of the previous proposition.
\begin{Corollary}\label{co rep Z}
Under  \HYP{f}-\HYP{r} the following representation holds true,
    \begin{align}\label{eq Repres Z}
 Z^{\Re}_t  = \EFp{t}{ \right. & \left. \px{}g^{a_T}(X_T) \Lambda^{a}_{t,T} D_tX_T
 -\sum_{j=1}^{N^{a}} \px{} c^{\alpha_{j-1}\alpha_j}(X_{\theta_j}) \Lambda^{a}_{t,{\theta_j}}D_t X_{\theta_j}
  \right. \nonumber \\
 & \left. 
 + \int_t^T \left(\px{} f^{a_s}(\Theta^\Re_s) \Lambda^{a}_{t,s} D_tX_s +\py{} f^{a_s} (\Theta^\Re_s) \Lambda^{a}_{t,s}D_t\tY^{\Re}_s\right)  \ud s } \;,
 \end{align}
  where  $a := \bar{a}^{\Re}$ is the optimal strategy associated with the representation
 in terms of switched BSDEs, recall Corollary \ref{co rep switch cor and dor}, and 
 for $\ell \in \cI$,
\begin{align}
\label{eq de Lambda}
\Lambda_{t,s}^{a} := \exp\left( \int_t^s \pz{a_r.} f^{a_r} (\Theta^\Re_u)\ud W_r - \frac{1}{2} \int_t^s |\pz{a_r.} f^{a_r}(\Theta^\Re_u)|^2 \ud r \right).
\end{align}

Moreover, under \HYP{f}, we have
 \begin{align}
  \label{eq co borne ZR}
  \abs{Z_t^{\Re}} &\leq \bar{L}(1+|X_t|),\quad \textrm{for all } t \in [0,T],
 \end{align}
 for some positive constant $\bar{L}$ that does not depend on the grid $\Re$.
\end{Corollary}

\proof
1. A version of $Z^\Re$ is given by $(D_t \tilde{Y}^\Re_t)_{0 \le t \le T}$. The expression
 of $D_tY$ is obtained directly by applying It\^o's formula, recall \eqref{eq rep DY}.

2. Under \HYP{f}-\HYP{r} the estimate \eqref{eq co borne ZR} follows from \eqref{eq estim DY} and \eqref{eq majo DX}. Under \HYP{f}, we can obtain the result by a standard kernel regularisation argument.
\eproof


\paragraph{Regularity of $(Y^{\Re},Z^{\Re})$.}
With the above results at hand, the study of the regularity of  $(Y^{\Re},Z^{\Re})$ follows
from ``classical'' arguments, see e.g. \cite{cha09, Chassagneux-Elie-Kharroubi-10}. For sake of completeness, we reproduce them below.

We consider a grid $\pi$ $:=$ $\{t_0=0,\ldots,t_n=T\}$ on the time interval $[0,T]$, with modulus $|\pi|:=\max_{0\leq i\leq n-1}|t_{i+1}-t_{i}|$, such that $\Re \subset \pi$.

\vspace{2mm} \noindent
We need to control the following quantities, representing the $\HP{2}$-regularity of $(\tY,Z)$:
\begin{align}
\esp{ \int_0^T|\tY^{\Re}_t - \tY^{\Re}_{\pi( t)}|^2 \ud t} \quad \text{ and }  \quad \esp{\int_0^T | Z^{\Re}_t - \bar{Z}^{\Re}_{\pi(t)} |^2  \ud t}\;,
\end{align}
 where $\pi(t):=\sup\{t_i\in\pi\,;\; t_i\le t\}$ is defined on $[0,T]$ as the projection to the closest previous grid point of $\pi$
and
\begin{align}\label{def bar Z}
\bar Z^{\Re}_{t_i} :=  \frac{1}{{t_{i+1}}-{t_i}} \esp{\int_{t_{i}}^{t_{i+1}} Z^{\Re}_s \ud s \,| \, \cF_{t_i}}\;,~~i\in\set{0,\ldots,n-1}.
\end{align}

\begin{Remark} \label{re best approx}{ \rm
Observe that $(\bar Z^{\Re}_s)_{ s\le T}:=(\bar Z^{\Re}_{\pi(s)})_{ s \le T}$ interprets as the best $\HP{2}$-approximation of the process $Z^{\Re}$ by adapted processes which are constant  on each interval $[\ti,\tip)$, for all $i<n$.
}
\end{Remark}

The first result is the regularity of the $Y$-component, which is a direct consequence of the bound \eqref{eq co borne ZR}.
\begin{Proposition}
\label{prop regularite YR}
Under $\HYP{f}$, the following holds
\begin{align*}
 \sup_{t \in [0,T]} \esp{|\tilde{Y}^{\Re}_t -\tilde{Y}^{\Re}_{\pi(t)} |^2}
 \le C_L|\pi| \;.
\end{align*}
\end{Proposition}

\proof
 We first observe that, for all $0\le t \le T$,
 \begin{align}
 \esp{ |\tY^{\Re}_t - \tY^{\Re}_{\pi(t)}|^2 } & \le     \esp{ \left| \int_{\pi(t)}^t f(X_s,\tY^{\Re}_s,Z^{\Re}_s) \ud s +  \int_{\pi(t)}^t Z^{\Re}_s \ud W_s \right |^2 } \; \nonumber
 \\
 & \le C_L \esp{\int_{\pi(t)}^t  \left(1+|X_s|^2 +|\tilde{Y}^{\Re}_s|^2 + |\tilde{Z}^{\Re}_s|^2 \right) \ud s} \,\nonumber
 \\
 &\le C_L \esp{|\pi| + \int_{\pi(t)}^t   |\tilde{Z}^{\Re}_s|^2 \ud s} \label{eq temp reg Y}
\end{align}
where we used \eqref{eq control X} and Proposition \ref{prop estimee Y Z K}. 
 From \eqref{eq co borne ZR}, we easily get
 $
  \esp{\int_{\pi(t)}^t \abs{Z_t^{\Re}}^2 \ud t } \leq C_L|\pi|
 $.
 Inserting the previous inequality into \eqref{eq temp reg Y} concludes the proof of this Proposition.
\eproof

The following Proposition gives us the regularity of $Z^{\Re}$. Its proof is postponed to the Appendix.

\begin{Proposition}
\label{prop regularite ZR}
Under $\HYP{f}$, the following holds
\begin{align*}
 \esp{\int_0^T  |Z^{\Re}_t -\bar{Z}^{\Re}_{t} |^2 \ud t}
 \le C_L\left( |\pi|^\frac12 + \kappa |\pi| \right).
\end{align*}
\end{Proposition}

 %
 %


\section{Study of the discrete-time approximation}
\label{section etude schema}

The aim of this section is to obtain a control on the error between the obliquely reflected backward scheme \eqref{schemeintro} and the discretely obliquely reflected BSDE \eqref{BSDEDORintro}. This is the purpose of Theorem \ref{th convergence scheme vers discret refl} in subsection \ref{section etude schema sous section convergence} below. In order to prove this key result, we start by interpreting the scheme in terms of the solution of a switching problem in subsection \ref{section etude schema sous section pb switching}. We then use this representation to obtain a general stability property for the scheme in subsection \ref{section etude schema sous section stabilite}. Subsection \ref{section etude schema sous section definitions} is devoted to preliminary definition and propositions.

\subsection{Definition and first estimates}
\label{section etude schema sous section definitions}

Given a grid $\pi$ of the interval $[0,T]$, we first consider an obliquely reflected backward scheme with a random generator and a random cost process $C^{\pi}$. For $t \in [0,T]$, we denote by $\mathcal{Q}_t^{\pi}$ the random closed convex set associated to $C^{\pi}_t$ and $\cP^\pi_t$ the projection onto $\cQ^\pi_t$, recall \eqref{eq def projection cout C} . The scheme is defined as follows.

\begin{Definition} $ $
\label{de scheme aleatoire general}
 \begin{enumerate}[(i)]
  \item The terminal condition $\cY^{\Re,\pi}_n$ is given by a random variable $\xi^{\pi} \in \cL^2(\mathcal{F}_T)$ valued in $\mathcal{Q}_T^{\pi}$ 
  \item for $0 \le i <n$,
  \begin{equation}
  \label{eq schema aleatoire general}
\begin{cases}
\widetilde \cY^{\Re,\pi}_{i} := \mathbb{E}[\cY^{\Re,\pi}_{i+1}\mid\mathcal{F}_{t_i
}] + h_i
F_i^{\pi}( \cZ^{\Re,\pi}_i),
\vspace*{2pt}\cr
\cZ^{\Re,\pi}_i := \mathbb{E}[\cY^{\Re,\pi}_{i+1}H_i \mid\mathcal{F}_{t_i}] ,
\vspace*{2pt}\cr
\cY^{\Re,\pi}_i := \widetilde \cY^{\Re,\pi}_i \mathbf{1}_{ \{ t_i\notin\Re\} } + \cP_{t_i}^\pi
(
\widetilde \cY^{\Re,\pi}_i)\mathbf{1}_{ \{ t_i\in\Re\} },
\end{cases}
\end{equation}
 \end{enumerate}
 with $(H_i)_{0 \le i <n}$ some $\mathbb{R}^{1 \times d}$ independent random vectors such that, for all $0 \le i < n$, $H_i$ is $\mathcal{F}_{t_{i+1}}$-measurable, $\mathbb{E}_{t_i}[H_i]=0$,
 \begin{equation}
  \label{eq structure H 1}
  \lambda_i I_{d \times d} = h_i \mathbb{E}[H_i^{\top}H_i] = h_i \mathbb{E}_{t_i}[H_i^{\top}H_i],
 \end{equation}
and 
\begin{equation}
 \label{eq structure H 2}
 \frac{\lambda}{d} \le \lambda_i \le \frac{\Lambda}{d},
\end{equation}
where $\lambda$ and $\Lambda$ are positive constants. 
\end{Definition}

\begin{Remark}
 Let us remark that \eqref{eq structure H 1} and \eqref{eq structure H 2} imply that
\begin{equation}
 \label{eq structure H 3}
 \lambda \le h_i \mathbb{E}[|H_i|^2] = h_i \mathbb{E}_{t_i} [|H_i|^2] \le \Lambda.
\end{equation}

\end{Remark}

\vspace{10pt}
In this section we use following assumptions, for some $p\geqslant 2$:

\HYP{Fd_p}
\begin{enumerate}[(i)]
 \item For all $i \in \set{0,...,n-1}$, $F_i^{\pi} : \Omega \times \cM^{d,d} \rightarrow \mathbb{R}^d$ is a $\mathcal{F}_{t_i} \otimes \mathcal{B}(\cM^{d,d})$-measurable function,
 \item the random cost process $C^\pi$ satisfies the structure condition \eqref{condition de structure C},
 \item $F_i^{\pi,j}(z) =  F_i^{\pi,j}(z^{j.})$ for all $j \in \cI$ and all $0 \le i \le n-1$,
 \item $|F_i^\pi(z)-F_i^\pi(z')| \le L^Z|z-z'|$ for all $z,z' \in \cM^{d,d}$,
 \item $\mathbb{E}\left[|\xi^\pi|^2+ \sum_{i=0}^{n-1} \abs{F_i^\pi(0)}^2 h_i+ \sup_{t_i\in \Re } |C^\pi_{t_i}|^p \right] \le C_p$, 
 \item $\sup_{0 \le i \le n-1} h_i\abs{H_i} L^Z \le 1.$
\end{enumerate}

\vspace{5pt}

\begin{Remark}
 \label{remarque existence unicite solution schema aleatoire general}
 \begin{enumerate}[i)]
  \item Under \HYP{Fd_2}, it is clear that the general scheme \eqref{eq schema aleatoire general} has a unique solution. 
  \item The weights $(H_i)_{0 \le i < n}$ depend also on the grid $\pi$ but we omit the script $\pi$ for ease of notation.
 \end{enumerate}
\end{Remark}

We observe that this obliquely reflected backward scheme can be rewritten equivalently for $i \in \llbracket 0,n \rrbracket$ as
\begin{equation}
\label{eq schema aleatoire general ecrit sans esperances conditionnelles}
\begin{cases}
 \widetilde \cY^{\Re,\pi}_{i} = \xi^\pi+\sum_{k=i}^{n-1} F_k^\pi( \cZ^{\Re,\pi}_k) h_k - \sum_{k=i}^{n-1} h_k\lambda_k^{-1} \cZ^{\Re,\pi}_k H_k^\top -\sum_{k=i}^{n-1} \Delta \cM_k+(\cK^{\Re,\pi}_n - \cK^{\Re,\pi}_i)\cr  
 \cK^{\Re,\pi}_k := \sum_{r=1}^k \Delta \cK^{\Re,\pi}_r \textrm{ with }\Delta \cK^{\Re,\pi}_r :=  \cY^{\Re,\pi}_r - \widetilde \cY^{\Re,\pi}_r,
\end{cases}
\end{equation}
with the definition $\widetilde \cY^{\Re,\pi}_{n}:=\xi^{\pi}$, with the convention $\sum_{k=n}^{n-1}...=0$,  where $(\lambda_k)$ are given by \eqref{eq structure H 1} and where, for all $k \in \llbracket 0,n-1 \rrbracket$, $\Delta \cM_k$ is an $\mathcal{F}_{t_{k+1}}$-measurable random vector satisfying 
\begin{align}\label{eq reec gen scheme}
\mathbb{E}_{t_k}[\Delta \cM_k]=0,\; \mathbb{E}_{t_k}[|\Delta \cM_k|^2] < \infty \;\text{ and }\; \mathbb{E}_{t_k}[\Delta \cM_k H_k]=0.
\end{align}

Following Corollary 2.5 in \cite{Chassagneux-Richou-14}, we know that assumption \HYP{Fd_p}(vi) is an essential ingredient to obtain a comparison result for classical time-discretized BSDE schemes. We are able to adapt this comparison result in the context of obliquely reflected backward scheme in the following proposition.

\begin{Proposition}
 \label{prop comparaison pour le schema}
 Let us consider two obliquely reflected backward schemes solutions $({}^1\widetilde{\cY}^{\Re,\pi},{}^1\cY^{\Re,\pi},{}^1\cZ^{\Re,\pi})$ and $({}^2\widetilde{\cY}^{\Re,\pi},{}^2\cY^{\Re,\pi},{}^2\cZ^{\Re,\pi})$, associated to generators $({}^1F_.^\pi)$, $({}^2F_.^\pi)$, terminal conditions ${}^1 \xi^\pi$, ${}^2 \xi^\pi$ and random cost processes $({}^1 \mathcal{C}^{\pi})$, $({}^2 C^{\pi})$ such that \HYP{Fd_2} is in force. If
 \begin{equation*}
{}^1\xi \preccurlyeq {}^2\xi, \quad \quad 
{}^1F_i({}^2\cZ_i^{\Re,\pi}) \preccurlyeq {}^2F_i({}^2\cZ_i^{\Re,\pi}), \quad \textrm{for all }~ 0 \le i \le n-1,
\end{equation*}
 \begin{equation*}
 \textrm{and} \quad ({}^1 C^{\pi}_{t_i})^{jk} \geqslant ({}^2 C^{\pi}_{t_i})^{jk}, \quad \textrm{for all }~ j,k \in \cI, \,\, t_i \in \Re,
 \end{equation*}
 then we have
 $${}^1\cY_i^{\Re,\pi} \preccurlyeq {}^2\cY_i^{\Re,\pi}\textrm{ and }~ {}^1\widetilde{\cY}_i^{\Re,\pi} \preccurlyeq {}^2 \widetilde{\cY}_i^{\Re,\pi}, \quad  \textrm{for all }~ 0 \le i \le n.$$
 Moreover, this comparison result stays true if these obliquely reflected backward schemes have two different reflection grids $\Re^1$ and $\Re^2$ with $\Re^1 \subset \Re^2$. In particular, we are allowed to have no projection for the first scheme, i.e. $\Re^1 = \emptyset$.
\end{Proposition}
\proof
We consider directly the scheme \eqref{eq schema aleatoire general}. Then, we just have to use the comparison theorem for backward schemes (Corollary 2.5 in \cite{Chassagneux-Richou-14}) and the monotonicity properties of $\cP$ (see Remark \ref{re monotonie proj}). {To be precise, in the Corollary 2.5 of \cite{Chassagneux-Richou-14} it is assumed that the inequality in assumption \HYP{Fd_p}(vi) is strict. Nevertheless, we can easily check that the result stays true when the inequality is large. }
\eproof

\begin{Proposition}
\label{prop estimee Y Z K schema general}
Assume that \HYP{Fd_2} is in force.
The unique solution $(\widetilde \cY^{\Re,\pi},\cY^{\Re,\pi},\cZ^{\Re,\pi})$ to \eqref{eq schema aleatoire general} satisfies 
$$\mathbb{E}\left[\sup_{0 \le i \le n} \abs{\widetilde \cY_i^{\Re,\pi}}^2+\sup_{0 \le i \le n} \abs{\cY_i^{\Re,\pi}}^2\right]+ \mathbb{E}\left[\sum_{i=0}^{n-1} h_i\abs{\cZ_i^{\Re,\pi}}^2\right]+ \mathbb{E}\left[ |\cK_n^{\Re,\pi}|^2\right] {+ \mathbb{E}\left[\sum_{i=0}^{n-1} \abs{\Delta \cM_i}^2\right]} \le C.$$
\end{Proposition}
\proof
The proof of uniform estimates (with respect to $n$ and $\kappa$)  divides, as usual, in two steps controlling separately $(\widetilde \cY^{\Re,\pi},\cY^{\Re,\pi})$ and $(\cZ^{\Re,\pi}, \cK^{\Re,\pi})$. It consists in transposing continuous time arguments, see e.g. proof of Theorem 2.4 in \cite{hamzha10}, in the discrete-time setting.
 
 \paragraph{Step 1. Control of $\widetilde \cY^{\Re,\pi}$ and $\cY^{\Re,\pi}$.}
We consider two non-reflected backward schemes bounding $\widetilde \cY^{\Re,\pi}$.

 Define the $\mathbb{R}^d$-valued random variable $\breve{\xi}$ and random maps $(\breve{F}_i)_{0 \le i \le n-1}$ by $(\breve{\xi})^j:=\sum_{k=1}^d \abs{(\xi^\pi)^k}$ and $(\breve{F}_i)^j(z):=\sum_{k=1}^d \abs{(F_i^\pi)^k(z)}$ for $1 \leq j \leq d$ and $0 \le i \le n-1$. We then denote by $(\breve{Y},\breve{Z})$ the unique solution of the following non-reflected backward scheme:
 \begin{equation*}
\begin{cases}
\breve{Y}_n = \breve{\xi}\cr
\breve{Z}_i = \mathbb{E}[\breve{Y}_{i+1}H_i \mid\mathcal{F}_{t_i}] , \cr
\breve{Y}_{i} = \mathbb{E}[\breve{Y}_{i+1}\mid\mathcal{F}_{t_i
}] + h_i
\breve{F}_i( \breve{Z}_i).
\end{cases}
\end{equation*}
Since all the components of $\breve{Y}$ are similar, $\breve{Y} \in \mathcal{Q}^\pi$. Thus the above backward scheme is an obliquely reflected backward scheme with same switching costs as in \eqref{eq schema aleatoire general}. We also introduce $(\mathring{Y},\mathring{Z})$ the solution of the following non-reflected backward scheme
 \begin{equation*}
\begin{cases}
\mathring{Y}_n = \xi^\pi \cr
\mathring{Z}_i = \mathbb{E}[\mathring{Y}_{i+1}H_i \mid\mathcal{F}_{t_i}] , \cr
\mathring{Y}_{i} = \mathbb{E}[\mathring{Y}_{i+1}\mid\mathcal{F}_{t_i
}] + h_i
{F}_i^\pi( \mathring{Z}_i).
\end{cases}
\end{equation*}
Using the comparison result given by Proposition \ref{prop comparaison pour le schema}, we straightforwardly deduce that $(\mathring{Y})^j \le (\widetilde \cY^{\Re,\pi})^{j} \le (\cY^{\Re,\pi})^{j} \le (\breve{Y})^j$, for all $j \in \cI$. Since $(\mathring{Y},\mathring{Z})$ and $(\breve{Y},\breve{Z})$ are solutions to standard backward schemes, {Proposition \ref{prop estimee schema non reflechi}}
and \HYP{Fd_2} lead to
\begin{align}
 \nonumber
 \mathbb{E}[\sup_{0 \le i \le n} | \widetilde \cY^{\Re,\pi}_i|^2+ \sup_{0 \le i \le n} | \cY^{\Re,\pi}_i|^2] &\le \mathbb{E}[\sup_{0 \le i \le n} | \mathring{Y}_i|^2 + \sup_{0 \le i \le n} | \breve{Y}_i|^2]\\
\label{eq estimee Ytilde dans Y Z K} & \le C\mathbb{E}\left[ \abs{\xi^\pi}^2 +\left(\sum_{i=0}^{n-1} \abs{F_i^\pi(0)}^2 h_i\right)\right]\\
\nonumber & \le C.
\end{align}

 \paragraph{Step 2. Control of $(\cZ^{\Re,\pi}, \cK^{\Re,\pi})$.}
Let us rewrite \eqref{eq schema aleatoire general ecrit sans esperances conditionnelles} for $\cY^{\Re,\pi}$ between $k$ and $k+1$ with $k \in \llbracket 0,n-1 \rrbracket$:
\begin{align*}
 \cY^{\Re,\pi}_k = \cY^{\Re,\pi}_{k+1} + F_k^\pi(\cZ_k^{\Re,\pi})h_k - h_k \lambda_k^{-1} \cZ_k^{\Re,\pi} H_k^{\top} - \Delta \cM_k + \Delta \cK_k^{\Re,\pi}.
\end{align*}
{Developing $|\cY^{\Re,\pi}_{k+1}|^2$ and taking the expectation, we have
\begin{align*}
 \mathbb{E}[|\cY^{\Re,\pi}_{k+1}|^2] =&\mathbb{E}[|\cY^{\Re,\pi}_{k}|^2] -2\mathbb{E}\left[\cY^{\Re,\pi}_{k}\left(F_k^\pi(\cZ_k^{\Re,\pi})h_k+  \Delta \cK_k^{\Re,\pi}\right)\right]\\
 & + \mathbb{E}\left[\abs{  h_k \lambda_k^{-1} \cZ_k^{\Re,\pi} H_k^\top  }^2\right]+ \mathbb{E}[|\Delta \cM_k|^2] +  \mathbb{E}\left[\abs{F_k^\pi(\cZ_k^{\Re,\pi})h_k +  \Delta \cK_k^{\Re,\pi}}^2 \right]
\end{align*}
 and, combining \HYP{Fd_2} with \eqref{eq structure H 1}-\eqref{eq structure H 2} and \eqref{eq reec gen scheme}, we get
}
\begin{align*}
 \mathbb{E}[|\cY^{\Re,\pi}_{k+1}|^2]
 &\geqslant \mathbb{E} [|\cY^{\Re,\pi}_{k}|^2]  -C\mathbb{E}\left[|\cY^{\Re,\pi}_{k}|\left(|F_k^\pi(0)|h_k + |\cZ_k^{\Re,\pi}|h_k\right)\right] -2\mathbb{E}\left[\cY^{\Re,\pi}_{k}\Delta \cK_k^{\Re,\pi}\right]\\
 &\quad + \mathbb{E}\left[h_k^2\lambda_k^{-2}\mathbb{E}_{t_k}\bigg[ \sum_{i,j\in \llbracket 1,d \rrbracket} ((\cZ_k^{\Re,\pi})^{\top}\cZ_k^{\Re,\pi})^{ij}(H_k)^{1i}(H_k)^{1j}  \bigg]\right]+ \mathbb{E}[|\Delta \cM_k|^2]\\
 & \geqslant \mathbb{E} [|\cY^{\Re,\pi}_{k}|^2]  -C\mathbb{E}\left[|\cY^{\Re,\pi}_{k}|\left(|F_k^\pi(0)|h_k + |\cZ_k^{\Re,\pi}|h_k\right)\right]-2\mathbb{E}\left[\cY^{\Re,\pi}_{k}\Delta \cK_k^{\Re,\pi}\right]\\
 & \quad + \frac{d}{\Lambda}\mathbb{E}\left[h_k|\cZ_k^{\Re,\pi}|^2\right]+ \mathbb{E}[|\Delta \cM_k|^2].
\end{align*}
Then we sum over $k \in \llbracket 0, n-1 \rrbracket$ and we compute, using Young's inequality with $\varepsilon>0$,
\begin{align}
  \nonumber \sum_{k=0}^{n-1} \mathbb{E}\left[h_k|\cZ_k^{\Re,\pi}|^2\right] +  \sum_{k=0}^{n-1} \mathbb{E}[|\Delta \cM_k|^2] &\le C_{\varepsilon}\mathbb{E}\left[ \sup_{0 \le k \le n} | \cY^{\Re,\pi}_k|^2 + \sum_{k=0}^{n-1} |F_k^\pi(0)|^2h_k\right]\\
  \nonumber &\quad +\varepsilon \sum_{k=0}^{n-1} \mathbb{E}\left[h_k|\cZ_k^{\Re,\pi}|^2\right]+ 2\mathbb{E}\left[ \sup_{0 \le k \le n} | \cY^{\Re,\pi}_k| |\cK_n^{\Re,\pi}|\right]\\
  \nonumber
  &\le C_{\varepsilon}\mathbb{E}\left[ \sup_{0 \le k \le n} | \cY^{\Re,\pi}_k|^2 + \sum_{k=0}^{n-1} |F_k^\pi(0)|^2h_k\right]\\
  \label{eq estimee Z} &\quad +\varepsilon \sum_{k=0}^{n-1} \mathbb{E}\left[h_k|\cZ_k^{\Re,\pi}|^2\right]+ \varepsilon\mathbb{E}\left[   |\cK_n^{\Re,\pi}|^2\right].
\end{align}
Moreover, we get from \eqref{eq schema aleatoire general ecrit sans esperances conditionnelles}
\begin{align}
 \nonumber
 \mathbb{E}\left[ |\cK_n^{\Re,\pi}|^2\right] &\le C\mathbb{E}\left[ \sup_{0 \le k \le n} | \widetilde \cY^{\Re,\pi}_k|^2 + \sum_{k=0}^{n-1} |F_k^\pi(0)|^2h_k\right]\\
\label{eq estimee K}
 & \quad +C\sum_{k=0}^{n-1} \mathbb{E}\left[h_k|\cZ_k^{\Re,\pi}|^2\right] +  C\sum_{k=0}^{n-1} \mathbb{E}[|\Delta \cM_k|^2].
\end{align}
Combining \eqref{eq estimee Z} with \eqref{eq estimee K}, and using \HYP{Fd_2} and \eqref{eq estimee Ytilde dans Y Z K}, classical calculations yield, {for $\varepsilon<(1+C)^{-1}$ with $C$ the constant appearing in \eqref{eq estimee K}},
\begin{align*}
   \sum_{k=0}^{n-1} \mathbb{E}\left[h_k|\cZ_k^{\Re,\pi}|^2\right] +  \sum_{k=0}^{n-1} \mathbb{E}[|\Delta \cM_k|^2] &\le C.
\end{align*}
Finally we can insert this last inequality into \eqref{eq estimee K} and use once again \HYP{Fd_2} and \eqref{eq estimee Ytilde dans Y Z K} to conclude the proof.
\eproof

\subsection{Optimal switching problem representation}
\label{section etude schema sous section pb switching}

We now introduce a discrete-time version of the switching problem, which will allow us to give a new representation of the scheme given in Definition \ref{de scheme aleatoire general}. To simplify notations, we start by adapting the definition of switching strategies to the discrete-time setting. A switching strategy $a$ is now a nondecreasing sequence of stopping times $(\theta_r)_{r \in \mathbb{N}}$ valued in $\mathbb{N}$, combined with a sequence of random variables $(\alpha_r)_{r \in \mathbb{N}}$ valued in $\cI$, such that  $\alpha_r$ is $\mathcal{F}_{t_{\theta_r}}$-measurable for any $r \in \mathbb{N}$. 

Then by mimicking Section \ref{section representation property}, we define classical objects related to switching strategies. For a switching strategy $a= (\theta_{r},\alpha_{r})_{r\in\N}$, we introduce $\cN^{a}$ the (random) number of switches before $n$:
\begin{align}
\cN^{a}  =  \#\{r\in\N^*~:~\theta_{r} \leq n \}\;. \label{eq def NaRpi}
\end{align}
To any switching strategy $a=(\theta_{r},\alpha_{r})_{r\in\N}$, we
associate the current state process $(a_i)_{i\in \llbracket 0,n\rrbracket}$ and the cumulative cost process $(\cA^{a}_i)_{i \in \llbracket0,n \rrbracket}$ defined respectively by
  \begin{align*}
  a_i \;:=\; \alpha_{0}\mathbbm{1}_{\{0\leq i<\theta_0\}}+\sum_{r=1}^{\cN^{a}}  \alpha_{r-1} \mathbbm{1}_{\{\theta_{r-1}\leq i < \theta_{r}\}}  \;\text{ and }\;  
  \cA^{a}_i \;:=\; \sum_{r=1}^{\cN^{a}} (C_{t_{\theta_{r}}}^\pi)^{\alpha_{r-1}\alpha_{r}} \mathbbm{1}_{\{\theta_r\le i \le n\}}\,, 
  \end{align*}
for $0\le i\le n$.
We denote by $\mathscr{A}^{\Re,\pi}$ the set of $\Re$-admissible strategies:
\begin{align*}
\mathscr{A}^{\Re,\pi} = \set{a=(\theta_{r},\alpha_{r})_{r \in \mathbb{N}} \textrm{ switching strategy } ~|~  t_{\theta_r} \in \Re \quad \forall r \in \llbracket{1,\cN^{a}\rrbracket} ,~\esp{|\mathcal{A}_{n}^{a}|^2}< \infty}\;.
\end{align*}
For $(i,j)\in\llbracket0,n\rrbracket \times\cI$, the set $\mathscr{A}^{\Re,\pi}_{i,j}$ of admissible strategies starting from $j$ at time $t_i$ is defined by
\begin{align*}
\mathscr{A}^{\Re,\pi}_{i,j} = \set{a=(\theta_{r},\alpha_{r})_{r \in \mathbb{N}}\in \mathscr{A}^{\Re,\pi} ~| \theta_{0}= i,~\alpha_{0}= j}\;.
\end{align*}

For a strategy $a \in \mathscr{A}^{\Re,\pi}_{i,j}$ we define the one dimensional $\Re$-switched backward scheme whose solution $(\cU^{\Re,\pi,a},\cV^{\Re,\pi,a})$ satisfies
\begin{align}
 \label{eq switched backward scheme}
\begin{cases}
 \cU^{\Re,\pi,a}_{n} = \xi^{\pi,a_{n}} \cr 
 \cV^{\Re,\pi,a}_{k} = \mathbb{E}[\cU^{\Re,\pi,a}_{k+1}H_k^\top \mid\mathcal{F}_{t_k}] ,\cr
\cU^{\Re,\pi,a}_{k} = \mathbb{E}[\cU^{\Re,\pi,a}_{k+1}\mid\mathcal{F}_{t_k
}] + h_k
F^{\pi,a_k}_k(\cV^{\Re,\pi,a}_{k})- \sum_{j=1}^{\cN^{a}} (C^{\pi}_{t_{\theta_j}})^{\alpha_{j-1}\alpha_{j}}\mathbbm{1}_{\theta_j\le k}, \quad i \le k <n.
\end{cases}
\end{align}

Similarly to equation \eqref{eq schema aleatoire general}, we observe that this obliquely reflected backward scheme can be rewritten equivalently for $k \in \llbracket i,n \rrbracket$ as
\begin{align}
 \nonumber
 \cU^{\Re,\pi,a}_{k} = &\xi^{\pi,a_n}+\sum_{m=k}^{n-1} F_m^{\pi,a_m}( \cV^{\Re,\pi,a}_m) h_m - \sum_{m=k}^{n-1} h_m\lambda_m^{-1} \cV^{\Re,\pi,a}_m H_m^\top -\sum_{m=k}^{n-1} \Delta \cM^a_m\\
 \label{eq switched backward scheme sans esp cond}
 &-\cA^{a}_n + \cA^{a}_k  
\end{align}
where $(\lambda_k)$ are given by \eqref{eq structure H 1} and, for all $k \in \llbracket 0,n-1 \rrbracket$, $\Delta \cM_k^a$ is an $\mathcal{F}_{t_{k+1}}$-measurable random variable satisfying 
\begin{align}\label{eq reec disc switch scheme}
\mathbb{E}_{t_k}[\Delta \cM_k^a]=0,\; \mathbb{E}_{t_k}[|\Delta \cM_k^a|^2] < \infty\; \text{ and } \mathbb{E}_{t_k}[\Delta \cM_k^a H_k]=0.
\end{align}

The next theorem is a Snell envelope representation of the obliquely reflected backward scheme.
\begin{Proposition}
\label{prop snell representation schema}
Assume than \HYP{Fd_2} is in force.
For any $j \in \cI$ and $0 \le i \le n$, the following hold:
\begin{enumerate}[(i)]
 \item The discrete process $\cY^{\Re,\pi}$ dominates any $\Re$-switched backward scheme {excluding a possible instantaneous initial switch}, that is,
 \begin{equation}
 \label{eq 1 prop snell enveloppe rep discret} 
 \cU^{\Re,\pi,a}_i {-\mathcal{A}_i^{a}} \le  ({\cY^{\Re,\pi}_i})^j, \quad \mathbb{P}\textrm{-a.s. for any } a \in \mathscr{A}^{\Re,\pi}_{i,j}.
 \end{equation}
 \item Define the strategy $\bar{a}^{\Re,\pi}=(\bar{\theta}_r,\bar{\alpha}_r)_{r \geqslant 0}$ recursively by $(\bar{\theta}_0,\bar{\alpha}_0):=(i,j)$ and, for $r \geqslant 1$,
 \begin{align*}
 \bar{\theta}_r & := \inf \big\{ k \in \llbracket\bar{\theta}_{r-1},n\rrbracket  ~\big{|} ~ t_k \in \Re, ~ (\widetilde{\cY}^{\Re,\pi}_k)^{\bar{\alpha}_{r-1}} \leqslant \max_{m \neq \bar{\alpha}_{r-1}} \{(\widetilde{\cY}^{\Re,\pi}_{k})^m -C_{t_k}^{\bar{\alpha}_{r-1} m} \} \big\},\\
  \bar{\alpha}_r & := \min \big\{ \ell \neq \bar{\alpha}_{r-1} ~\big|~ (\widetilde{\cY}^{\Re,\pi}_{\theta_r})^\ell -C_{t_{\theta_r}}^{\bar{\alpha}_{r-1} \ell} = \max_{m \neq \bar{\alpha}_{r-1}} \{(\widetilde{\cY}^{\Re,\pi}_{\theta_r})^m -C_{t_{\theta_r}}^{\bar{\alpha}_{r-1} m} \}  \big\}
 \end{align*}
 Then we have $\bar{a}^{\Re,\pi} \in \mathscr{A}^{\Re,\pi}_{i,j}$ and
 \begin{equation}
 \label{eq 2 prop snell enveloppe rep discret}
 ({\cY^{\Re,\pi}_i})^j = \cU^{\Re,\pi,\bar{a}^{\Re,\pi}}_i {-\mathcal{A}_i^{\bar{a}^{\Re,\pi}}} \quad \mathbb{P}\textrm{-a.s.}
 \end{equation}
 \item The following ``Snell envelope'' representation holds:
 \begin{equation}
  \label{eq 3 prop snell enveloppe rep discret}
({\cY^{\Re,\pi}_i})^j = \esssup_{a \in \mathscr{A}^{\Re,\pi}_{i,j}} \left( \cU^{\Re,\pi,a}_i  {-\mathcal{A}_i^{{a}}} \right) \quad \mathbb{P}\textrm{-a.s.}
\end{equation}
\end{enumerate}
\end{Proposition}

\proof
{We will adapt the proof of Theorem 3.1 in \cite{huytan10} and Theorem 2.1 in \cite{Chassagneux-Elie-Kharroubi-10} to the discrete time setting. As we have already remark before Proposition \ref{th rep switch generic}, there are some mistakes in these two results that we must take into account to obtain the correct formulation.} Observe first that assertion (iii) is a direct consequence of (i) and (ii). 

Let us fix $i \in \llbracket 0,n \rrbracket$ and $j \in \cI$.

\paragraph{Step 1.}We first prove (i). 
{ Set $a =(\theta_r,\alpha_r)_{r \geqslant 0} \in \mathscr{A}_{i,j}^{\Re,\pi}$. Since the sequence $(\theta_r)_{r \geqslant 0}$ is just a non decreasing sequence, we introduce a subsequence $a^\#$ that is increasing. We recursively define $(\theta_r^\#,\alpha_r^\#)_{r \geqslant 0}$ by
\begin{align*}
 &\theta^\#_0:=\theta_{\underline{r}_0}\quad \text{with} \quad \underline{r}_0:=0\\
 &\alpha_0^{\#} := \alpha_{\overline{r}_0}\quad \text{with} \quad \overline{r}_0:=\sup\{k \in \mathbb{N} | \theta_k=\theta^\#_0\}\\
 &\theta^\#_{r+1}:=\theta_{\underline{r}_{r+1}} \quad \text{with} \quad \underline{r}_{r+1}:=\inf \{k\in \N | \theta_k>\theta^\#_r\} \\
 &\alpha^\#_{r+1}:=\alpha_{\overline{r}_{r+1}} \quad \text{with} \quad \overline{r}_{r+1}:= \sup\{k \in \mathbb{N} | \theta_k = \theta^\#_{r+1}\}.
\end{align*}
Let us observe that this new strategy $a^\#$ only keeps track of the last state when instantaneous switches occur in the strategy $a$.
We also introduce the process $(\widetilde{\cY}^a,\cZ^a)$ defined, for $k \in \llbracket i,n \rrbracket$, by
\begin{equation}
\label{def Ya Za}
\begin{cases}
 \widetilde{\cY}^a_k &:= \sum_{r \geqslant 0} (\widetilde \cY^{\Re,\pi}_k)^{\alpha_r} \mathbbm{1}_{\theta_r \leqslant k < \theta_{r+1}} + \xi^{\pi,a_n} \mathbbm{1}_{k=n}\\
 \cZ^a_k &:= \sum_{r\geqslant 0} (\cZ_k^{\Re,\pi})^{\alpha_r} \mathbbm{1}_{\theta_r\leqslant k < \theta_{r+1}},
\end{cases}
 \end{equation}}
 {
and we can remark that $(\widetilde{\cY}^a,\cZ^a)=(\widetilde{\cY}^{a^\#},\cZ^{a^\#})$.
Observe that these processes jump between the components of the obliquely reflected backward scheme \eqref{eq schema aleatoire general ecrit sans esperances conditionnelles} according to the strategy $a$, and, between two jumps of $a^\#$, we have
\begin{align}
 \nonumber
 \widetilde{\cY}^a_{\theta_r^\#} &= ({\cY}^{\Re,\pi}_{\theta_{r+1}^\#})^{\alpha_r^\#}+ \sum_{k=\theta_r^\#}^{\theta_{r+1}^\#-1} F_k^{\pi,\alpha_r^\#}( (\cZ^{\Re,\pi}_k)^{\alpha_r^\#}) h_k - \sum_{k=\theta_r^\#}^{\theta_{r+1}^\#-1} h_k\lambda_k^{-1} (\cZ^{\Re,\pi}_k)^{\alpha_r^\#} H_k^\top \\
 \nonumber
 &\quad \quad-\sum_{k=\theta_r^\#}^{\theta_{r+1}^\#-1} \Delta (\cM_k)^{\alpha_r^\#}+(\cK^{\Re,\pi}_{\theta_{r+1}^\#-1})^{\alpha_r^\#} - (\cK^{\Re,\pi}_{\theta_r^\#})^{\alpha_r^\#}\\
\nonumber
 &=\widetilde{\cY}^a_{\theta_{r+1}^\#}+\sum_{k=\theta_r^\#}^{\theta_{r+1}^\#-1} F_k^{\pi,a_k}( \cZ^{a}_k) h_k - \sum_{k=\theta_r^\#}^{\theta_{r+1}^\#-1} h_k\lambda_k^{-1} \cZ^{a}_k H_k^\top -\sum_{k=\theta_r^\#}^{\theta_{r+1}^\#-1} \Delta \cM_k^{a} \\
\label{eq Ya entre deux switch}
 &\quad \quad+ (\cK^{\Re,\pi}_{\theta_{r+1}^\#-1})^{\alpha_r^\#} - (\cK^{\Re,\pi}_{\theta_r^\#})^{\alpha_r^\#}+ \left(({\cY}^{\Re,\pi}_{\theta_{r+1}^\#})^{\alpha_r^\#}-(\widetilde{\cY}^{\Re,\pi}_{\theta_{r+1}^\#})^{\alpha_{r+1}^\#}\right), \quad \quad r \geqslant 0.
\end{align}
Introducing
\begin{align*}
 \cK^a_k :=& \sum_{r=0}^{\cN^{a^{\#}}-1} \left[ \sum_{m=(\theta_r^\#+1)\wedge k}^{(\theta_{r+1}^\#-1)\wedge k} (\Delta \cK^{\Re,\pi}_m)^{\alpha_r^\#}\right.\\
 &\quad \quad + \mathbbm{1}_{\theta_{r+1}^\# \le k} \left( ({\cY}^{\Re,\pi}_{\theta_{r+1}^\#})^{\alpha_r^\#}-(\widetilde{\cY}^{\Re,\pi}_{\theta_{r+1}^\#})^{\alpha_{r+1}^\#}+(C_{t_{\theta_{r+1}^\#}}^\pi)^{\alpha_{r}^\#\alpha_{r+1}^\#}\right)\\
 &\left.\quad \quad +\mathbbm{1}_{\theta_{r+1}^\# \le k} \left( \sum_{m=\underline{r}_{r+1}}^{\overline{r}_{r+1}} (C_{t_{\theta_{r+1}^\#}}^\pi)^{\alpha_{m-1} \alpha_{m}}  -(C_{t_{\theta_{r+1}^\#}}^\pi)^{\alpha_{r}^\#\alpha_{r+1}^\#}\right)\right]
\end{align*}
for $k \in \llbracket i,n \rrbracket$, and summing up \eqref{eq Ya entre deux switch} over $r$, we get, for $k \in \llbracket i,n\rrbracket$,
\begin{align*}
 \widetilde{\cY}^a_k =& \xi^{\pi,a_n}+\sum_{m=k}^{n-1} F_m^{\pi,a_m}( \cZ^{a}_m) h_m - \sum_{m=k}^{n-1} h_m\lambda_m^{-1} \cZ^{a}_m H_m^\top -\sum_{m=k}^{n-1} \Delta (\cM_m)^{a}\\
 & -\cA^{a}_n+\cA^{a}_k + \cK^a_n-\cK^a_k.
\end{align*}}
{
Using the relation $\cY^{\Re,\pi}_{\theta_r^\#} = \cP_{t_{\theta_r^\#}}^\pi (\widetilde{\cY}^{\Re,\pi}_{\theta_r^\#})$ for all $r \in \llbracket 1,\cN^{a^\#} \rrbracket$ and the structure condition \ref{condition de structure C}, we easily check that $\cK^a$ is an increasing process. Since $\cU^{\Re,\pi,a}$ solves \eqref{eq switched backward scheme sans esp cond}, we deduce by a comparison argument (see Corollary 2.5 in \cite{Chassagneux-Richou-14}) that $\cU^{\Re,\pi,a}_i \le \widetilde{\cY}^a_i$. If $t_i \notin \Re$, then $\widetilde{\cY}^a_i = (\widetilde{\cY}^{\Re,\pi}_{i})^{j} = ({\cY}^{\Re,\pi}_{i})^{j}$ and $\mathcal{A}_{i}^a=0$ which implies that 
$$\cU^{\Re,\pi,a}_i {-\mathcal{A}_i^{a}} \le  ({\cY}^{\Re,\pi}_i)^j.$$
When $t_i \in \Re$,  we can remark that $ \widetilde{\cY}^a_i=(\widetilde{\cY}^{\Re,\pi}_{i})^{\alpha_0^\#}$ and so we have
\begin{align} \nonumber
 ({\cY}^{\Re,\pi}_{i})^{j} +\mathcal{A}_i^a =&  (\widetilde{\cY}^{\Re,\pi}_{i})^{\alpha_0^\#} + \left( ({\cY}^{\Re,\pi}_{i})^{j} - (\widetilde{\cY}^{\Re,\pi}_{i})^{\alpha_0^\#} + (C^{\pi}_{t_i})^{j \alpha^\#_0}\right) + \left( \sum_{m=1}^{\overline{r}_0} (C_{t_i}^{\pi})^{\alpha_{m-1} \alpha_m} - (C^{\pi}_{t_i})^{j \alpha^\#_0}\right)\\
\label{ineg:point initial}
 \geqslant& (\widetilde{\cY}^{\Re,\pi}_{i})^{\alpha_0^\#}  \geqslant \cU^{\Re,\pi,a}_i,
\end{align}
by using the relation $\cY^{\Re,\pi}_{i} = \cP_{t_i}^\pi (\widetilde{\cY}^{\Re,\pi}_{i})$ and the structure condition \ref{condition de structure C}.
Since $a$ is arbitrary in $\mathscr{A}^{\Re,\pi}_{i,j}$, we deduce \eqref{eq 1 prop snell enveloppe rep discret}.}

\paragraph{Step 2.} We now prove (ii).

Consider the strategy $\bar{a}^{\Re,\pi}$ given above as well as the associated process $(\widetilde{\cY}^{\bar{a}^{\Re,\pi}}, \cZ^{\bar{a}^{\Re,\pi}})$ defined as in \eqref{def Ya Za}. By definition of $ \bar{a}^{\Re,\pi}$, we have
{
\begin{equation}
\label{eq switch discret 1}
(\cY^{\Re,\pi}_{\bar{\theta}_{r+1}^\#})^{\bar{\alpha}_r^\#} = (\cP_{t_{\bar{\theta}_{r+1}^\#}}^\pi (\widetilde{\cY}^{\Re,\pi}_{\bar{\theta}_{r+1}^\#}))^{\bar{\alpha}_r^\#} = (\widetilde{\cY}^{\Re,\pi}_{\bar{\theta}_{r+1}^\#})^{\bar{\alpha}_{r+1}^\#}-(C_{t_{\bar{\alpha}_{r+1}^\#}}^\pi)^{\bar{\alpha}_r^\#\bar{\alpha}_{r+1}^\#}, \quad r \geqslant 0
\end{equation}
and
$$(\cY^{\Re,\pi}_{k})^{\bar{\alpha}_r^\#} > \max_{m \neq \bar{\alpha}_{r}^\#} \{(\widetilde{\cY}^{\Re,\pi}_{k})^m -C_{t_k}^{\bar{\alpha}_{r}^\# m} \},  \quad r \geqslant 0,\, k \in \rrbracket \bar{\theta}_{r}^\#, \bar{\theta}_{r+1}^\# \llbracket, \, t_k \in \Re,$$}
{which imply that
\begin{equation}
\label{eq switch discret 2}
(\cY^{\Re,\pi}_{k})^{\bar{\alpha}_r} = (\widetilde{\cY}^{\Re,\pi}_k)^{\bar{\alpha}_r}, \quad r \geqslant 0, \,  k \in \rrbracket \bar{\theta}_{r}, \bar{\theta}_{r+1} \llbracket.
\end{equation}
Moreover, the structure condition \eqref{condition de structure C} implies that the optimal switching strategy is allowed to switch at most one time per date and so we have, for this strategy,
\begin{equation}
\label{eq switch discret 3}
\overline{r}_{r+1} = \underline{r}_{r+1}, \quad \forall r \ge 0.
\end{equation}
We can remark that equations \eqref{eq switch discret 1}, \eqref{eq switch discret 2} and \eqref{eq switch discret 3} give us that
 $\cK^{\bar{a}^{\Re,\pi}}=0$ and then, for all  $k \in \llbracket i,n\rrbracket$,
\begin{align}
 \nonumber
 \widetilde{\cY}^{\bar{a}^{\Re,\pi}}_k =& \xi^{\pi,{\bar{a}^{\Re,\pi}}_n}+\sum_{m=k}^{n-1} F_m^{\pi,{\bar{a}^{\Re,\pi}}_m}( \cZ^{{\bar{a}^{\Re,\pi}}}_m) h_m - \sum_{m=k}^{n-1} h_m\lambda_m^{-1} \cZ^{{\bar{a}^{\Re,\pi}}}_m H_m^\top\\ 
 \label{eq etape 2 snell schema}
 & -\sum_{m= k}^{n-1} (\Delta \cM_m)^{a} -\cA^{{\bar{a}^{\Re,\pi}}}_n+\cA^{{\bar{a}^{\Re,\pi}}}_k.
\end{align}
Hence, $(\widetilde{\cY}^{\bar{a}^{\Re,\pi}},\cZ^{\bar{a}^{\Re,\pi}})$ and $(\cU^{\Re,\pi,\bar{a}^{\Re,\pi}},\cV^{\Re,\pi,\bar{a}^{\Re,\pi}})$ are solutions of the same backward scheme and $(\widetilde{\cY}_i^{\Re,\pi})^{\bar{\alpha}_0^\#} = \cU^{\Re,\pi,\bar{a}^{\Re,\pi}}_i$. If $t_i \notin \Re$, then we get directly \eqref{eq 2 prop snell enveloppe rep discret}. When $t_i \notin \Re$, we use the definition of $\bar{a}^{\Re,\pi}$ to obtain that the inequality in \eqref{ineg:point initial} becomes an equality for this optimal strategy.}
{To complete the proof, we only need to check that ${\bar{a}^{\Re,\pi}} \in \mathscr{A}^{\Re,\pi}$, that is $\mathbb{E}[|\cA^{\bar{a}^{\Re,\pi}}_{n}|^2] <\infty$. By definition of ${\bar{a}^{\Re,\pi}}$ on $\llbracket i, n \rrbracket$ and the structure condition on costs \eqref{condition de structure C}, we get $|\cA^{\bar{a}^{\Re,\pi}}_i| \le \max_{k \neq j} |C^{jk}_{t_i}|$ since we have at most one jump at date $i$. Then, Assumption \HYP{Fd 2} gives us $\mathbb{E}[|\cA^{\bar{a}^{\Re,\pi}}_i|^2] \le C$.  Combining \eqref{eq etape 2 snell schema}, for $k=i$, with the Lipschitz property of $F^\pi$ we get that
\begin{align*}
 \mathbb{E}\left[|\cA^{\bar{a}^{\Re,\pi}}_n|^2 \right] \leqslant& C\left(1+\mathbb{E}\left[\sup_{i \leqslant i \leqslant n}  |\widetilde{\cY}^{\bar{a}^{\Re,\pi}}_k|^2\right] + \mathbb{E}\left[\sum_{k=i}^{n-1} h_k\abs{\cZ_k^{\Re,\pi}}^2\right]+\mathbb{E}[|\cA^{\bar{a}^{\Re,\pi}}_i|^2]\right)\\
 &+ Cn\mathbb{E}\left[ \sum_{k=i}^{n-1} |\Delta \cM_k|^2\right].
\end{align*}
Finally, we just have to apply estimates in Proposition \ref{prop estimee Y Z K schema general} to obtain the square integrability of $\cA^{\bar{a}^{\Re,\pi}}_n$ and the proof is complete.}
\eproof

\begin{Proposition}
\label{prop estimee U V N schema}
Assume than \HYP{Fd_2} is in force. For all $0 \le i \le n$, $j \in \cI$, we have 
$$ \mathbb{E}\left[\sup_{i \le k \le n} \abs{\cU^{\Re,\pi,\bar{a}^{\Re,\pi}}_k}^2 +\sum_{k=i}^{n-1} h_k\abs{\cV^{\Re,\pi,\bar{a}^{\Re,\pi}}_k}^2 +\abs{\cA^{\bar{a}^{\Re,\pi}}_n}^2 +\abs{\cN^{\bar{a}^{\Re,\pi}}}^2\right] \le C,$$
for the optimal strategy $\bar{a}^{\Re,\pi} \in \mathscr{A}^\Re_{i,j}$.
\end{Proposition}
\proof
Fix $(i,j) \in \llbracket 0, n \rrbracket \times \cI$. According to the identification of $(\cU^{\Re,\pi,\bar{a}^{\Re,\pi}},\cV^{\Re,\pi,\bar{a}^{\Re,\pi}})$ with $(\widetilde{\cY}^{\bar{a}^{\Re,\pi}},\cZ^{\bar{a}^{\Re,\pi}})$ obtained in the proof of Proposition \ref{prop snell representation schema}, we deduce from Proposition \ref{prop estimee Y Z K schema general} expected controls on $\cU^{\Re,\pi,\bar{a}^{\Re,\pi}}$ and $\cV^{\Re,\pi,\bar{a}^{\Re,\pi}}$. 

By taking conditional expectation in \eqref{eq etape 2 snell schema}, we have
\begin{align*}
 \mathbb{E}_{t_i} [\cA^{{\bar{a}^{\Re,\pi}}}_n] =&  \mathbb{E}_{t_i} \left[\xi^{\pi,{\bar{a}^{\Re,\pi}}_n}  -\widetilde{\cY}^{\bar{a}^{\Re,\pi}}_i + \sum_{m=i}^{n-1} F_m^{\pi,{\bar{a}^{\Re,\pi}}_m}( \cZ^{{\bar{a}^{\Re,\pi}}}_m) h_m +\cA^{{\bar{a}^{\Re,\pi}}}_i\right].
\end{align*}
Thus, using standard inequalities and the growth of $F^\pi$, we easily obtain
\begin{align*}
 \mathbb{E}[|\cA^{{\bar{a}^{\Re,\pi}}}_n|^2] \le&  C\mathbb{E} \left[\sup_{i \le k \le n} |\widetilde{\cY}^{\bar{a}^{\Re,\pi}}_k|^2 + \sum_{m=i}^{n-1} |\cZ^{\bar{a}^{\Re,\pi}}_m|^2 h_m +|\cA^{{\bar{a}^{\Re,\pi}}}_i|^2\right].
\end{align*}

We have already noticed in the proof of Proposition \ref{prop snell representation schema} that we have $|\cA^{\bar{a}^{\Re,\pi}}_i| \le \max_{k \neq j} |C^{jk}_i|$, which inserted into the previous inequality leads to $\mathbb{E}[|\cA^{{\bar{a}^{\Re,\pi}}}_n|^2] \le C$.

We finally complete the proof, observing from the structure condition \eqref{condition de structure C} that 
$$\mathbb{E}[|\cN^{{\bar{a}^{\Re,\pi}}}|^2] \le C\mathbb{E}[|\cA^{{\bar{a}^{\Re,\pi}}}_n|^2].$$
\eproof

\subsection{Stability of obliquely reflected backward schemes}
\label{section etude schema sous section stabilite}

We now consider two obliquely reflected backward schemes, with different parameters but the same reflection grid $\Re$. For $\ell \in \{1,2\}$, we consider an $\mathcal{F}_T$-measurable random terminal condition ${}^{\ell} \xi$, a random generator $z \mapsto { }^{\ell}F(.,z)$ and random cost processes $(^{\ell}C^{ij})_{1 \le i ,j \le d}$ satisfying the structural condition \eqref{condition de structure C}. As in Subsection \ref{section etude schema sous section pb switching}, terminal conditions, generators and cost processes are allowed to depend on $\pi$ but we omit the script $\pi$ for reading convenience. We denote by $({}^{\ell} \widetilde \cY^{\Re,\pi},  {}^{\ell}  \cY^{\Re,\pi}, {}^{\ell} \cZ^{\Re,\pi})$ the solution of the associated obliquely reflected backward scheme.

Defining $\delta \cY^{\Re,\pi} := {}^{1} \cY^{\Re,\pi}-{}^{2} \cY^{\Re,\pi}$, $\delta \widetilde \cY^{\Re,\pi} := {}^{1} \widetilde \cY^{\Re,\pi}-{}^{2} \widetilde \cY^{\Re,\pi}$, $\delta \cZ^{\Re,\pi} := {}^{1} \cZ^{\Re,\pi}-{}^{2} \cZ^{\Re,\pi}$, $\delta \xi := {}^{1} \xi- {}^{2} \xi$ together with 
\begin{align*}
\abs{\delta C_t}_{\infty} &:= \max_{i,j \in \cI} \abs{{}^1 C_t^{ij} - {}^2 C_t^{ij}},\\
\abs{\delta F_k}_{\infty} &:= \max_{i \in \cI} \sup_{z \in \cM^{d,d}} \abs{{}^1 F^i_k-{}^2 F^i_k}(z),
\end{align*}
for $0 \le k \le n-1$, we prove the following stability result.

\begin{Proposition}
\label{prop stabilite schema gene aleatoire}
Assume that \HYP{Fd_p} is in force for some given $p \geqslant 2$. Then we have, for any $i \in \llbracket 0,n \rrbracket$,
 \begin{align*}
  &\sup_{i \le k \le n} \mathbb{E}\left[\abs{\delta \cY^{\Re,\pi}_k}^2 + \abs{ \delta \widetilde \cY^{\Re,\pi}_k}^2 \right] +\frac{1}{\kappa}\mathbb{E}\left[\sum_{k=i}^{n-1} h_k \abs{\delta \cZ^{\Re,\pi}_k}^2\right] \\
  & \leq C   \mathbb{E}\left[\sum_{k=i}^{n-1}\abs{\delta F_k}_{\infty}^2h_k +\abs{\delta \xi}^2 \right] + C_p\kappa^{4/p} \mathbb{E}\left[\sup_{0 \le k \le n, t_k \in \Re} \abs{\delta C_{t_k}}^{p}_{\infty} \right]^{2/p} .
 \end{align*}
\end{Proposition}
{
\begin{Remark}
 The fact that the number of reflection times $\kappa$ appears only in front of the costs term justify why we assume $L^p$ integrability only on costs in \HYP{Fd_p}. Indeed, more $p$ is large, more $\kappa^{4/p}$ slowly increase (in $\kappa$). However, assuming some $L^p$ integrability on terminal conditions and generators would have no significant impact on our study. 
\end{Remark}
}

\proof
We adapt to our setting the proof of Proposition 2.3 in \cite{Chassagneux-Elie-Kharroubi-10}. The proof is divided into three steps and relies heavily on the reinterpretation in terms of switching problems. We first introduce a convenient dominating process and then provide successively the controls on $\delta \cY^{\Re,\pi}$ and $\delta \cZ^{\Re,\pi}$ terms.

\paragraph{Step 1. Introduction of an auxiliary backward scheme.}
Let us define $F:={}^1 F \vee {}^2 F$, $\xi := {}^1 \xi \vee {}^2 \xi$ and $C$ by $C^{ij}:= {}^1 C^{ij} \vee {}^2 C^{ij}$. Observe \HYP{Fd_p} holds for the data $(C,F,\xi)$ and $C$ satisfies the structure condition \eqref{condition de structure C}. We denote by $(\widetilde \cY^{\Re,\pi}, \cY^{\Re,\pi},\cZ^{\Re,\pi})$ the solution of the discretely obliquely reflected backward scheme with generator $F$, terminal condition $\xi$, reflection grid $\Re$ and cost process $C$. 

Using Proposition \ref{prop comparaison pour le schema} and the definition of $F$, $\xi$ and $C$, we obtain that
\begin{equation}
 \label{eq comparaison processus dans preuve stabilite}
 {\cY}^{\Re,\pi} \succcurlyeq {}^1{\cY}^{\Re,\pi} \vee {}^2{\cY}^{\Re,\pi}.
\end{equation}
Using Proposition \ref{prop snell representation schema}, we introduce switched backward schemes associated to ${}^1 \cY^{\Re,\pi}$, ${}^2 \cY^{\Re,\pi}$ and $\cY^{\Re,\pi}$ and denote by $\check{a}=(\check{\theta}_r,\check{\alpha}_r)_{r \geqslant 0}$ the optimal strategy related to $\cY^{\Re,\pi}$ starting from a fixed $(i,j) \in \llbracket 0,n \rrbracket \times \cI$. Therefore, we have
{
\begin{align}
 \nonumber
 ({\cY}^{\Re,\pi}_i)^j &=  \cU^{\Re,\pi,\check{a}}_{i}-\cA^{\check{a}}_i\\
 &= \xi^{\check{a}_n}+\sum_{k=i}^{n-1} F_k^{\check{a}_k}( \cV^{\Re,\pi,\check{a}}_k) h_k - \sum_{k=i}^{n-1} h_k\lambda_k^{-1} \cV^{\Re,\pi,\check{a}}_k H_k^\top -\sum_{k=i}^{n-1} \Delta \cM^{\check{a}}_k
 \label{eq switched backward scheme check sans esp cond preuve stabilite}
 -\cA^{\check{a}}_n   
\end{align}}

\paragraph{Step 2. Stability of the $Y$ component.}
Since $\check{a} \in \mathscr{A}^{\Re,\pi}_{i,j}$, we deduce from Proposition \ref{prop snell representation schema} (i) that, for $\ell \in \set{1,2}$,
{
\begin{align}
\nonumber
({}^\ell {\cY}^{\Re,\pi}_i)^j &\geqslant  {}^\ell\cU^{\Re,\pi,\check{a}}_{i}-{}^\ell\cA^{\check{a}}_i\\
&= {}^\ell\xi^{\check{a}_n}+\sum_{k=i}^{n-1} {}^\ell F_k^{\check{a}_k}( {}^\ell\cV^{\Re,\pi,\check{a}}_k) h_k - \sum_{k=i}^{n-1} h_k\lambda_k^{-1} {}^\ell\cV^{\Re,\pi,\check{a}}_k H_k^\top -\sum_{k=i}^{n-1} \Delta {}^\ell\cM^{\check{a}}_k
\label{eq U check 1 et 2} -{}^\ell\cA^{\check{a}}_n, 
\end{align}}
where ${}^\ell\cA^{\check{a}}$ is the process of cumulative costs $({}^\ell C^{ij})_{i,j \in \cI}$ associated to the strategy $\check{a}$. Combining this estimate with \eqref{eq comparaison processus dans preuve stabilite} and \eqref{eq switched backward scheme check sans esp cond preuve stabilite}, we derive 
{
\begin{equation}
 \label{eq comparaison processus dans preuve stabilite 2}
|({}^1 {\cY}^{\Re,\pi}_i)^j-({}^2 {\cY}^{\Re,\pi}_i)^j| \le | \cU^{\Re,\pi,\check{a}}_{i}-\cA^{\check{a}}_i - ({}^1\cU^{\Re,\pi,\check{a}}_{i}-{}^1\cA^{\check{a}}_i)|+| \cU^{\Re,\pi,\check{a}}_{i}-\cA^{\check{a}}_i - ({}^2\cU^{\Re,\pi,\check{a}}_{i}-{}^2\cA^{\check{a}}_i)|. 
\end{equation}}
Since both terms on the right-hand side of \eqref{eq comparaison processus dans preuve stabilite 2} are treated similarly, we focus on the first one and introduce discrete processes $\Gamma^{\check{a}} := \cU^{\Re,\pi,\check{a}} - \cA^{\check{a}}$ and ${}^1\Gamma^{\check{a}} := {}^1\cU^{\Re,\pi,\check{a}} - {}^1\cA^{\check{a}}$. Rewriting \eqref{eq switched backward scheme check sans esp cond preuve stabilite} and \eqref{eq U check 1 et 2} between $k$ and $k+1$ for $k \in \llbracket i, n-1 \rrbracket$, we get
\begin{align*}
 \Gamma^{\check{a}}_k -{}^1\Gamma^{\check{a}}_k &= \Gamma^{\check{a}}_{k+1} -{}^1\Gamma^{\check{a}}_{k+1} + [F_k^{\check{a}_k}(\cV^{\Re,\pi,\check{a}}_k)-{}^1F_k^{\check{a}_k}({}^1\cV^{\Re,\pi,\check{a}}_k)]h_k\\
 & \quad - h_k\lambda_k^{-1} [\cV^{\Re,\pi,\check{a}}_k-{}^1\cV^{\Re,\pi,\check{a}}_k] H_k^\top -[\Delta \cM^{\check{a}}_k - \Delta {}^1\cM^{\check{a}}_k].
\end{align*}
Using the identity $|y|^2=|x|^2 +2x(y-x)+|x-y|^2$, we obtain,
\begin{align*}
 &\quad \mathbb{E}_{t_k} [|\Gamma^{\check{a}}_{k+1} -{}^1\Gamma^{\check{a}}_{k+1}|^2]\\
 &= |\Gamma^{\check{a}}_k -{}^1\Gamma^{\check{a}}_k|^2 -2(\Gamma^{\check{a}}_k -{}^1\Gamma^{\check{a}}_k)(F_k^{\check{a}_k}(\cV^{\Re,\pi,\check{a}}_k)-{}^1F_k^{\check{a}_k}({}^1\cV^{\Re,\pi,\check{a}}_k))h_k\\
 &\quad + \mathbb{E}_{t_k} \left[\left|[F_k^{\check{a}_k}(\cV^{\Re,\pi,\check{a}}_k)-{}^1F_k^{\check{a}_k}({}^1\cV^{\Re,\pi,\check{a}}_k)]h_k - h_k\lambda_k^{-1} [\cV^{\Re,\pi,\check{a}}_k-{}^1\cV^{\Re,\pi,\check{a}}_k] H_k^\top -[\Delta \cM^{\check{a}}_k - \Delta {}^1\cM^{\check{a}}_k]\right|^2\right].
\end{align*}
Then, by the same reasoning as in the step 2 of the proof of Proposition \ref{prop estimee Y Z K schema general}, previous equality becomes 
\begin{align*}
 \quad \mathbb{E}_{t_k} [|\Gamma^{\check{a}}_{k+1} -{}^1\Gamma^{\check{a}}_{k+1}|^2] &\geqslant |\Gamma^{\check{a}}_k -{}^1\Gamma^{\check{a}}_k|^2 -2(\Gamma^{\check{a}}_k -{}^1\Gamma^{\check{a}}_k)(F_k^{\check{a}_k}(\cV^{\Re,\pi,\check{a}}_k)-{}^1F_k^{\check{a}_k}({}^1\cV^{\Re,\pi,\check{a}}_k))h_k\\
 &\quad + \frac{d}{\Lambda} h_k |\cV^{\Re,\pi,\check{a}}_k-{}^1\cV^{\Re,\pi,\check{a}}_k|^2,
\end{align*}
and we obtain, by summing over $k$ and taking expectation,
\begin{align*}
 & \quad \mathbb{E} \left[|\Gamma^{\check{a}}_{i} -{}^1\Gamma^{\check{a}}_{i}|^2 + \sum_{k=i}^{n-1} h_k |\cV^{\Re,\pi,\check{a}}_k-{}^1\cV^{\Re,\pi,\check{a}}_k|^2\right] \\
 &\le C\mathbb{E} \bigg[|\Gamma^{\check{a}}_n -{}^1\Gamma^{\check{a}}_n|^2 + \sum_{k=i}^{n-1} |\Gamma^{\check{a}}_k -{}^1\Gamma^{\check{a}}_k| |F_k^{\check{a}_k}(\cV^{\Re,\pi,\check{a}}_k)-{}^1F_k^{\check{a}_k}({}^1\cV^{\Re,\pi,\check{a}}_k)|h_k\bigg].
\end{align*}
Since $F = {}^1 F \vee {}^2 F$ and ${}^1 F$ is a Lipschitz function, we also get
\begin{align*}
 |F_k^{\check{a}_k}(\cV^{\Re,\pi,\check{a}}_k)-{}^1F_k^{\check{a}_k}({}^1\cV^{\Re,\pi,\check{a}}_k)| &\le |\delta F_k|_{\infty} +C|\cV^{\Re,\pi,\check{a}}_k-{}^1\cV^{\Re,\pi,\check{a}}_k|,
\end{align*}
and then, by using Young's inequality and discrete Gronwall's lemma (see \cite{Clark-87}), we deduce from the last and the penultimate inequalities that
{
\begin{align}
 \label{eq U-1U}
 \mathbb{E} \left[|\Gamma^{\check{a}}_{i} -{}^1\Gamma^{\check{a}}_{i}|^2\right] &\le C \left(  \mathbb{E}\left[|\delta \xi|^2+ \sum_{k=i}^{n-1} \abs{\delta F_k}_{\infty}^2h_k + |{}\cA^{\check{a}}_n-{}^1\cA^{\check{a}}_n|^2
 \right] \right).
\end{align}}
Moreover we compute
\begin{align*}
 \mathbb{E}[|{}\cA^{\check{a}}_n-{}^1\cA^{\check{a}}_n|^2] &\le \mathbb{E} \left[|\cN^{\check{a}}|^2 \sup_{0 \le m \le n, t_m \in \Re} |\delta C_{t_m}|^2_{\infty}\right].
\end{align*}
If $p=2$, then $\cN^{\check{a}} \le \kappa$ yields
 $$\mathbb{E}[|{}\cA^{\check{a}}_n-{}^1\cA^{\check{a}}_n|^2] \le \kappa^2 \mathbb{E} \left[ \sup_{0 \le m \le n, t_m \in \Re} |\delta C_{t_m}|^2_{\infty}\right].$$
Otherwise, from Proposition \ref{prop estimee U V N schema}, H\"older's inequality and the fact that $\cN^{\check{a}} \le \kappa$, we deduce
\begin{align*}
 \mathbb{E}[|{}\cA^{\check{a}}_n-{}^1\cA^{\check{a}}_n|^2] &\le \mathbb{E} \left[|\cN^{\check{a}}|^{\frac{2p}{p-2}}\right]^{\frac{p-2}{p}} \mathbb{E}\left[ \sup_{0 \le m \le n, t_m \in \Re} |\delta C_{t_m}|^p_{\infty}\right]^{2/p}\\
 &\le \mathbb{E} \left[\kappa^{\frac{2p}{p-2}-2}|\cN^{\check{a}}|^{2}\right]^{\frac{p-2}{p}} \mathbb{E}\left[ \sup_{0 \le m \le n, t_m \in \Re} |\delta C_{t_m}|^p_{\infty}\right]^{2/p}\\
 &\le C_p \kappa^{4/p}  \mathbb{E}\left[ \sup_{0 \le m \le n, t_m \in \Re} |\delta C_{t_m}|^p_{\infty}\right]^{2/p}.
\end{align*}
Inserting the last estimate into \eqref{eq U-1U}, we get
{
\begin{align*}
 \mathbb{E} \left[|\Gamma^{\check{a}}_{i} -{}^1\Gamma^{\check{a}}_{i}|^2\right] &\le C \mathbb{E}\left[|\delta \xi|^2+ \sum_{k=i}^{n-1} \abs{\delta F_k}_{\infty}^2h_k \right]  + C_p\kappa^{4/p}  \mathbb{E}\left[ \sup_{0 \le m \le n, t_m \in \Re} |\delta C_{t_m}|^p_{\infty}\right]^{2/p}.
\end{align*}}

By symmetry, we have the same estimate for {$\mathbb{E} \left[|\Gamma^{\check{a}}_{i} -{}^2\Gamma^{\check{a}}_{i}|^2\right]$}. Therefore, from \eqref{eq comparaison processus dans preuve stabilite 2} and the fact that $j$ is arbitrary, we deduce the wanted estimate for {$\mathbb{E}[|\delta  \cY^{\Re,\pi}_i|^2]$}. 
{
Concerning the upper bound of $\mathbb{E}[|\delta \widetilde \cY^{\Re,\pi}_i|^2]$, we just have to use the scheme \eqref{de scheme aleatoire general} when $0 \le i \le n-1$:
$$\delta \widetilde \cY^{\Re,\pi}_i = \mathbb{E}\left[\delta \cY^{\Re,\pi}_{i+1} |\mathcal{F}_{t_i}\right] + h_i\left( {}^1 F_i^{\pi}({}^1\cZ_i^{\Re,\pi})-{}^2 F_i^{\pi}({}^2\cZ_i^{\Re,\pi}) \right),$$
so we get, by using the Lipschitz regularity of ${}^1 F_i^{\pi}$, H\"older inequality and \ref{eq structure H 3}, 
\begin{align*}
 \mathbb{E}\left[|\delta \widetilde \cY^{\Re,\pi}_i|^2\right] &\le C\left( \mathbb{E}\left[|\delta  \cY^{\Re,\pi}_{i+1}|^2\right]+h_i|\delta F_i|^2_{\infty}+h_i\mathbb{E}\left[\mathbb{E}_{t_i}\left[|\delta  \cY^{\Re,\pi}_{i+1}|^2\right]\mathbb{E}_{t_i}\left[|H_i|^2\right]\right]\right)\\ 
 &\le C\left(\mathbb{E}\left[|\delta  \cY^{\Re,\pi}_{i+1}|^2\right]+|\delta F_i|^2_{\infty}h_i\right).
\end{align*}
Then we just have to use the estimate already obtained for $\mathbb{E}[|\delta  \cY^{\Re,\pi}_{i+1}|^2]$ to conclude.
}

\paragraph{Step 3. Stability of the $Z$ component.}

Observing that $\delta \cZ_k^{\Re,\pi} = \mathbb{E}_{t_k}[(\delta \cY^{\Re,\pi}_{k+1} - \mathbb{E}_{t_{k}}[\delta \cY^{\Re,\pi}_{k+1}])H_k^{\top}]$, one computes
\begin{align} \label{eq stab Z 2}
h_k |\delta \cZ^{\Re,\pi}_k|^2 \le C\EFp{\tk{}}{\left|\delta {\cY}^{\Re,\pi}_{k+1}\right|^2 - \left|\EFp{\tk{}}{\delta {\cY}^{\Re,\pi}_{k+1}}\right|^2}
\end{align}
From the scheme's definition, we have
\begin{align*}
\left|\EFp{\tk{}}{\delta {\cY}^{\Re,\pi}_{k+1}}\right|^2 \ge \left|\delta \widetilde{\cY}^{\Re,\pi}_{k}\right|^2 -2\left|\delta \widetilde{\cY}^{\Re,\pi}_{k}\left\{{}^1F_k({}^1\cZ_k^{\Re,\pi})-{}^2F_k({}^2\cZ_k^{\Re,\pi})\right\}h_k\right|\;.
\end{align*}
Inserting the last estimate into \eqref{eq stab Z 2} and using \HYP{Fd_p}, we obtain,
\begin{align*}
h_k |\delta \cZ^{\Re,\pi}_k|^2 \le C\left(\EFp{\tk{}}{|\delta {\cY}^{\Re,\pi}_{k+1}|^2}-|\delta \widetilde{\cY}^{\Re,\pi}_{k}|^2 +  Ch_k| \delta \widetilde{\cY}^{\Re,\pi}_{k}|^2 + \frac{1}{2} h_k |\delta \cZ^{\Re,\pi}_k|^2 + h_k |\delta F_k|^2_{\infty} \right)\;.
\end{align*}
Taking expectation on both sides and summing over $k$,  we get
\begin{align*}
& \quad \frac{1}{2}\sum_{k=i}^{n-1} h_k\esp{ |\delta \cZ^{\Re,\pi}_k|^2}\\
&\le C\left( \esp{|\delta {\cY}^{\Re,\pi}_{n}|^2} + \sum_{\substack{k=i\\t_k \in \Re}}^{n-1} \esp{|\delta {\cY}^{\Re,\pi}_{k}|^2-|\delta \widetilde{\cY}^{\Re,\pi}_{k}|^2} 
+  \max_{i \le k \le n-1}\esp{|\delta {\widetilde{\cY}}^{\Re,\pi}_{k}|^2} 
 + \sum_{k=i}^{n-1} h_k |\delta F_k|^2_{\infty} \right)\;,
 \\
 &\le C\left( \esp{|\delta {\cY}^{\Re,\pi}_{n}|^2} +  \kappa \max_{i \le k \le n-1}\esp{|\delta {\cY}^{\Re,\pi}_{k}|^2+|\delta \widetilde{\cY}^{\Re,\pi}_{k}|^2} 
 + \sum_{k=i}^{n-1} h_k |\delta F_k|^2_{\infty} \right)\;.
\end{align*}
The proof is concluded using estimates on $\delta \widetilde{\cY}^{\Re,\pi}$ and $\delta {\cY}^{\Re,\pi}$ already obtained in the first part of the proof. 
\eproof

\vspace{10pt}
We will now use this general stability result on obliquely reflected backward schemes to obtain a $L^2$-stability result for the scheme \eqref{schemeintro} (see \cite{Chassagneux-Crisan-13} for a general definition of $L^2$-stability for backward schemes). Firstly, we introduce a perturbed version of the scheme given in \eqref{schemeintro}.
\begin{Definition}
\label{de scheme perturbe}
 \begin{enumerate}[(i)]
  \item The terminal condition is given by a $\mathcal{F}_T$-measurable random variable $\bar{Y}_n \in \mathscr{L}^2$;
  \item for $0 \le i <n$,
\begin{equation}\label{eq scheme perturbe}
\begin{cases}
\bar{Z}_i^{\Re,\pi} := \mathbb{E}[\bar{Y}^{\Re,\pi}_{i+1}H_i\mid\mathcal{F}_{t_i}] ,
\vspace*{2pt}\cr
\widetilde{\bar{Y}}^{\Re,\pi}_{i} := \mathbb{E}[\bar{Y}^{\Re,\pi}_{i+1}\mid\mathcal{F}_{t_i
}] + h_i
f(\Xp_{t_i},\widetilde{\bar{Y}}^{\Re,\pi}_{i}, \bar{Z}_i^{\Re,\pi}) +\zeta^f_i,
\vspace*{2pt}\cr
\bar{Y}_i^{\Re,\pi} := \widetilde{\bar{Y}}_i^{\Re,\pi} \mathbf{1}_{ \{ t_i\notin\Re\} } + \bar{\cP}_{t_i}
(\Xp_{t_i},
\widetilde{\bar{Y}}_i^{\Re,\pi})\mathbf{1}_{ \{ t_i\in\Re\} } ,
\end{cases}
\end{equation}
with $\bar{\cP}$ the oblique projection 
\begin{equation*}
 \bar{\cP}_{t_i} : (x,y) \in \mathbb{R}^d \times \mathbb{R}^d \mapsto \left( \max_{j \in \cI} \set{ y^j -\bar{c}_{t_i}^{kj}(x) } \right)_{1 \le k \le d},
\end{equation*}
associated to costs $\bar{c}_{t_i}(x) := c(x)+\zeta^c_{t_i}$. Perturbations $\zeta^Y_i:=(\zeta^f_i,\zeta^c_{t_i})$ are $\mathcal{F}_{t_i}$-measurable and square integrable random variables. Moreover we assume that the random costs $(\bar{c}_{t_i}(X_{t_i}))_{0 \le i \le n}$ satisfy the structure conditions \eqref{condition de structure C}.
\end{enumerate}
\end{Definition}

Setting $\delta Y_i = {Y}_i^{\Re,\pi} - \bar{Y}_i^{\Re,\pi}$, $\delta \widetilde Y_i = \widetilde{Y}_i^{\Re,\pi} - \widetilde{\bar{Y}}_i^{\Re,\pi}$ and $\delta Z_i = {Z}_i^{\Re,\pi} - \bar{Z}_i^{\Re,\pi}$, we obtain the following $L^2$-stability result for the scheme \eqref{schemeintro}.
\begin{Proposition}
 \label{prop stabilite L2 schema}
Assume that \HYP{f} is in force and, for all $p \ge 2$,  
\begin{align}
\label{eq condition iltegrabilite zetaf}
\esp{|\bar{Y}_n|^2+ \sum_{i=0}^{n-1} |\zeta^f_i|^2 + \sup_{0 \le i \le n} \abs{\zeta^c_{t_i}}^p} \le C.
\end{align}
We also assume that $|\pi|L^Y <1$ and
\begin{align}
\label{eq condition Hi pour stabilite L2 du schema}
 \left(\sup_{0 \le i \le n-1} h_i\abs{H_i}\right) L^Z \le 1.
\end{align}
Then schemes \eqref{schemeintro} and \eqref{eq scheme perturbe} are well defined and the following $L^2$-stability holds true, for all $p\ge 2$,
\begin{eqnarray}
 & & \sup_{0 \le i \le n} \esp{|\delta Y_i|^2 + |\delta \widetilde Y_i|^2} + \frac1{\kappa}\sum_{i=0}^{n-1} h_i\esp{ |\delta Z_i|^2} \nonumber  \\
 \label{eq stab YZ pr}
&\le& C\left(\esp{|\delta Y_n|^2}+\sum_{i=0}^{n-1}\frac{1}{h_i}\esp{|\zeta_i^f|^2}\right)+C_p\kappa^{4/p} \esp{\sup_{\substack{0 \le i \le n\\t_i \in \Re}} |\zeta_i^c|^p}^{2/p}.
\end{eqnarray}
 \end{Proposition}

\proof
Since we have assumed $|\pi| L^Y<1$, then a simple fixed point argument shows that schemes \eqref{schemeintro} and \eqref{eq scheme perturbe} are well defined, i.e. there exists a unique solution to each scheme.

For the $L^2$-stability, we want to apply Proposition \ref{prop stabilite schema gene aleatoire} with ${}^1\xi = g(X_T^{\pi})$, ${}^2\xi= \bar{Y}_n$, ${}^1 F_i(z)=f(\Xp_{t_i},\widetilde{{Y}}^{\Re,\pi}_{i},z)$, ${}^2 F_i(z)=f(\Xp_{t_i},\widetilde{\bar{Y}}^{\Re,\pi}_{i},z) + \zeta^f_i$, ${}^1 C_{t_i} = c(X^{\pi}_{t_i})$ and ${}^2 C_{t_i} = c(X^{\pi}_{t_i})+ \zeta^c_{t_i}$. To do this, we have to check that assumption \HYP{Fd_p} is fulfilled for these two obliquely reflected backward schemes. Firstly, we have assumed
$$ \left(\sup_{0 \le i \le n-1} h_i\abs{H_i}\right) L^Z \le 1.$$
Moreover, hypothesis \HYP{f}, assumption \eqref{eq condition iltegrabilite zetaf} and classical estimates for processes $X$ and $X^{\pi}$ leads to 
\begin{align*}
&\esp{|{}^1\xi|^2+|{}^2\xi|^2 + \sum_{i=0}^{n-1} [|{}^1F_i(0)|^2+|{}^2F_i(0)|^2]h_i +\sup_{0 \le i\le n} [|{}^1 C_{t_i}|^p+|{}^2 C_{t_i}|^p]}\\
&\le C_p + C\esp{\sup_{0 \le i\le n} [|\widetilde{{Y}}^{\Re,\pi}_{i}|^2+|\widetilde{\bar{Y}}^{\Re,\pi}_{i}|^2]}.
\end{align*}
To estimate quantities $\esp{\sup_{0 \le i\le n} |\widetilde{{Y}}^{\Re,\pi}_{i}|^2}$ and $\esp{\sup_{0 \le i\le n} |\widetilde{\bar{Y}}^{\Re,\pi}_{i}|^2}$, we just have to rewrite slightly the first step of the proof of Proposition \ref{prop estimee Y Z K schema general}. The beginning of the proof stays true: \eqref{eq estimee Ytilde dans Y Z K} yields, for all $i \in \llbracket 0,n\rrbracket$,
\begin{align*}
 \esp{\sup_{i \le k \le n} \left[|\widetilde{{Y}}^{\Re,\pi}_{k}|^2 + |\widetilde{\bar{Y}}^{\Re,\pi}_{k}|^2\right]} &\le C\esp{ |{}^1\xi|^2 + |{}^2\xi|^2 + \sum_{k=i}^{n-1} \left[ |{}^1F_k(0)|^2+|{}^2F_k(0)|^2\right]h_k}\\
 &\le C\left(1+\sum_{k=i}^{n-1} \esp{\sup_{k \le m \le n} \left[|\widetilde{{Y}}^{\Re,\pi}_{m}|^2 + |\widetilde{\bar{Y}}^{\Re,\pi}_{m}|^2\right]}h_k \right). 
\end{align*}
Thus, the discrete Gronwall lemma (see \cite{Clark-87}) allows to conclude that 
\begin{align*}
 \esp{\sup_{0 \le k \le n} \left[|\widetilde{{Y}}^{\Re,\pi}_{k}|^2 + |\widetilde{\bar{Y}}^{\Re,\pi}_{k}|^2\right]} &\le C
\end{align*}
and then assumption \HYP{Fd_p} is fulfilled. Proposition \ref{prop stabilite schema gene aleatoire} and \HYP{f} imply, for all $i \in \llbracket 0,n\rrbracket$,
\begin{align*} 
& \sup_{i \le k \le n} \esp{|\delta Y_k|^2 + |\delta \widetilde Y_k|^2} + \frac1{\kappa}\sum_{k=i}^{n-1} h_k\esp{ |\delta Z_k|^2} \\
\le& C\left(\esp{|\delta Y_n|^2}+\sum_{k=i}^{n-1}|\delta F_k|_{\infty}^2\right)+C_p\kappa^{4/p} \esp{\sup_{\substack{0 \le k \le n\\t_k \in \Re}} |\zeta_{t_k}^c|^p}^{2/p}\\
\le&  C\left(\esp{|\delta Y_n|^2}+\sum_{k=0}^{n-1}\frac{1}{h_k}\esp{|\zeta_k^f|^2} +\sum_{k=i}^{n-1} \sup_{k \le m \le n} \esp{|\delta Y_m|^2 + |\delta \widetilde Y_m|^2}h_k\right)\\
&+C_p\kappa^{4/p} \esp{\sup_{\substack{0 \le k \le n\\t_k \in \Re}} |\zeta_{t_k}^c|^p}^{2/p}.
\end{align*}
Applying the discrete Gronwall lemma to the last inequality completes the proof.
\eproof

\subsection{Convergence analysis of the discrete-time approximation}
\label{section etude schema sous section convergence}
We will give now the main result of this section that provides an upper bound for the error between the obliquely reflected backward scheme \eqref{schemeintro} and the discretely obliquely reflected BSDE \eqref{BSDEDORintro}.
\begin{Theorem}
\label{th convergence scheme vers discret refl}
 Assume that \HYP{f} is in force. We also assume that $|\pi|L^Y<1$ and weights $(H_i)_{0 \le i \le n-1}$ are given by
\begin{align}
\label{eq def HiR}
 (H_i)^{\ell} = \frac{-R}{h_i} \vee \frac{W^\ell_{t_{i+1}} -W^\ell_{t_i}}{h_i} \wedge \frac{R}{h_i}, \quad 1 \le \ell \le d,
\end{align}
 with $R$ a positive parameter such that $RL^Z \le 1$.
 Then the following holds:
 \begin{align*}
  \sup_{0 \le i \le n} \esp{|\widetilde{Y}^{\Re}_{t_i} -\widetilde{Y}^{\Re,\pi}_i|^2 + |{Y}^{\Re}_{t_i} -{Y}^{\Re,\pi}_i|^2} + \frac{1}{\kappa}\esp{ \sum_{i=0}^{n-1}\int_{t_i}^{t_{i+1}} |Z_s^{\Re} - Z_i^{\Re,\pi}|^2\ud s}\le C_R \left( |\pi|^{1/2} +\kappa |\pi|\right).
 \end{align*}
\end{Theorem}
\proof
\paragraph{Step 1. Expression of the perturbing error.}
Since we want to apply Proposition \ref{prop stabilite L2 schema}, we first observe that $(Y^{\Re},Z^{\Re})$ can be rewritten as a perturbed obliquely reflected backward scheme. Namely, setting $\bar{Y}_i := Y^{\Re}_{t_i}$ and $\widetilde{\bar{Y}}_i := \widetilde{Y}^{\Re}_{t_i}$, for all $i \in \llbracket 0,n\rrbracket$, we have

\begin{equation}\label{eq edsr discretem reflechie comme schema perturbe }
\begin{cases}
\bar{Z}_i^{} := \mathbb{E}[\bar{Y}^{}_{i+1}H_i\mid\mathcal{F}_{t_i}] ,
\vspace*{2pt}\cr
\widetilde{\bar{Y}}^{}_{i} := \mathbb{E}[\bar{Y}^{}_{i+1}\mid\mathcal{F}_{t_i
}] + h_i
f(\Xp_{t_i},\widetilde{\bar{Y}}^{}_{i}, \bar{Z}_i) +\zeta^f_i,
\vspace*{2pt}\cr
\bar{Y}_i^{} := \widetilde{\bar{Y}}_i^{} \mathbf{1}_{ \{ t_i\notin\Re\} } + \bar{\cP}_{t_i}
(\Xp_{t_i},
\widetilde{\bar{Y}}_i^{})\mathbf{1}_{ \{ t_i\in\Re\} } ,
\end{cases}
\end{equation}
with
\begin{align*}
 \zeta_i^f= \mathbb{E}_{t_i}\left[ \int_{t_i}^{t_{i+1}} \left(f(X_s,\widetilde{Y}_s^{\Re},Z_s^{\Re})-f(\Xp_{t_i},\widetilde{{Y}}^{\Re}_{t_i}, \bar{Z}_i) \right) \ud s\right]\quad \textrm{and} \quad \zeta^c_{t_i}= c(X_{t_i})-c(X^{\pi}_{t_i}).
\end{align*}
Let us check that \eqref{eq condition iltegrabilite zetaf} is fulfilled for all $p \ge 2$: using \HYP{f}, Proposition \ref{prop estimee Y Z K} and classical estimates for $X$ and $X^{\pi}$, we get
\begin{align*}
&\esp{|\bar{Y}_n|^2+ \sum_{i=0}^{n-1} |\zeta^f_i|^2 + \sup_{0 \le i \le n} \abs{\zeta^c_{t_i}}^p}\\
&\le C_p\esp{1+\sup_{s \in [0,T]} |X_s|^p + \sup_{i \in \llbracket 0,n \rrbracket} |X_{t_i}^{\pi}|^{p} + \sup_{s \in [0,T]} |\tilde{Y}^{\Re}_s|^2 + \int_0^T |Z_s^{\Re}|^2\ud s +\sum_{i=0}^{n-1} |\tilde{Y}^{\Re}_{t_{i+1}} H_i h_i|^2 }\\
&\le C_p\left( 1+ \esp{\sup_{s \in [0,T]} |\tilde{Y}^{\Re}_s|^4 + \left( \sum_{i=0}^{n-1} | H_i h_i|^2\right)^2}\right).
\end{align*}
Applying Burkholder-Davis-Gundy inequality, we have $\esp{\left( \sum_{i=0}^{n-1} | H_i h_i|^2\right)^2} \le C$ and so \eqref{eq condition iltegrabilite zetaf} is fulfilled. 
Finally, we easily check that \eqref{eq def HiR} implies \eqref{eq condition Hi pour stabilite L2 du schema}.
\paragraph{Step 2. Discretization error for the $Y$ component.}
Setting $p=4$, we apply Proposition \ref{prop stabilite L2 schema} and  get by direct calculations
\begin{align*}
 &\sup_{0 \le i \le n} \esp{|\widetilde{Y}^{\Re}_{t_i} -\widetilde{Y}^{\Re,\pi}_i|^2 + |{Y}^{\Re}_{t_i} -{Y}^{\Re,\pi}_i|^2}\\
 &\le C\left(\esp{|g(X_T)-g(X_T^{\pi})|^2}+\sum_{i=0}^{n-1}\frac{1}{h_i}\esp{|\zeta_i^f|^2}\right)+C_p\kappa \esp{\sup_{\substack{0 \le i \le n\\t_i \in \Re}} |\zeta_{t_i}^c|^4}^{1/2}\\
 &\le  C\esp{|X_T-X^{\pi}_T|^2} + C\sup_{0 \le i \le n-1} \esp{\sup_{s \in [t_i,t_{i+1}]} |X_s-X_{t_i}^{\pi}|^2} + C\esp{\int_0^T |\widetilde{Y}^{\Re}_s -\widetilde{Y}^{\Re}_{\pi(s)}|^2\ud s}\\
 &\quad + C\esp{\sum_{i=0}^{n-1} \int_{t_i}^{t_{i+1}} |{Z}^{\Re}_s -\bar{Z}^{\Re}_{t_i}|^2\ud s}+ C\esp{\sum_{i=0}^{n-1} |\bar{Z}^{\Re}_{t_i} -\bar{Z}_{i}|^2 h_i}+ C\kappa \esp{\sup_{0 \le i \le n-1} |X_{t_i}-X_{t_i}^{\pi}|^4}^{1/2} .
\end{align*}
Classical estimations on the Euler scheme for SDEs, see e.g. \cite{klopla92}, yield
$$\esp{|X_T-X^{\pi}_T|^2} + \sup_{0 \le i \le n-1} \esp{\sup_{s \in [t_i,t_{i+1}]} |X_s-X_{t_i}^{\pi}|^2} +\kappa \esp{\sup_{0 \le i \le n-1} |X_{t_i}-X_{t_i}^{\pi}|^4}^{1/2}  \le C\kappa|\pi|.$$
Applying Proposition \ref{prop regularite YR} and Proposition \ref{prop regularite ZR}, we obtain
$$\esp{\int_0^T |\widetilde{Y}^{\Re}_s -\widetilde{Y}^{\Re}_{\pi(s)}|^2ds}+\esp{\sum_{i=0}^{n-1} \int_{t_i}^{t_{i+1}} |{Z}^{\Re}_s -\bar{Z}^{\Re}_{t_i}|^2\ud s} \le C(|\pi|+|\pi|^{1/2}+\kappa |\pi|).$$
It remains to bound the term:
\begin{align*}
 &\esp{\sum_{i=0}^{n-1} |\bar{Z}^{\Re}_{t_i} -\bar{Z}_{i}|^2 h_i}\\ 
 \le& 2\underbrace{\esp{\sum_{i=0}^{n-1} \left|\bar{Z}^{\Re}_{t_i} -\mathbb{E}_{t_i}\left[Y_{t_{i+1}}^{\Re}\frac{\Delta W_i}{h_i}\right]\right|^2 h_i}}_{:=A} + 2\underbrace{\esp{\sum_{i=0}^{n-1} \left|\mathbb{E}_{t_i}\left[Y_{t_{i+1}}^{\Re}\left(\frac{\Delta W_i}{h_i}-H_i\right)\right]\right|^2 h_i}}_{:=B}.
\end{align*}
By remarking that $\bar{Z}^{\Re}_{t_i} = \mathbb{E}_{t_i}\left[\int_{t_i}^{t_{i+1}} Z^{\Re}_s dW_s \frac{\Delta W_i}{h_i}\right]$, we have
\begin{align*}
 A &=\esp{\sum_{i=0}^{n-1} \left|\mathbb{E}_{t_i}\left[\int_{t_i}^{t_{i+1}}f(X_s,\widetilde{Y}_s^{\Re},Z_s^{\Re})\ud s\frac{\Delta W_i}{h_i}\right]\right|^2 h_i} \\
 &\le \esp{\sum_{i=0}^{n-1} h_i \int_{t_i}^{t_{i+1}}|f(X_s,\widetilde{Y}_s^{\Re},Z_s^{\Re})|^2\ud s} \le |\pi|\esp{ \int_{0}^{T}|f(X_s,\widetilde{Y}_s^{\Re},Z_s^{\Re})|^2\ud s} \le C|\pi|.
\end{align*}
Finally, we also get by standard calculations, Proposition \ref{prop estimee Y Z K} and classical results about Gaussian distribution tails
\begin{align*}
 B &\le \sup_{0 \le i \le n-1} \esp{|Y^{\Re}_{t_{i+1}}|^2} \times \sup_{0 \le i \le n-1} \esp{\left|\frac{\Delta W_i}{h_i} -H_i\right|^2}\le C\left(\frac{2d}{h_i} \int_{Rh_i^{-1}}^{+\infty} x^2\frac{e^{-x^2/2}}{\sqrt{2\pi}} \ud x \right)\\
 &\le \frac{C}{h_i}\left(Rh_i^{-1}e^{-R^2h_i^{-2}/2}\right) \le \frac{CR}{h_i^2}\left(\frac{2h_i^2}{R^2}\right)^{3/2}\le C_R |\pi|.
\end{align*}

\paragraph{Step 3. Discretization error for the $Z$ component.}
Let us remark that we have
\begin{align*}
\frac{1}{\kappa}\esp{ \sum_{i=0}^{n-1}\int_{t_i}^{t_{i+1}} |Z_s^{\Re} - Z_i^{\Re,\pi}|^2\ud s}\le& \frac{1}{\kappa}\esp{ \sum_{i=0}^{n-1}\int_{t_i}^{t_{i+1}} |Z_s^{\Re} - \bar{Z}_{t_i}^{\Re}|^2\ud s}+\frac{1}{\kappa}\esp{ \sum_{i=0}^{n-1}|\bar{Z}_{t_i}^{\Re} - \bar{Z}_i|^2h_i}\\
&+\frac{1}{\kappa}\esp{ \sum_{i=0}^{n-1} |\bar{Z}_i - Z_i^{\Re,\pi}|^2h_i}.
\end{align*}
Previous calculations already yield
$$\frac{1}{\kappa}\esp{ \sum_{i=0}^{n-1}\int_{t_i}^{t_{i+1}} |Z_s^{\Re} - \bar{Z}_{t_i}^{\Re}|^2\ud s}+\frac{1}{\kappa}\esp{ \sum_{i=0}^{n-1}|\bar{Z}_{t_i}^{\Re} - \bar{Z}_i|^2h_i} \le \frac{C}{\kappa} \left( |\pi|^{1/2} +\kappa |\pi|\right).$$
Moreover, we apply Proposition \ref{prop stabilite L2 schema} to obtain
$$\frac{1}{\kappa}\esp{ \sum_{i=0}^{n-1} |\bar{Z}_i - Z_i^{\Re,\pi}|^2h_i} \le C \left( |\pi|^{1/2} +\kappa |\pi|\right),$$
thanks to estimates obtained in step 2.
\eproof

 %
 %


\section{Application to continuously reflected BSDEs}
This section is devoted to the study of the error between the scheme \eqref{schemeintro} and the continuously obliquely reflected BSDEs \eqref{eqBSDECORIntro}. An upper bound of this error is stated in Subsection \ref{section continuous time subsection convergence} while Subsection \ref{section continuous time subsection discretely2continuous} is devoted to the error between the continuously obliquely reflected BSDEs \eqref{eqBSDECORIntro} and the discretely obliquely reflected BSDEs \eqref{BSDEDORintro}. Before these results, we start by giving some classical estimates on the solution of \eqref{eqBSDECORIntro}.

 \begin{Proposition}
 \label{prop estimee Y Z K continu}
 Assume that \HYP{f} is in force. There exists a unique solution $(Y,Z,K) \in \mathscr{S}_2 \times \mathscr{H}_2 \times \mathscr{K}_2$ to \eqref{eqBSDECORIntro} and it satisfies, for all $p \geq 2$, 
 $$|Y|_{\mathscr{S}_p}+ |Z|_{\mathscr{H}_p}+ |K_T|_{\mathscr{L}^p} \leq C_p.$$
 \end{Proposition}
 \proof
 The existence and uniqueness result comes from \cite{Chassagneux-Elie-Kharroubi-11}.
 Concerning estimates, we want to apply Proposition \ref{prop existence unicite EDSRs obliqu refle cont et discr avec generateur aleatoire} with terminal condition $\xi = g(X_T)$, random generator $F(s,z)=f(X_s,Y_s,z)$ and costs $C^{ij}_s=c^{ij}(X_s)$. So, we just have to show that \HYP{F_p} is in force. Thus, using the fact that $f$ is a Lipschitz function with respect to $y$, it is sufficient to control $\widetilde Y^{\Re}$ in $\mathscr{S}^p$ to conclude. We are able to obtain estimates on $|\widetilde Y^{\Re}|_{\mathscr{S}^p}$ by a direct adaptation to the continuous time setting of the proof of Proposition \ref{prop estimee Y Z K}.
 \eproof

\subsection{Error between discretely and continuously reflected BSDEs}
\label{section continuous time subsection discretely2continuous}
We show here that the error between the continuously reflected BSDE \eqref{eqBSDECORIntro} and the discretely reflected BSDE \eqref{BSDEDORintro} is controlled in a convenient way. We start by introducing a temporary assumption.

\HYP{z} For all $(x,y,z) \in \R^d \times \R^d \times \cM^{d,d}$, 
$|f(x,y,z)| \leq C(1+|x|+|y|).$

\vspace{2mm}

\begin{Proposition}
\label{pr theorem 5.2}
Assume that \HYP{f} and \HYP{z} are in force, then
\begin{align*}
 \esp{\sup_{t \in [0,T]}|Y_t-Y^\Re_t|^2 + \sup_{t \in [0,T]}|Y_t-\tilde{Y}^\Re_t|^2} \le C |\Re|\log(2T/|\Re|).
\end{align*}
Moreover, if the cost functions are constant, we  obtain a better rate of convergence, namely,
\begin{align*}
 \esp{\sup_{t \in [0,T]}|Y_t-Y^\Re_t|^2 + \sup_{t \in [0,T]}|Y_t-\tilde{Y}^\Re_t|^2} \le C |\Re|.
\end{align*}

\end{Proposition}

\proof
1. We denote $(\ckY,\ckZ,\ckK)$ the solution of an auxiliary continuously obliquely reflected BSDE with cost functions $c$, with terminal condition $\xi := g(X_T)$ and whose random generator is given by
\begin{align*}
 \ckf(s,z) = f(X_s,Y_s,z) \vee f(X_s,\tYR_s,z) \;.
\end{align*}
We also denote $(\tilde{Y},\tilde{Z},\tilde{K})$ the solution of the continuously obliquely reflected BSDE with cost functions $c$, with terminal condition $\xi := g(X_T)$ and with random driver $\tilde{f}(s,z) = f(X_s,\tYR_s,z)$.
From Proposition \ref{th rep switch generic}, we know that each component of $\ckY$, $Y$ and $\tilde{Y}$ can be represented as optimal values of some control problem namely
{
\begin{align}
\nonumber
(\ckY_t)^i &= \esssup_{a \in \mathscr{A}_{i,t}} \left(\ckU^a_t-A^a_t\right) =\ckU^{\cka}_t-A^{\cka}_t,\\ 
(Y_t)^i &= \esssup_{a \in \mathscr{A}_{i,t}} \left(U^a_t-A^a_t\right), \quad (\tilde{Y}_t)^i = \esssup_{a \in \mathscr{A}_{i,t}} \left({U}^{\Re,a}_t-A^a_t\right),
\label{eq representation ckY Y YR}
\end{align}}
with $t \in [0,T]$, $i \in \cI$, $\ckU^a$, $U^a$ and $U^{\Re,a}$ solutions to following ``switched'' BSDEs:
\begin{align}
\label{eq U check}
\ckU^a_t &= \xi^{a_T} + \int_t^T\check{f}^{a_s}(s,\ckV^a_s)\ud s -\int_t^T\ckV^a_s \ud W_s - A^a_T+ A^a_t\;,\\
\label{eq U pas check}
U^a_t &= \xi^{a_T} + \int_t^T f^{a_s}(X_s,Y_s,V^a_s)\ud s -\int_t^T V^a_s \ud W_s - A^a_T+ A^a_t\;,\\
U^{\Re,a}_t &= \xi^{a_T} + \int_t^T f^{a_s}(X_s,Y^{\Re}_s,V^{\Re,a}_s)\ud s -\int_t^T V^{\Re,a}_s \ud W_s - A^a_T+ A^a_t\;,
\end{align} 
and $\check{a}$ the optimal strategy given by Proposition \ref{th rep switch generic}. {As previously, we denote $N^a$ the number of switches associated to the strategy $a \in \mathcal{A}_{i,t}$.}
Using a comparison argument, we easily check that $\ckU^a \ge U^a \vee U^{\Re,a}$, for any strategy $a \in \mathscr{A}_{i,t}$. This estimate combined with \eqref{eq representation ckY Y YR} leads to 
\begin{align*}
\ckY^\ell \ge Y^\ell \vee (\tilde{Y})^\ell \; \text{ for all }\; \ell \in \set{1,\dots,d}\;.
\end{align*}
Moreover, Corollary \ref{co rep switch cor and dor} and \eqref{eq representation ckY Y YR} give us that
{
\begin{align*}
 (\YR_t)^i =  \esssup_{a \in \mathscr{A}_{i,t}^\Re} \left(U^{\Re,a}_t-A_t^a\right) \leq \esssup_{a \in \mathscr{A}_{i,t}} \left(U^{\Re,a}_t-A_t^a\right) = (\tilde{Y}_t)^i.
\end{align*}}
Then, we finally obtain
{
\begin{align} \label{eq dom 1}
\ckY^\ell \ge Y^\ell \vee (\YR)^\ell \; \text{ for all }\; \ell \in \set{1,\dots,d}\;.
\end{align}}
Furthermore we observe that, for all $\ell \in \set{1,\dots,d}$ and all $t \in  [0,T]$,
{
\begin{align} \label{eq main ineq}
|(Y_t)^\ell - (\YR_t)^\ell | \le |(\ckY_t)^\ell - (Y_t)^\ell | + |(\ckY_t)^\ell - (\YR_t)^\ell |. 
\end{align}}
We will now deal separately with the two terms in the right hand side of the above inequality.
\\

2.a We start by studying the first term. From the representation in terms of switched BSDEs given in \eqref{eq representation ckY Y YR}, we know that $(\ckY_t)^\ell = \ckU^{\cka}_t {-A_t^{\cka}}$ and $(Y_t)^\ell \ge U^{\cka}_t{-A_t^{\cka}}$ with $U^a$ solution to \eqref{eq U pas check}.
Indeed $\cka \in \mathscr{A}_{\ell,t}$ is the optimal strategy associated to the driver $\ckf$ and is \emph{a priori} sub-optimal for the driver $f$.
Combining this with \eqref{eq dom 1}, we obtain that
\begin{align*}
0 \le (\ckY_t)^\ell - (Y_t)^\ell \le \ckU^{\cka}_t - U^{\cka}_t
\end{align*}
and we only need now to control the right hand  inequality. 
By applying Itô's formula to the process $e^{\beta t}|\ckU^{\cka}_t - U^{\cka}_t|^2$ and by using assumption \HYP{f}, usual computations lead to, for some $\beta >0$, 
\begin{align*}
 &e^{\beta t}|\ckU^{\cka}_t - U^{\cka}_t|^2 + \mathbb{E}_t\left[\int_t^T e^{\beta s} |\ckV^{\cka}_s-V^{\cka}_s|^2\ud s\right]\\
 \le& \mathbb{E}_t\left[ \int_t^T e^{\beta s} \left[2C|\ckU^{\cka}_s - U^{\cka}_s|(|\ckV^{\cka}_s-V^{\cka}_s|+|Y_s - \tYR_s|)-\beta |\ckU^{\cka}_s - U^{\cka}_s|^2 \right]\ud s \right]\\
 \le & \mathbb{E}_t\left[ \int_t^T e^{\beta s} \left[(2C^2-\beta)|\ckU^{\cka}_s - U^{\cka}_s|^2+ |\ckV^{\cka}_s-V^{\cka}_s|^2+|Y_s - \tYR_s|^2\right]\ud s \right],
\end{align*}
and then, for any $\beta$ large enough,
\begin{align} \label{eq step 2.a final}
e^{\beta t}|(\ckY_t)^\ell - (Y_t)^\ell |^2 \le  \EFp{t}{\int_t^T e^{\beta_s}|Y_s - \tYR_s|^2 \ud s}.
\end{align}

2.b We now study the second term in the right hand side of \eqref{eq main ineq}. Combining \eqref{eq dom 1} and
the representation in term of ``switched BSDEs'' given by \eqref{eq representation ckY Y YR}, we have, for all $t \in [0,T]$, $\ell \in \set{1, \dots, d}$,
{
\begin{align}
\label{eq temp ineq 2.a}
 0 \le (\ckY_t)^\ell - (\YR_t)^\ell \le \ckU^{\cka}_t-A^{\cka}_t - (\YR_t)^\ell
\end{align}}
for some $\cka \in \mathscr{A}_{t,\ell}$. We now introduce the strategy $a$, standing for the projection of $\cka=(\check{\theta}_k,\check{\alpha}_k)$ on
the grid $\Re$, namely: $a:=(\theta_k,\alpha_k) \in \mathscr{A}^\Re_{t,\ell}$ defined by
\begin{align*} 
 \theta_k = \inf \set{r \ge \check{\theta}_k\;,\; r \in \Re} \; \text{ and }\; \alpha_k = \check{\alpha}_k.
\end{align*}
Note that, if the optimal strategy $\cka$ has many times of switching on $(r_j,r_{j+1}]$, where $r_j$ and $r_{j+1}$ belong to the grid $\Re$, the projected
strategy $a$ will have many instantaneous switches at $r_{j+1}$, see also Remark \ref{re swithing strategy} .
\\
From Corollary \ref{co rep switch cor and dor}, we have {$\YR_t \ge U^{\Re,a}_t-A^a_t$} which, combined to  \eqref{eq temp ineq 2.a}, leads to
{
\begin{align} \label{eq interm 2.b.0}
 |(\ckY_t)^\ell - (\YR_t)^\ell| \le |\ckU^{\cka}_t-A_t^{\cka} - (U^{\Re,a}_t-A_t^a)|\,.
\end{align}}
{We introduce continuous processes $\ckG^{} := \ckU^{\cka^{}}-A^{\cka^{}}$ and $\Gamma = U^{\Re,a} - A^a$.
We then have, for all $s \in [t,T]$,
\begin{align*}
 \ckG_s^{} - \Gamma_s = \ckG_T^{} - \Gamma_T + \int_s^T \set{\ckf^{\cka^{}_s}(u,\ckV^{\cka^{}}_u) - f^{a_s}(X_u,\tYR_u,V^{\Re,a}_u) } \ud u 
 -\int_s^T (\ckV^{\cka^{}}_u - V^{\Re,a}_u)\ud W_u \;.
\end{align*}
}
{By applying Itô's formula to the process $e^{\beta s}|\ckG_s^{} - \Gamma_s|^2$ and by using assumption \HYP{f}, usual computations lead to, for $\beta >0$ large enough, 
\begin{align}
\nonumber &e^{\beta t}|\ckG_t^{} - \Gamma_t|^2 \\
\nonumber \le& \mathbb{E}_t\left[ \int_t^T e^{\beta s} \left[2C|\ckG_s^{} - \Gamma_s|\{|\ckf^{\cka^{}_s}(s,\ckV^{\cka^{}}_s)- \ckf^{a_s}(s,\ckV^{\cka^{}}_s)| +|\ckV^{\cka^{}}_s - V^{\Re,a}_s|+|Y_s - \tYR_s|\}\right]\ud s\right]\\
\nonumber &-\beta \mathbb{E}_t\left[ \int_t^T e^{\beta s} \left[|\ckG_s^{} - \Gamma_s|^2 \right]\ud s \right]- \mathbb{E}_t\left[\int_t^T e^{\beta s} |\ckV^{\cka^{}}_s - V^{\Re,a}_s|^2\ud s\right]+\EFp{t}{e^{\beta T}|\ckG_T^{} - \Gamma_T|^2}\\
 \le&\EFp{t}{e^{\beta T}|\ckG_T^{} - \Gamma_T|^2 
 + Ce^{\beta T}\int_t^T  |\ckf^{\cka^{}_s}(s,\ckV^{\cka^{}}_s)- \ckf^{a_s}(s,\ckV^{\cka^{}}_s)|^2\ud s +  \int_t^T e^{\beta s}|Y_s -\tYR_s |^2\ud s}.
 \label{eq interm 2.b 1}
\end{align}}

{
On one hand, using \HYP{z} we compute that
\begin{align}
\nonumber
\int_t^T |\ckf^{\cka^{}_s}(s,\ckV^{\cka^{}}_s)- \ckf^{a_s}(s,\ckV^{\cka^{}}_s)|^2 \ud s
 &= \int_t^T \left| \sum_{k=1}^{N^{\cka^{}}} \ckf^{\check{\alpha}_{k-1}}(s,\ckV^{\cka^{}}_s)
     (\1_{\set{\check{\theta}^{}_{k-1} \le s < \check{\theta}^{}_{k} }}
      - \1_{\set{ {\theta}_{k-1} \le s < {\theta}_{k} }} ) \right|^2 \ud s
 \\
 & \le C |N^{\cka}|^2 \sup_{s \in [0,T]}(1 + |X_s|^2 + |Y_s|^2 + |\tYR_s|^2)|\Re|.
 \label{eq interm 2.b.2}
\end{align}
On the other hand, by using \HYP{f} we obtain
\begin{align}
\label{eq interm 2.b.3}
 |\ckG_T^{} - \Gamma_T|^2 = |A_T^{\cka^{}}-A_T^a|^2 \le C |N^{\cka}|^2 \sup_{1 \le k \le \kappa} \sup_{r \in [r_{k-1},r_k]}|X_r-X_{r_k}|^2.
\end{align}
}

{Combining \eqref{eq interm 2.b.2} and \eqref{eq interm 2.b.3} with \eqref{eq interm 2.b.0} and \eqref{eq interm 2.b 1}, we get
\begin{align} \nonumber
  e^{\beta t} |(\ckY_t)^\ell - (\YR_t)^\ell|^2 &\le e^{\beta t}|\ckG_t^{} - \Gamma_t|^2 \\
\label{eq step 2.b final}
&\le 
  \EFp{t}{C_\beta\cE(\Re) + 2\int_t^T e^{\beta u}|Y_u -\tYR_u |^2\ud u},
\end{align}
with 
$$
\cE(\Re) := |N^{\cka}|^2 \sup_{s \in [0,T]}(1 + |X_s|^2 + |Y_s|^2 + |\tYR_s|^2) |\Re| + |N^{\cka}|^2 \sup_{1 \le k \le \kappa} \sup_{r \in [r_{k-1},r_k]}|X_r-X_{r_k}|^2. 
$$}

2.c {Combining \eqref{eq step 2.a final}  and \eqref{eq step 2.b final}  with \eqref{eq main ineq}, we obtain, for all $t \le s \le T$,
\begin{align*}
\EFp{t}{e^{\beta s} |Y_s - \YR_s|^2} &\le 
  C_\beta\EFp{t}{\cE(\Re)} + 2\int_s^T \EFp{t}{e^{\beta u}|Y_u -\tYR_u |^2}\ud u\\
  &\le 
  C_\beta\EFp{t}{\cE(\Re)} + 2\int_s^T \EFp{t}{e^{\beta u}|Y_u -\YR_u |^2}\ud u
\end{align*}
since $\YR_u = \tYR_u$, $\ud u$ a.e. Then, a direct application of Gronwall's lemma gives us 
\begin{align*}
|Y_t - \YR_t|^2 \le \EFp{t}{e^{\beta t} |Y_t - \YR_t|^2} \le 
 C_\beta \EFp{t}{\cE(\Re)}.
\end{align*}
Using Jensen's inequality, Doob's maximal inequality and Cauchy-Schwarz inequality, the previous inequality allows us to obtain
\begin{align*}
\EFp{}{\sup_{t \in [0,T]} |Y_t - \YR_t|^2} \le& C \EFp{}{\cE(\Re)^2}^{1/2}\\
 \le& C\EFp{}{|N^{\cka}|^8}^{1/4}\EFp{}{\sup_{s \in [0,T]}(1 + |X_s|^8 + |Y_s|^8 + |\tYR_s|^8)}^{1/4} |\Re|\\
&+ C\EFp{}{|N^{\cka}|^8}^{1/4} \EFp{}{\sup_{1 \le k \le \kappa} \sup_{r \in [r_{k-1},r_k]}|X_r-X_{r_k}|^8}^{1/4}.
\end{align*}
Finally, we just have to apply estimates of Proposition \ref{prop estimee Y Z K continu}, Proposition \ref{prop estimee Y Z K}, classical estimate for $X$, and Theorem 1 in \cite{Fischer-Nappo-09} to get
\begin{align*}
 \EFp{}{\sup_{t \in [0,T]} |Y_t - \YR_t|^2} \le C|\Re|+C|\Re|\log (2T/|\Re|)
\end{align*}
and
\begin{align*}
 \EFp{}{\sup_{t \in [0,T]} |Y_t - \tYR_t|^2} = \EFp{}{\sup_{t \in [0,T]} |Y_{t^+} - \YR_{t^+}|^2} \le C|\Re|+C|\Re|\log (2T/|\Re|).
\end{align*}}
To conclude the proof, we just have to remark that the term $\sup_{1 \le k \le \kappa} \sup_{r \in [r_{k-1},r_k]}|X_r-X_{r_k}|^2$ does not appear in $\cE(\Re)$ when cost functions are constant.
\eproof

\begin{Proposition} \label{pr theorem 5.3}
 Let us assume that \HYP{z} and \HYP{f} are in force, then the following holds:
 \begin{equation*}
   \mathbb{E}\left[\int_0^T \abs{Z_s-Z_s^{\Re}}^2\ud s\right]\leq C \abs{\Re}^{1/2}\sqrt{\log(2T/|\Re|)}.
 \end{equation*}
If cost functions are constant, the previous inequality holds true without the term $\sqrt{\log(2T/|\Re|)}$.
\end{Proposition}
\proof
Introduce $\delta \widetilde{Y} := Y - \widetilde{Y}^{\Re}$, $\delta Y : = Y - Y^{\Re}$, $\delta Z : = Z - Z^{\Re}$ and $\delta f:=f(X,Y,Z)-f(X,\widetilde{Y}^{\Re},Z^\Re)$. Applying It\^o's formula to the \emph{c\`adl\`ag} process $|\delta \widetilde{Y}|^2$, we get
\begin{align*}
 |\delta \widetilde{Y}_0|^2 + \int_0^T |\delta Z_s|^2\ud s =  |\delta \widetilde{Y}_T|^2- 2\int_{0}^T \delta \widetilde{Y}_{s^-} \ud \delta \widetilde{Y}_s - \sum_{0 < s \le T} |\delta \widetilde{Y}_s -\delta {Y}_s|^2.
\end{align*}
Recalling that $\delta \widetilde{Y}_{s^-}=\delta Y_s$, $\int_0^T \delta Y_s \ud K_s^{\Re} \ge 0$ and the Lipschitz property of $f$, standard arguments lead to
\begin{align*}
 \EFp{}{|\delta \widetilde{Y}_0|^2} +\EFp{}{\int_0^T |\delta Z_s|^2\ud s} \le C\EFp{}{\int_0^T \delta Y_s \ud K_s } \le C\EFp{}{\sup_{0 \le  t \le T} |\delta Y_t|^2}^{1/2} \EFp{}{K_T^2}^{1/2}.  
\end{align*}
Then, using Proposition \ref{prop estimee Y Z K continu} and Proposition \ref{pr theorem 5.2} concludes the proof.
\eproof

As a by-product we get a strong estimate on $Z$.
\begin{Corollary} \label{co control Z}
Let us assume that assumption \HYP{f} is in force. Then we have
\begin{align*}
|Z_t| \le \bar{L}(1 + |X_t|) \quad \ud \P\otimes \ud t \; \textrm{a.e.}
\end{align*}
where $\bar{L}$ is the constant that appears in \eqref{eq co borne ZR}.
\end{Corollary}

\proof Let us introduce a new generator $\hat{f}(x,y,z):=f(x,y,\rho_x(z))$ with $\rho_x$ the projection 
on the Euclidean ball of radius $\bar{L}(1+\abs{x})$ where $\bar{L}$ comes from 
the estimate on $\ZR$ given in \eqref{eq co borne ZR}. We easily have that $\hat{f}$ is a Lipschitz function such that
$$\abs{\hat{f}(x,y,z)} \leq C(1+\abs{x}+\abs{y}).$$
We denote $(\hat{Y},\hat{Z},\hat{K})$ the solution of the obliquely reflected BSDE with generator $\hat{f}$.
Since \HYP{f} is in force, we can use \eqref{eq co borne ZR} for the discretely reflected BSDE
with generator $\hat{f}$ and we get that 
\begin{align*}
 |\hat{Z}^\Re_t| \le \bar{L}(1 +|X_t|) \quad \ud \P\otimes \ud t \; \textrm{a.e.}
\end{align*}
Using Proposition \ref{pr theorem 5.3}, we take $|\Re| \rightarrow 0$ and we obtain that 
$$ |\hat{Z}_t| \le \bar{L}(1 +|X_t|)  \quad \ud \P\otimes \ud t \; \textrm{a.e.}$$ 
and then 
$$\hat{f}(t,X_t,\hat{Y}_t,\hat{Z}_t) = {f}(t,X_t,\hat{Y}_t,\hat{Z}_t) \quad \ud \P\otimes \ud t \; \textrm{a.e.}$$
Thus, by uniqueness of the solution to the obliquely reflected BSDE, we have that $\hat{Z}=Z$, concluding the
proof of the Corollary.
\eproof

\begin{Theorem}
\label{th convergence continuement ref vers discret refl}
We assume that \HYP{f} is in force. Then we have
 \begin{align}
 \esp{\sup_{t \in [0,T]}|Y_t-Y^\Re_t|^2 + 
 \sup_{t \in [0,T]}|Y_t-\tilde{Y}^\Re_t|^2 } \le C |\Re|\log (2T/|\Re|),
\end{align}
 and
  \begin{equation}
  \label{convergence ZR vers Z}
   \mathbb{E}\left[\int_0^T \abs{Z_s-Z_s^{\Re}}^2\ud s\right]\leq C\sqrt{\abs{\Re}\log (2T/|\Re|)}.
 \end{equation}
 If, furthermore, cost functions are constant, previous estimates hold true without the $\log (2T/|\Re|)$ term.
\end{Theorem}

\proof 
Thanks to  Corollary \ref{co control Z}, we can replace the generator $f$ by  $\hat{f}(x,y,z):=f(x,y,\rho_x(z))$ with $\rho_x$ the projection 
on the Euclidean ball of radius $\bar{L}(1+\abs{x})$ without modifying our BSDEs. Since \HYP{z} is in force for the generator $\hat{f}$, we can apply Proposition \ref{pr theorem 5.2} and Proposition \ref{pr theorem 5.3} and the theorem is proved.
\eproof

\subsection{Proof of Theorem \ref{th convergence schema vers solution continuement reflechie}}
\label{section continuous time subsection convergence}
Combining the previous results with the control of the error between the discrete-time scheme and the discretely obliquely reflected BSDE derived in Section \ref{section etude schema}, we obtain the convergence of the discrete time scheme to the solution of the continuously obliquely reflected BSDE. Namely, we just have to put together Theorem \ref{th convergence continuement ref vers discret refl} and Theorem \ref{th convergence scheme vers discret refl}. {We emphasize the fact that taking $|\Re| \sim |\pi|^{1/2}$ allows us to minimize our error upper bound on $Y$, while the error upper bound on $Z$ does not converge to $0$. On the other hand, taking $|\Re| \sim |\pi|^{1/3}$ allows us to minimize our error upper bound on $Z$ and to obtain a non optimal error upper bound on $Y$.}


 \section*{Acknowledgements} 
 Authors would like to thanks an anonymous referee and the associate editor for their very careful reading and helpful comments that have deeply improved the manuscript. They also would like to thanks Sa\"id Hamad\`ene for having identify an important error in the preliminary version.

 %
 %

 \appendix 
 \section{Appendix}

 \subsection{Proof of Proposition \ref{prop estimee Y Z K}}

 Observing that on each interval $[r_j,r_{j+1})$, $(\tilde{Y}^{\Re},Z^{\Re})$ solves a standard BSDE, existence and uniqueness follow from a concatenation procedure and \cite{parpen90}.

 Concerning estimates, we cannot apply directly Proposition 2.1 in \cite{Chassagneux-Elie-Kharroubi-10} since we have a generator $f$ with a coupling in $y$. Our strategy is to apply Proposition \ref{prop existence unicite EDSRs obliqu refle cont et discr avec generateur aleatoire} with terminal condition $\xi = g(X_T)$, random generator $F(s,z)=f(X_s,\tilde{Y}^\Re_s,z)$ and costs $C^{ij}_s=c^{ij}(X_s)$. So, we just have to show that \HYP{F_p} is in force. Thus, using the fact that $f$ is a Lipschitz function with respect to $y$, it is sufficient to control $\widetilde Y^{\Re}$ in $\mathscr{S}^p$.

 As in the proof of Theorem 2.4 in \cite{hamzha10}, we consider two nonreflected BSDEs bounding $\tY^{\Re}$.
 Define the $\mathbb{R}^d$-valued random variable $\breve{g}(X_T)$ and the random map $\breve{f}$ by $(\breve{g})^j(x):=\sum_{i=1}^d \abs{(g)^i}$ and $(\breve{f})^j(\omega,t,z):=\sum_{i=1}^d \abs{(f)^i(X_t(\omega),\tY^{\Re}_t(\omega),z)}$ for $1 \leq j \leq d$. We then denote by $(\breve{Y},\breve{Z}) \in (\mathscr{S}^{p} \times \mathscr{H}^p)$ the solution of the following nonreflected BSDE:
 \begin{equation*}
  \label{equation 2 prop estimee Y Z K}
  \breve{Y}_t = \breve{g}(X_T) + \int_t^T \breve{f}(s,\breve{Z}_s)\ud s - \int_t^T \breve{Z}_s \ud W_s, \quad 0 \leq t \leq T.
 \end{equation*}
Since all the components of $\breve{Y}$ are similar, $\breve{Y} \in \mathcal{Q}$. We also introduce $(\mathring{Y},\mathring{Z})$ the solution of the BSDE
$$\mathring{Y}_t = g(X_T) + \int_t^T f(X_s,\tY^{\Re}_s,\mathring{Z}_s)\ud s -\int_t^T \mathring{Z}_s \ud W_s, \quad 0 \leq t \leq T.$$
Using a comparison argument on each interval $[r_j,r_{j+1})$ and the monotonicity property of $\mathcal{P}$, we straightforwardly deduce $(\mathring{Y})^i \leq (\tY^{\Re})^{i} \leq (\breve{Y})^i$, for all $1 \leq i \leq d$. Since $(\mathring{Y},\breve{Y})$ are solutions to standard non-reflected BSDEs, classical estimates (see e.g. \cite{Briand-Delyon-Hu-Pardoux-Stoica-03}) lead to
\textcolor{black}{
\begin{align*}
\EFp{t}{\sup_{t \leq s \leq T} | \tY^{\Re}_s|^p} &\le \EFp{t}{\sup_{t \leq s \leq T} | \mathring{Y}_s|^p + \sup_{t \leq s \leq T} | \breve{Y}_s|^p}\\
 & \le C_p\EFp{t}{ \abs{g(X_T)}^p +\int_t^T \abs{f(X_s,\tY^{\Re}_s,0)}^p\ud s}\\
 & \le C_p\EFp{t}{ 1 + \sup_{s \in [t,T]}|X_s|^p+ \int_t^T \sup_{s \leq u \leq T} \abs{\tY^{\Re}_u}^p \ud s}.
\end{align*}
Finally, using Gronwall's lemma we get
\begin{align*}
\EFp{t}{\sup_{t \leq s \leq T} | \tY^{\Re}_s|^p} \le C_p\, \EFp{t}{1 + \sup_{s \in [t,T]}|X_s|^p}
\end{align*}
which leads to, recall \eqref{eq control X},
\begin{align}
\EFp{t}{\sup_{t \leq s \leq T} | \tY^{\Re}_s|^p} \le C_p\, (1+|X_t|^p )\;, \label{eq control YRe markov}
\end{align}
and in particular to
 $
 |\tilde{Y}^{\Re}|_{\mathscr{S}_p}  \le C_p.
$
}

\subsection{A priori estimates}
\label{se app apriori switched rep}
In this section, we prove a generic estimate for a process 
that can be represented by using switched BSDEs. This result is tailor-made for the solution of 
obliquely reflected BSDEs.
For a positive process $\beta \in \mathscr{S}^2$, we denote by  $\bar{\mathscr{A}}$ the set of strategies $a \in \mathscr{A}$, satisfying
\begin{align}\label{eq app estim Na}
\EFp{t}{|\mathcal{N}^a|^2}^\frac12 \le \beta_t\,,\quad  \text{ for } t \le T. 
\end{align}
We consider a process $\fX \in \mathscr{S}^p$, for all $p \ge 2$, and 
for $a \in \bar{\mathscr{A}}$, we define
\begin{align*}
\fA^{a}_t \;:=\; \sum_{j=1}^{\mathcal{N}^a} \gamma^a_{\theta_j} \fX_{\theta_j} \1_{\{\theta_j\le t \le T\}} \;,
\end{align*}
where $\gamma$ is a process in $\mathscr{S}^2$ essentially bounded by a constant $\Lambda$. 
We also consider a process $\fY \in \mathscr{S}^2$ which is given
by $_t = (\fY^{i}_t)_{1 \le i \le d}$  s.t. $\fY^{i}_t = \fU^a_t$ for some $a \in \bar{\mathscr{A}}\cap \mathscr{A}_{i,t}$ where, for $t \le r \le T$,
 \begin{align}
 \label{definition fU}
  \fU^a_r =  \nu^a\fX_T + \int_r^T F^a(s,\fX_s,\fU^a_s,\fV^a_s, \fY_s) \ud s - \int_r^T \fV^a_s \ud W_s + \fA^a_T - \fA^a_r \;
 \end{align}
 with $\nu^a$ a $\cF_T$-measurable random variable essentially bounded by $\Lambda$ and $F$ a progressively measurable map satisfying
 \begin{align} \label{eq ass Fa}
 |F^a(s,x,u,v,y)| \le 
    \Lambda(|x|+|u|+|v|+|y|) \;.
 \end{align}
{Since $F^a$ depends on $\fY$ in \eqref{definition fU}, the definition of $\fY$ is implicite. We emphasize that we are not interested in the existence of the process $\fY \in \mathscr{S}^2$, we just want to obtain some estimates on this process.}

\begin{Proposition} \label{pr estim gene}
\begin{align*}
 |\fY_r|^2 \le C_\Lambda (1+\beta_r)\EFp{r}{\sup_{r \le s \le T}|\fX_s|^4}^\frac12, \quad r \in [0,T].
\end{align*}
\end{Proposition}
\proof Let us introduce
 $
 \fG^a = \fU^a + \fA^a \,.
 $
%
Applying It\^o's formula, we obtain for all $r \le t \le u \le T$,
\begin{align*}
 \EFp{r}{|\fG^a_u|^2 + \int_u^T|\fV^a_s|^2 \ud s}
 \le
 \EFp{r}{|\fG^a_T|^2 + 2 \int_u^T\fG^a_sF^a(s,\fX_s,\fU^a_s,\fV^a_s, \fY_s) \ud s}\;.
\end{align*}
Using classical arguments and the assumption on $F$, we obtain
\begin{align}\label{eq estim gene 1}
 \EFp{r}{|\fG^a_u|^2 }
 \le C_\Lambda
 \EFp{r}{\sup_{t \le s \le T}|\fX_s|^2 + \int_u^T  |\fY_s|^2 \ud s} 
  +\sup_{t \le s \le T}\EFp{r}{|\fA^a_s|^2}\;.
\end{align}
We observe that, for $t \le s \le T$,
\begin{align*}
\EFp{r}{|\fA^{a}_s|^2} &= \EFp{r}{ \left| \sum_{j=1}^{\mathcal{N}^a} \gamma^a_{\theta_j} \fX_{\theta_j} \1_{\{\theta_j\le s \le T\}} 
\right|^2
}\;
\\
&\le \Lambda \EFp{r}{\mathcal{N}^a\sup_{t \le s \le T }|\fX_s|^2} \le \Lambda \beta_r\EFp{r}{\sup_{t \le s \le T}|\fX_s|^4}^\frac12\;.
\end{align*}
Inserting the previous inequality into \eqref{eq estim gene 1}, we obtain,
\begin{align}\label{eq estim gene 2}
 \EFp{r}{|\fG^a_u|^2 }
 \le C_\Lambda
 (1+\beta_r)\EFp{r}{\sup_{r \le s \le T}|\fX_s|^4}^\frac12 + C_\Lambda \EFp{r}{\int_u^T  |\fY_s|^2 \ud s }\;.
\end{align}
In particular, for all $r \le t \le T$, we compute
\begin{align*}
 \EFp{r}{|\fY_t|^2 } = \sum_{i=1}^d \EFp{r}{|\fY^i_t|^2 }
 \le C_\Lambda
 (1+\beta_r) \EFp{r}{\sup_{r \le s \le T}|\fX_s|^4}^\frac12 + C_\Lambda \EFp{r}{\int_t^T  |\fY_s|^2 \ud s } \;.
\end{align*}
Using Gronwall's Lemma, we get
\begin{align*}
 |\fY_r|^2 \le C_\Lambda (1+\beta_r) \EFp{r}{\sup_{r \le s \le T}|\fX_s|^4}^\frac12 \;.
\end{align*}
\eproof

 \subsection{Proof of Proposition \ref{prop regularite ZR}}
 
Before starting the proof, let us state the following estimates on the $\Lambda$-process appearing in the representation \eqref{eq Repres Z}.
\begin{align}\label{eq lambda sup}
\sup_{a \in \mathscr{A}^{\Re}} \| \sup_{t\le s\le T} \Lambda_{t,s}^{a}\|_{\mathscr{L}^p }\le C^p_L\,, \quad 0\leq t \leq T,~~ p\ge 2\,.
 \end{align}
We also compute from the dynamics of $\Lambda$ that
\begin{align}
\sup_{a \in \mathscr{A}^{\Re}} \left( \|  \Lambda^{a}_{t,t} - \Lambda^{a}_{t,u} \|_{\mathscr{L}^p} +  \|\sup_{t \le s \le T } | \Lambda^a_{u,s} - \Lambda^a_{t,s} | \; \|_{\mathscr{L}^p} \right)
 &\le C^p_L \sqrt{t-u}\;,\quad u\le  t\le T  \,, \quad p\ge 2\,. \label{eq lambda control 1}
\end{align}
 
The proof of Proposition \ref{prop regularite ZR} follows from the same arguments as in the proof of Theorem 3.1 in \cite{Chassagneux-Elie-Kharroubi-10}. The novelty comes from the term $DY$ but the  estimates \eqref{eq estim DY}-\eqref{eq estim DY-DY} allow to control it without
any difficulty.
From Remark \ref{re best approx}, it is clear that
\begin{align}\label{eq majo control}
 \esp{\int_0^T | Z^{\Re}_s - \bar Z^{\Re}_{s} |^2  \ud s } & \le   \esp{\int_0^T | Z^{\Re}_s -  Z^{\Re}_{\pi(s)} |^2  \ud s } \;.
\end{align}
For $s\le T$ and $a = (\alpha_{k}, \theta_{k})_{k\ge 0} \in \mathscr{A}^{\Re}_{s,\ell}$ with $\ell \in \cI$, we define $(V^{a}_{s,t})_{s\le t \le T}$ by
\begin{align*}
 V^{a}_{s,t} := \E_{t} \Big[ & \px{}g^{a_T}(X_T) \Lambda^{a}_{s,T} D_sX_T  - \sum_{k=1}^{N^{a}} \px{} c^{ \alpha_{j-1} , \alpha_j } (X_{\theta_{k}})  \Lambda^{a}_{s,\theta_k}D_sX_{ \theta_k}
  \nonumber\\
 &    
  +
  \int_s^T \left ( \px{} f^{a_u}(\Theta^\Re_u) \Lambda_{s,u}^{a} D_sX_u +\py{}f^a (\Theta^\Re_u) \Lambda^{a}_{s,u}D_s\tY^{\Re}_u) \right) \ud u
 \Big] \;.
 \end{align*}
\newcommand{\besta}[1]{\!a^{#1}}

\vspace{2mm} \noindent
We now fix $\ell \in \cI$ and denote by $\besta{u} \in \mathscr{A}^{\Re}_{u,\ell}$, for $u \le T$, the optimal strategy associated to the representation of $(\tY^{\Re}_{u})^{\ell}$, recalling (ii) in Corollary \ref{co rep switch cor and dor}.

Observe that, by definition, we have
\begin{align}\label{eq strategy}
N^{\besta{t}} = N^{\besta{u}} \text{ and }\;  \besta{t} = \besta{u} , \textrm{ for all } \quad r_j \le t \le u < r_{j+1} \quad \text{and}\quad    \;j< \kappa  \;.
\end{align}
 Fix $i<n$,  and deduce from  \eqref{eq Repres Z} and \eqref{eq strategy} that
 \begin{align}\label{eq control Z 1}
 \esp{|(Z^{\Re}_t)^\ell - (Z^{\Re}_{t_i})^\ell|^2} &=  \esp{|V^{\besta{t}}_{t,t} - V^{\besta{{t_i}}}_{t_i,t_i}|^2} \; \le \; 2\left( \esp{|V^{\besta{{t_i}}}_{t,t} - V^{\besta{{t_i}}}_{t_i,t}|^2} + \esp{|V^{ \besta{{t_i}}}_{t_i,t} - V^{\besta{{t_i}}}_{t_i,t_i}|^2}\right)\,,\qquad\qquad
 \end{align}
 for $t\in [{t_i},{t_{i+1}})$. Combining \HYP{r}, \eqref{eq majo DX},  \eqref{Control Regu DsXt}, \eqref{eq lambda sup}, \eqref{eq lambda control 1}  and Cauchy-Schwartz inequality with the definition of $V^a$, we deduce
\begin{align} \label{eq control Z 2 bis}
\esp{|V^{\besta{{t_i}}}_{t,t} - V^{\besta{{t_i}}}_{t_i,t}|^2} \le C_{L} | \pi |^{\frac12}\;, \quad {t_i}\le t \le {t_{i+1}}\,, \quad i\le n\,.
\end{align}

Since $V^{ \besta{{t_i}}}_{t_i, .}$ is a martingale on $[{t_i},{t_{i+1}}]$, we obtain
 \begin{align}\label{eq control Z 3}
  \esp{|V^{\besta{{t_i}}}_{{t_i},t} - V^{\besta{{t_i}}}_{{t_i},{t_i}}|^2}
  & \le
  \esp{|V^{\besta{{t_i}}}_{{t_i},{t_{i+1}}} - V^{\besta{{t_i}}}_{{t_i},{t_i}}|^2} \nonumber\\
  & \le
  \esp{|V^{\besta{{t_i}}}_{{t_{i+1}},{t_{i+1}}}|^2 - |V^{\besta{{t_i}}}_{{t_i},{t_i}}|^2}
  +
  \esp{|V^{\besta{{t_i}}}_{{t_i},{t_{i+1}}}|^2 - |V^{\besta{{t_i}}}_{{t_{i+1}},{t_{i+1}}}|^2} \nonumber \\
  & \le   \esp{|V^{\besta{{t_i}}}_{{t_{i+1}},{t_{i+1}}}|^2 - |V^{\besta{{t_i}}}_{{t_i},{t_i}}|^2}
  + C_{L} |\pi|^{\frac12}\;, \qquad {t_i}\le t\le {t_{i+1}}\,,
\end{align}

where the last inequality follows from \eqref{eq control Z 2 bis}.
Combining \eqref{eq control Z 1}, \eqref{eq control Z 2 bis}, \eqref{eq control Z 3} and summing up over $i$, we obtain
\begin{align*}
 \esp{ \int_{0}^{T} |(Z^{\Re}_t)^\ell - (Z^{\Re}_{ \pi(t)})^\ell|^2 \ud t} \le C_{L}|\pi|^{\frac12} + |\pi| \Big(  \esp{|V^{\besta{r_{\kappa-1}}}_{T,T}|^2 \!-\! |V^{\besta{0}}_{0,0}|^2 } + \sum_{j=1}^{\kappa-1} (|V^{\besta{\rjm}}_{\rj,\rj}|^2 \!-\! |V^{\besta{\rj}}_{\rj,\rj}|^2) \Big).
\end{align*}
Combined with \eqref{eq majo DX} and \eqref{eq lambda sup}, this concludes the proof since $\ell$ is arbitrary.

 {
 \subsection{Estimate on the scheme solution}
 \label{section estimee schema non reflechi}
 In this subsection we prove a moment estimate for the solution of a classical non reflected scheme. 
More precisely we consider the following scheme
 \begin{enumerate}[(i)]
  \item The terminal condition $\cY^{}_n$ is given by a random variable $\xi^{} \in \cL^2(\mathcal{F}_T)$ 
  \item for $0 \le i <n$,
  \begin{equation*}
\begin{cases}
\cY_{i} := \mathbb{E}[\cY^{}_{i+1}\mid\mathcal{F}_{t_i
}] + h_i
F_i^{}( \cZ^{}_i),
\vspace*{2pt}\cr
\cZ^{}_i := \mathbb{E}[\cY^{}_{i+1}H_i \mid\mathcal{F}_{t_i}] ,
\end{cases}
\end{equation*}
 \end{enumerate}
 with $(H_i)_{0 \le i <n}$ some $\mathbb{R}^{1 \times d}$ independent random vectors such that, for all $0 \le i < n$, $H_i$ is $\mathcal{F}_{t_{i+1}}$-measurable, $\mathbb{E}_{t_i}[H_i]=0$ and \eqref{eq structure H 1}-\eqref{eq structure H 2} are fulfilled.
%
%
%
\begin{Proposition}
 \label{prop estimee schema non reflechi}
We assume that 
\begin{enumerate}[(i)]
 \item For all $i \in \set{0,...,n-1}$, $F_i^{} : \Omega \times \cM^{d,d} \rightarrow \mathbb{R}^d$ is a $\mathcal{F}_{t_i} \otimes \mathcal{B}(\cM^{d,d})$-measurable function,
 \item $|F_i(z)| \le C(|F_i(0)|+|z|)$ for all $z \in \cM^{d,d}$,
 \item $\mathbb{E}\left[|\xi|^2+ \sum_{i=0}^{n-1} \abs{F_i(0)}^2 h_i\right] \le C$. 
\end{enumerate}
Then, we have
$$\mathbb{E}\left[\sup_{0 \le i \le n} \abs{\cY_i^{}}^2\right]+ \mathbb{E}\left[\sum_{i=0}^{n-1} h_i\abs{\cZ_i^{}}^2\right]\le C\mathbb{E}\left[|\xi|^2+ \sum_{i=0}^{n-1} \abs{F_i(0)}^2 h_i\right].$$
\end{Proposition}
\proof
Since calculations are quite standard we only give the sketch of the proof. We observe that the backward scheme can be rewritten equivalently for $i \in \llbracket 0,n \rrbracket$ as
\begin{equation}
\label{eq schema aleatoire general non reflechi ecrit sans esperances conditionnelles}
\cY^{}_{i} = \xi+\sum_{k=i}^{n-1} F_k( \cZ^{}_k) h_k - \sum_{k=i}^{n-1} h_k\lambda_k^{-1} \cZ^{}_k H_k^\top -\sum_{k=i}^{n-1} \Delta \cM_k  
\end{equation}
with the convention $\sum_{k=n}^{n-1}...=0$,  where $(\lambda_k)$ are given by \eqref{eq structure H 1} and where, for all $k \in \llbracket 0,n-1 \rrbracket$, $\Delta \cM_k$ is an $\mathcal{F}_{t_{k+1}}$-measurable random vector satisfying 
\begin{align*}
\mathbb{E}_{t_k}[\Delta \cM_k]=0,\; \mathbb{E}_{t_k}[|\Delta \cM_k|^2] < \infty \;\text{ and }\; \mathbb{E}_{t_k}[\Delta \cM_k H_k]=0.
\end{align*}
We rewrite \eqref{eq schema aleatoire general non reflechi ecrit sans esperances conditionnelles} for $\cY$ between $k$ and $k+1$ with $k \in \llbracket 0,n-1\rrbracket$ and we develop $\abs{\cY_{k+1}}^2$ to get
\begin{align}
\nonumber
|\cY_{k+1}|^2 =& |\cY_{k}|^2-2\cY_k\left(F_k( \cZ^{}_k) h_k -h_k\lambda_k^{-1} \cZ^{}_k H_k^\top - \Delta \cM_k \right)\\
& + |F_k( \cZ^{}_k) h_k- h_k\lambda_k^{-1} \cZ^{}_k H_k^\top-\Delta \cM_k|^2.
 \label{eq ito discret carre schema gene non reflechi} 
\end{align}
By taking the expectation and doing same kind of classical computations than in Step 2 of the proof of Proposition \ref{prop estimee Y Z K schema general} we get
\begin{equation}
\label{eq estimee 1}
\sup_{0 \le i \le n} \mathbb{E}\left[ \abs{\cY_i^{}}^2\right]+ \mathbb{E}\left[\sum_{i=0}^{n-1} h_i\abs{\cZ_i^{}}^2\right]+\mathbb{E}\left[\sum_{i=0}^{n-1} \abs{\Delta \cM_i^{}}^2\right]\le C\mathbb{E}\left[|\xi|^2+ \sum_{i=0}^{n-1} \abs{F_i(0)}^2 h_i\right].
\end{equation}
We can get back to equation \eqref{eq ito discret carre schema gene non reflechi} to have
$$|\cY_{k}|^2 \le |\cY_{k+1}|^2 + 2\cY_k F_k( \cZ^{}_k) h_k-2\cY_k h_k\lambda_k^{-1} \cZ^{}_k H_k^\top-2\cY_k \Delta \cM_k$$
and, by summing up to $n$,
$$|\cY_{k}|^2 \le |\xi|^2 +2\sum_{i=k}^{n-1} \cY_i F_i( \cZ^{}_i) h_i+ 2 \left|\sum_{i=k}^{n-1} \cY_i h_i\lambda_i^{-1} \cZ^{}_i H_i^\top\right| + 2 \left|\sum_{i=k}^{n-1}\cY_i \Delta \cM_i\right|.$$
Thus we obtain
\begin{align*}
 \sup_{0 \le k \le n} |\cY_{k}|^2 \le |\xi|^2+2\sum_{i=0}^{n-1} \cY_i F_i( \cZ^{}_i) h_i+ 2\sup_{0 \le k \le n} \left|\sum_{i=k}^{n-1} \cY_i h_i\lambda_i^{-1} \cZ^{}_i H_i^\top\right| + 2\sup_{0 \le k \le n} \left|\sum_{i=k}^{n-1}\cY_i \Delta \cM_i\right|,
\end{align*}
and, by using assumptions on $F$, Burkholder-Davis-Gundy inequality, \eqref{eq structure H 1}-\eqref{eq structure H 2} and Young inequality, we get
\begin{align*}
 \mathbb{E}\left[\sup_{0 \le k \le n} |\cY_{k}|^2\right] \le& \mathbb{E} \left[ |\xi|^2+\sum_{i=0}^{n-1} \abs{F_i(0)}^2 h_i\right] + C \sup_{0 \le k \le n} \mathbb{E}\left[|\cY_{k}|^2\right] + C\mathbb{E}\left[\sum_{i=0}^{n-1} h_i\abs{\cZ_i^{}}^2\right]\\
 &+ \frac{1}{2}\mathbb{E}\left[ \sup_{0 \le k \le n} |\cY_{k}|^2\right] + C \mathbb{E}\left[\sum_{i=0}^{n-1} \abs{\Delta \cM_i^{}}^2\right].
\end{align*}
Finally, we just have to put \eqref{eq estimee 1} into the last estimate to conclude.
}

 \bibliographystyle{siam}

\begin{thebibliography}{10}

\bibitem{boucha08}
{\sc B.~Bouchard and J.-F. Chassagneux}, {\em Discrete-time approximation for
  continuously and discretely reflected {BSDE}s}, Stochastic Processes and
  their Applications, 118 (2008), pp.~2269--2293.

\bibitem{Briand-Delyon-Hu-Pardoux-Stoica-03}
{\sc P.~Briand, B.~Delyon, Y.~Hu, {\'E}.~Pardoux, and L.~Stoica}, {\em {$L\sp
  p$} solutions of backward stochastic differential equations}, Stochastic
  Process. Appl., {\bf 108} (2003), pp.~109--129.

\bibitem{carlud08}
{\sc R.~Carmona and M.~Ludkovski}, {\em Pricing asset scheduling flexibility
  using optimal switching}, Applied Mathematical Finance, 15 (2008),
  pp.~405--447.

\bibitem{Chassagneux-Richou-14}
{\sc J.~Chassagneux and A.~Richou}, {\em Numerical simulation of quadratic
  {BSDE}s}, Ann. Appl. Probab., 26 (2016), pp.~262--304.

\bibitem{cha09}
{\sc J.-F. Chassagneux}, {\em A discrete-time approximation for doubly
  reflected {BSDE}s}, Advances in Applied Probability, 41 (2009), pp.~101--130.

\bibitem{Chassagneux-Crisan-13}
{\sc J.-F. Chassagneux and D.~Crisan}, {\em Runge-kutta schemes for backward
  stochastic differential equations}, The Annals of Applied Probability, 24
  (2014), pp.~679--720.

\bibitem{Chassagneux-Elie-Kharroubi-11}
{\sc J.-F. Chassagneux, R.~Elie, and I.~Kharroubi}, {\em A note on existence
  and uniqueness for solutions of multidimensional reflected {BSDE}s},
  Electron. Commun. Probab., 16 (2011), pp.~120--128.

\bibitem{Chassagneux-Elie-Kharroubi-10}
{\sc J.~F. Chassagneux, R.~Elie, and I.~Kharroubi}, {\em Discrete-time
  approximation of multidimensional {BSDE}s with oblique reflections}, Ann.
  Appl. Probab., 22 (2012), pp.~971--1007.

\bibitem{Clark-87}
{\sc D.~S. Clark}, {\em Short proof of a discrete {G}ronwall inequality},
  Discrete Appl. Math., 16 (1987), pp.~279--281.

\bibitem{djeham09}
{\sc B.~Djehiche, S.~Hamadene, and A.~Popier}, {\em A finite horizon optimal
  multiple switching problem}, SIAM Journal on Control and Optimization, 48
  (2009), pp.~2751--2770.

\bibitem{ElKaroui-Peng-Quenez-97}
{\sc N.~El~Karoui, S.~Peng, and M.~C. Quenez}, {\em Backward stochastic
  differential equations in finance}, Math. Finance, {\bf 7} (1997), pp.~1--71.

\bibitem{elie2010probabilistic}
{\sc R.~Elie and I.~Kharroubi}, {\em Probabilistic representation and
  approximation for coupled systems of variational inequalities}, Statistics \&
  probability letters, 80 (2010), pp.~1388--1396.

\bibitem{Fischer-Nappo-09}
{\sc M.~Fischer and G.~Nappo}, {\em On the moments of the modulus of continuity
  of {I}tô processes}, Stochastic Analysis and Applications, 28 (2009),
  pp.~103--122.

\bibitem{Gegout-Petit-Pardoux-96}
{\sc A.~G{\'e}gout-Petit and {\'E}.~Pardoux}, {\em \'{E}quations
  diff\'erentielles stochastiques r\'etrogrades r\'efl\'echies dans un
  convexe}, Stochastics Stochastics Rep., {\bf 57} (1996), pp.~111--128.

\bibitem{hamjea07}
{\sc S.~Hamad{\`e}ne and M.~Jeanblanc}, {\em On the starting and stopping
  problem: application in reversible investments}, Mathematics of Operations
  Research, 32 (2007), pp.~182--192.

\bibitem{hamadene2013viscosity}
{\sc S.~Hamad{\`e}ne and M.~Morlais}, {\em Viscosity solutions of systems of
  {PDE}s with interconnected obstacles and switching problem}, Applied
  Mathematics \& Optimization, 67 (2013), pp.~163--196.

\bibitem{hamzha10}
{\sc S.~Hamad{\`e}ne and J.~Zhang}, {\em Switching problem and related system
  of reflected backward {SDE}s}, Stochastic Processes and their applications,
  120 (2010), pp.~403--426.

\bibitem{huytan10}
{\sc Y.~Hu and S.~Tang}, {\em Multi-dimensional {BSDE} with oblique reflection
  and optimal switching}, Probability Theory and Related Fields, 147 (2010),
  pp.~89--121.

\bibitem{klopla92}
{\sc P.~E. Kloeden and E.~Platen}, {\em Numerical solution of stochastic
  differential equations}, vol.~23, Springer, 1992.

\bibitem{majzha05}
{\sc J.~Ma and J.~Zhang}, {\em Representations and regularities for solutions
  to {BSDE}s with reflections}, Stochastic processes and their applications,
  115 (2005), pp.~539--569.

\bibitem{nua96}
{\sc D.~Nualart}, {\em The Malliavin Calculus and Related Topics}, Probability
  and Its Applications, Springer, 2006.

\bibitem{parpen90}
{\sc {\'E}.~Pardoux and S.~G. Peng}, {\em Adapted solution of a backward
  stochastic differential equation}, Systems Control Lett., {\bf 14} (1990),
  pp.~55--61.

\bibitem{portou09}
{\sc A.~Porchet, N.~Touzi, and X.~Warin}, {\em Valuation of a power plant under
  production constraints and market incompleteness}, Mathematical Methods of
  Operations research, 70 (2009), pp.~47--75.

\end{thebibliography}

\def\cprime{$'$}

\end{document}